\documentclass[12pt]{article}

\usepackage{graphicx}
\usepackage{amsfonts}
\usepackage{amsmath}
\usepackage{amssymb}
\usepackage{color}
\usepackage{subfigure}
\usepackage{algorithmic,algorithm}

\def\J{{\mathcal{J}}}
\def\C{{\mathcal{C}}}
\def\A{{\mathcal{A}}}

\def\L{{\mathcal{L}}}

\def\RR{{\mathbb{R}}}

\def\bU{{\bf U}}
\def\bP{{\bf P}}
\def\bB{{\bf B}}

\def\bK{{\bf K}}
\def\bV{\boldsymbol{\V}}
\def\tC{\tilde{\C}}

\def\g{{\mathbf g}}
\def\m{{\mathbf m}}
\def\n{{\mathbf n}}
\def\x{{\mathbf x}}
\def\q{{\mathbf q}}

\def\t{{\mathbf t}}

\def\s{{\mathbf s}}

\def\f{{\mathbf f}}
\def\b{{\mathbf b}}
\def\V{{\mathbf V}}
\def\bnabla{\boldsymbol{\nabla}}
\def\bDelta{\boldsymbol{\Delta}}

\def\Dpartial#1#2{ {\partial #1 \over \partial #2} }

\def\Dpartialn#1#2#3{ {\partial^{#3} #1 \over \partial #2^{#3}} }

\def\Bmp#1{ \begin{minipage}{#1} }
\def\Emp{ \end{minipage} }
\def\Bmpc#1{ \begin{minipage}[c]{#1} }
\def\Bmpt#1{ \begin{minipage}[t]{#1} }
\def\Bmpb#1{ \begin{minipage}[b]{#1} }

\def\ubar{{\overline{u}}}

\newcommand{\argmin}{\operatorname{argmin}}

\newcommand{\eps}{\varepsilon}
\newcommand{\To}{\longrightarrow}

\newtheorem{assumption}{Assumptions}
\newtheorem{lemma}{Lemma}

\makeatletter
\newcommand{\rmnum}[1]{\romannumeral #1}
\newcommand{\Rmnum}[1]{\expandafter\@slowromancap\romannumeral #1@}
\makeatother

\title{A Method for Geometry Optimization in a Simple Model of
  Two-Dimensional Heat Transfer}

\author{X.~Peng\footnotemark[2], \and {K.~Niakhai}\footnotemark[2] \and B.~Protas\footnotemark[2]}

\begin{document}

\maketitle
\renewcommand{\thefootnote}{\fnsymbol{footnote}}

\footnotetext[2]{Department of Mathematics and Statistics, McMaster
  University, Hamilton, ON, Canada}

\renewcommand{\thefootnote}{\arabic{footnote}}

\begin{abstract}
  This investigation is motivated by the problem of optimal design of
  cooling elements in modern battery systems. We consider a simple
  model of two-dimensional steady-state heat conduction
  {described by elliptic partial differential equations} and
  involving a one-dimensional cooling element represented by a contour
  {on which interface boundary conditions are specified}.  The
  problem consists in finding an optimal shape of the cooling element
  which will ensure that the solution in a given region is close (in
  the least squares sense) to some prescribed target distribution. We
  formulate this problem as PDE-constrained optimization and the
  locally optimal contour shapes are found using a gradient-based
  descent algorithm in which the Sobolev shape gradients are obtained
  using methods of the shape-differential calculus. {The main
    novelty of this work is an accurate and efficient approach to the
    evaluation of the shape gradients based on a boundary-integral
    formulation which exploits certain analytical properties of the
    solution and does not require grids adapted to the contour. This
    approach is thoroughly validated and} optimization results
  obtained in different test problems {exhibit} nontrivial shapes
  of the computed optimal contours.
\end{abstract}

{\bf{Keywords:}} 
heat transfer, adjoint-based optimization, shape calculus, 
Sobolev gradients, boundary integral equations 

\bigskip
{\bf{AMS subject classifications:}}
80M50, 
35Q93, 
49Q10, 
49Q12, 
65N38  

\pagestyle{myheadings}
\thispagestyle{plain}
\markboth{X.~PENG, K.~NIAKHAI AND B.~PROTAS}{Optimal Geometry in Two-Dimensional Heat Transfer}

\section{Introduction}
\label{sec:intro}

\subsection{Motivation}
\label{sec:motiv}

\begin{figure}
\begin{center}
\mbox{
\subfigure[]{
\includegraphics[width=0.5\textwidth]{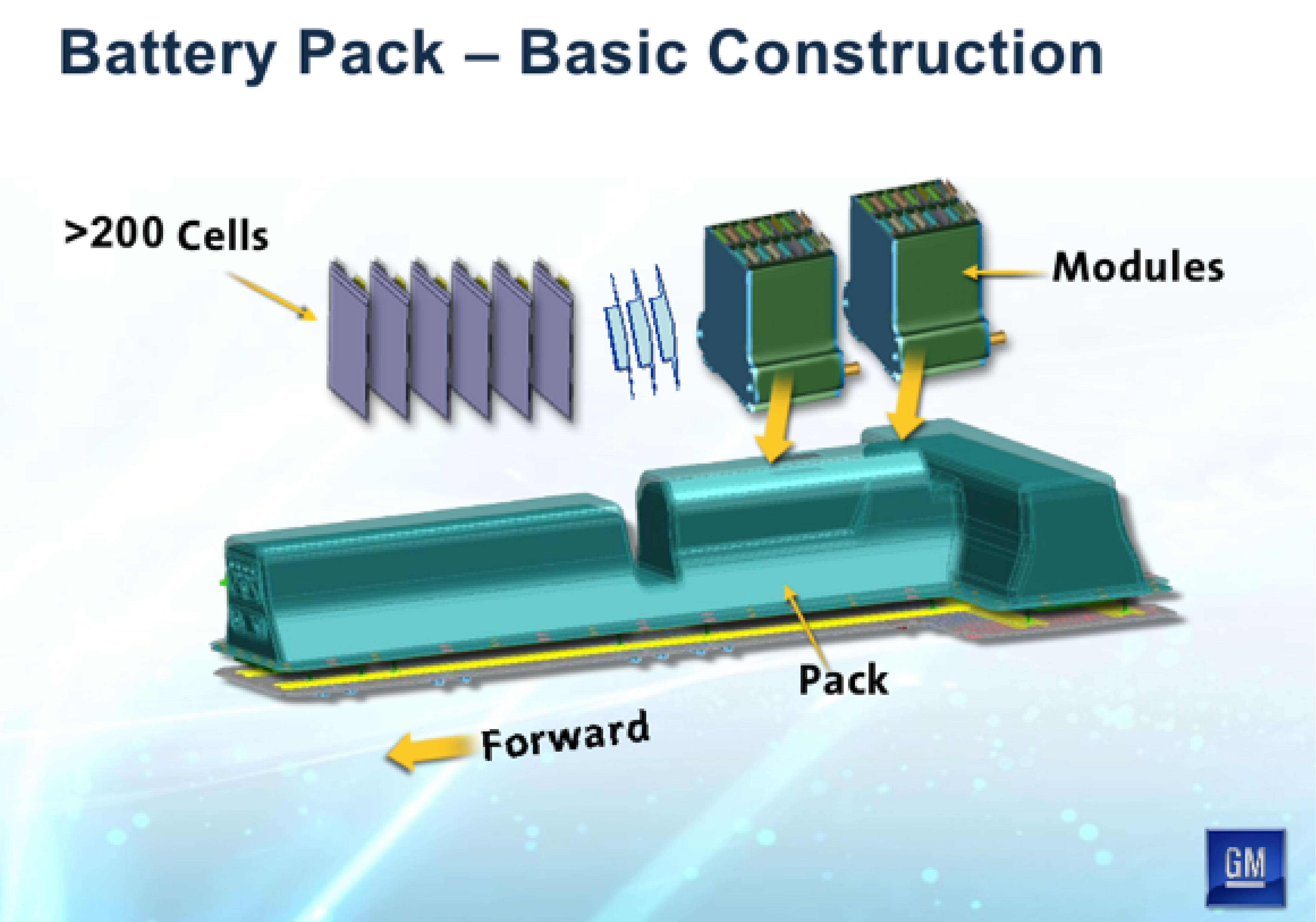}} 
\qquad
\subfigure[]{
\includegraphics[width=0.4\textwidth]{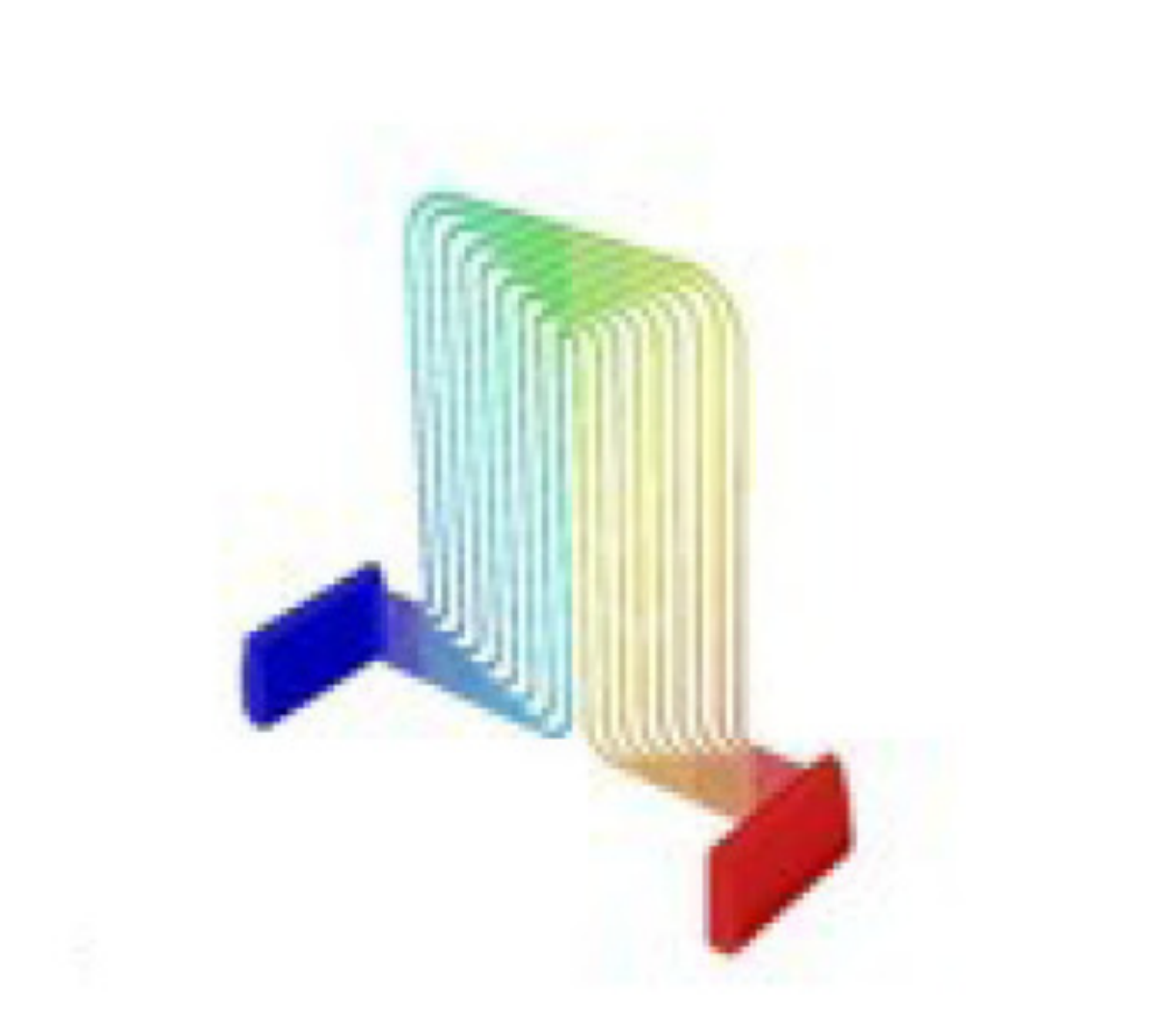}}} 
\caption{(a) Battery system used in hybrid-electric vehicles and (b)
  possible design of the cooling elements (courtesy of General Motors
  of Canada).}
\label{fig:battery_pack}
\end{center}
\end{figure}

The goal of this investigation is to develop and validate a
computational method for optimization of the shape of cooling elements
in general steady heat transfer problems. The motivation for this work
comes from problems encountered in the design of battery systems for
hybrid-electric (HEV) and electric vehicles (EV) \cite{s08} {in
  which a central role is played by} methods of the thermal battery management (TMB)
{ensuring} that the battery operates in a suitable thermal
environment \cite{chkchak10}. A typical
battery system used in automotive applications is shown in Figure
\ref{fig:battery_pack}a, whereas in Figure \ref{fig:battery_pack}b we
present a possible design of the channels with the coolant fluid
acting as the heat-exchange elements. {In these
  applications a key issue is} optimization of the shape of the
cooling elements, so that the temperature distribution {is} as
close as possible to prescribed profiles in some selected regions of
the battery system. Assuming a known distribution of the heat sources
representing the heat generation in the battery, {mathematical
  models of such problems lead to systems of elliptic boundary-value
  problems defined on irregular domains and subject to some
  rather complicated boundary conditions. Optimization of geometry of
  the cooling elements thus leads to shape-optimization problems for
  such systems of equations, and in this study we propose an approach
  based on the continuous (i.e., infinite-dimensional, or
  ``optimize-then-differentiate'', \cite{g03}) formulation and the
  methods of the shape-differential calculus. The main novel
  contribution is the development and validation of an accurate and
  efficient technique based on the boundary-integral formulation for
  the evaluation of the shape gradients which is a key enabler of the
  proposed optimization strategy.}

In the literature {devoted to heat transfer and the related field
  of fluid mechanics} most of the works concerning shape optimization,
or equivalently shape identification, concern problems formulated in
the ``discretize-then-differentiate'' setting, where a
finite-dimensional optimization problem is {set up} based on a
discrete version of the governing equations, somewhat limiting the
flexibility in dealing with different geometries. {Such
  approaches were pursued, for example, in
  \cite{jk11,lchw01,chch03a,chch03b,chch03c,ars02}, and we also
  mention the monograph \cite{mp09}. Approaches based on continuous
  adjoint formulations usually rely on the shape-differential calculus
  to determine the shape sensitivities. The shape calculus, reviewed
  in the monographs \cite{sz92,dz01a,hm03}, is a general suite of
  techniques derived from differential geometry which allow one to
  differentiate solutions of partial differential equations (PDEs) and
  functionals defined on these solutions with respect to variations of
  the domains on which these PDEs are defined. Applications of various
  continuous shape-optimization approaches to problems involving heat
  transfer, fluid flow and phase transformations were investigated in
  \cite{ss10,ws10,yz98,hz07a,vp09,vplg09,bh11,rmp11}.} We add that, as
regards {the numerical representation of free boundaries in PDE
  problems}, there are two main computational approaches, namely,
{the} ``interface capturing'' methods based on the use of
suitable implicit functions, such as the level set formulation
\cite{of02}, and {the} ``interface tracking'' methods which rely
on explicit representations of the boundary. {Since the model
  problem considered here is described by elliptic PDEs, we survey
  below the state-of-the-art numerical techniques used for the
  solution of optimization problems for such system based on the
  continuous formulation.}

\subsection{Review of Computational Methods for Shape Optimization of Elliptic PDEs}
\label{sec:review}

{In both paradigms, i.e., in the approaches relying on the level
  sets to capture the interface and in the methods based on explicit
  interface tracking, optimization problems are typically solved using
  discrete (with respect to some pseudo-time) forms of gradient flows
  in which suitably defined shape gradients are used as the descent
  directions.  Starting with the seminal work \cite{s96}, most
  attention has recently been focused on level-set-based techniques in
  which the level-set function is evolved using the Hamilton-Jacobi
  equation with the velocity field given as an extension of the shape
  gradient away from the interface. Their advantage is that they do
  not require interface-fitted domain discretizations and perform well
  on simple Cartesian grids. The governing and adjoint problems can be
  solved using the immersed interface method \cite{li06}, as was done
  for example in \cite{ikl01,i03}, or with a penalization technique
  \cite{chbgi09}. Regularization aspects of such approaches were
  investigated in \cite{b01}, whereas the study \cite{b03} explored
  formulations resulting from different definitions of the inner
  products for the shape variations. Limitations of such methods arise
  when the boundary conditions and/or the shape gradients defined on
  the interface have a more complicated form (e.g., include
  derivatives), as then they tend to be difficult to evaluate
  accurately on Cartesian grids. On the other hand, shape optimization
  techniques based on explicit interface tracking typically require
  interface-fitted discretization of the domains on which the
  governing and adjoint systems are solved. This discretization then
  needs to be updated during iterations which can be a complicated
  process. Such approaches were reviewed in \cite{dmnv07}, whereas
  some applications to image processing are discussed in
  \cite{sym07,kd12}}

{The approach proposed here is based on explicit interface
  tracking combined with suitably chosen Sobolev gradients. While both
  the shape-differentiation and Sobolev gradients are well-known
  techniques, the main novelty of the proposed approach is a method
  for the evaluation of shape gradients which is based on a
  boundary-integral formulation coupled with an elliptic solver
  constructed using a Cartesian grid. In comparison to the approaches
  described above, it offers the following advantages
\begin{itemize}
\item it is characterized by a high (in principle spectral) accuracy
  in approximating complex interface boundary conditions and
  expressions for the shape gradients, so that only modest resolution
  is required to discretize the contour,
\item as boundary-fitted grids need not be constructed, it can deal
  with fairly complicated contour shapes at a low computational cost.
\end{itemize}
The proposed implementation takes advantage of the analytic structure
of the governing equations. While boundary-integral techniques have
been used for shape optimization of elliptic PDEs, this was typically
done in the discrete setting (i.e., ``discretize-then-differentiate'')
with or without the adjoint equations used to evaluate the shape
sensitivities as in \cite{hch97,hh99,hs06,hl10,nr00,hv02}. In
\cite{afa06} the optimized shape was described in terms of a graph of
a function, so that determination of the gradients did not require
methods of the shape-differential calculus. A boundary-integral
formulation for a time-dependent (parabolic) shape optimization
problem was devised in \cite{ht13}. We add that all of these
approaches relied on the standard techniques of the boundary-element
method (BEM) to evaluate the resulting integral expressions. Finally,
we also mention \cite{rg07} and some references cited therein where
the shape sensitivities were expressed in terms of hypersingular
integral equations (obtained via shape-differentiation of the standard
boundary-integral formulations). We will comment on this interesting
alternative approach at the end of the paper.}

The structure of the paper is as follows: in the next Section we
introduce the mathematical model of the system and state the
optimization problem, in the following Section we {briefly
  describe} a gradient-based descent algorithm based on
shape-differentiation and smoothed (Sobolev) gradients; the
{proposed} computational method {for the solution of} the
governing and adjoint system {and evaluation of the sensitivities
  is presented in detail} in Section \ref{sec:implement}, whereas
validation tests and results demonstrating application of the method
to some selected shape optimization problems are presented in Section
\ref{sec:results}; discussion and conclusions are deferred to Section
\ref{sec:final}.

\section{Mathematical Model and Optimization Problem}
\label{sec:model}

We will consider a simplified model of the problem based on the
following set of assumptions
\begin{assumption} \Bmp{0.1cm}\Emp 
\begin{enumerate}
\renewcommand{\theenumi}{\roman{enumi}}
\item heat transfer is independent of time and occurs via conduction
  only with $k > 0$ representing the constant thermal conductivity,

\item the battery pack is treated as a 2D square region $\Omega
  \subset \RR^2$ with isolated boundary $\partial\Omega$ (i.e., the
  heat flux vanishes on $\partial \Omega$),

\item the distribution of the heat sources in the battery is
  given by the function $q \; : \; \Omega \rightarrow \RR$ which we
  will assume to be square-integrable, i.e., $q \in L_2(\Omega)$; the
  corresponding temperature distribution will be denoted by $u \; : \;
  \Omega \rightarrow \RR$,

\item the cooling element is represented by a $C^1$ curve $\C$ of
  total length $L = \oint_{\C}\, ds$ and characterized by the
  reference temperature $u_0$; the density $w$ of the heat flux
  absorbed by the cooling element {at a point $\x_{\C} \in \C$}
  is modelled using Newton's law of cooling as $w = \gamma (u|_{\C} -
  u_0)$, where $\gamma > 0$ is a constant heat transfer coefficient
  and the temperature field $u$ is continuous across the contour $\C$,
\label{ass_newton}

\item given an arbitrary subdomain $\A \subseteq \Omega$, the target
  temperature distribution is given by $\ubar \; : \; \A \To \RR$.
\end{enumerate}
\end{assumption}
\begin{figure}
\begin{center}
\mbox{
\subfigure[]{
\includegraphics[width=0.35\textwidth]{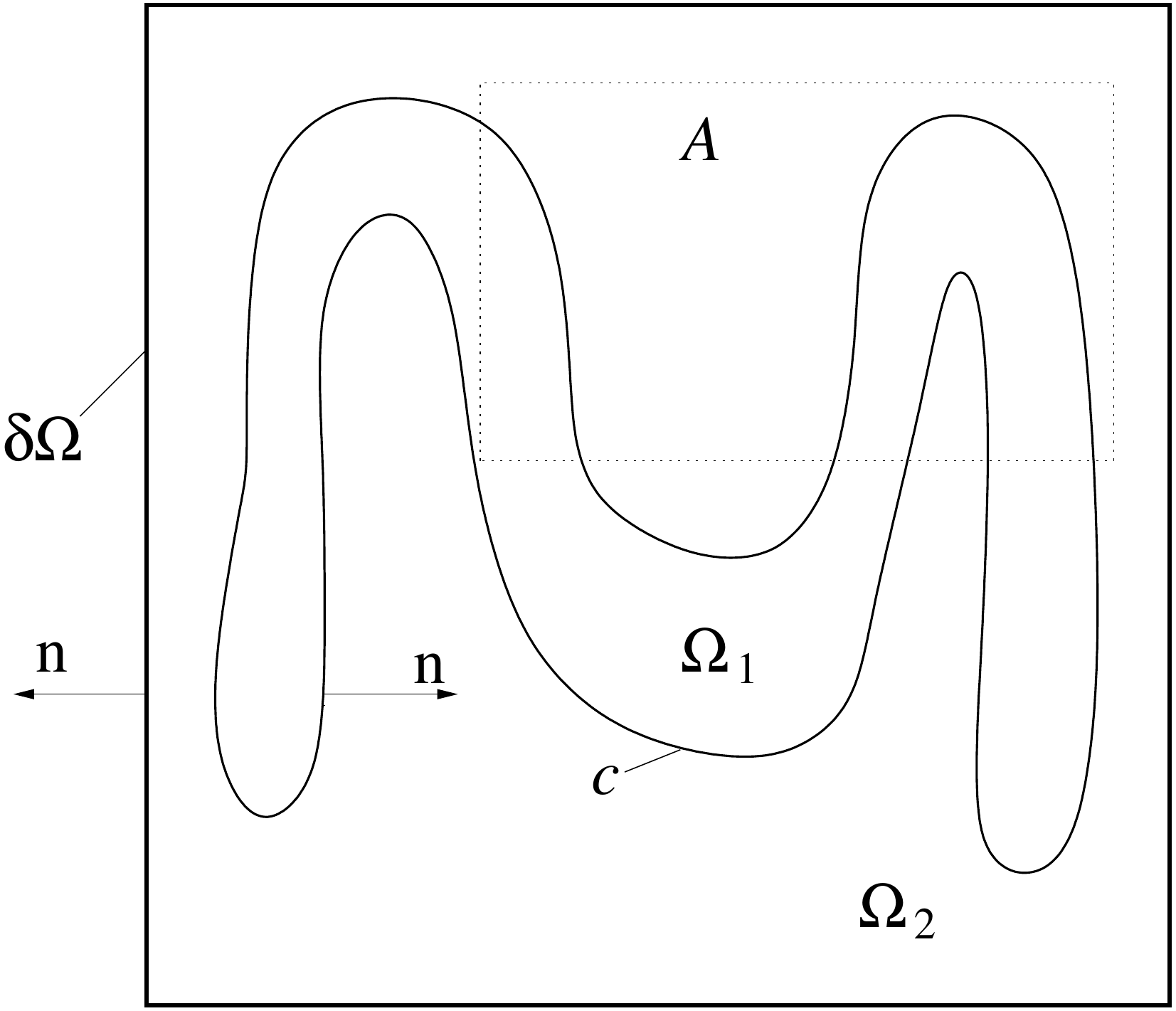}} 
\qquad\qquad
\subfigure[]{
\includegraphics[width=0.35\textwidth]{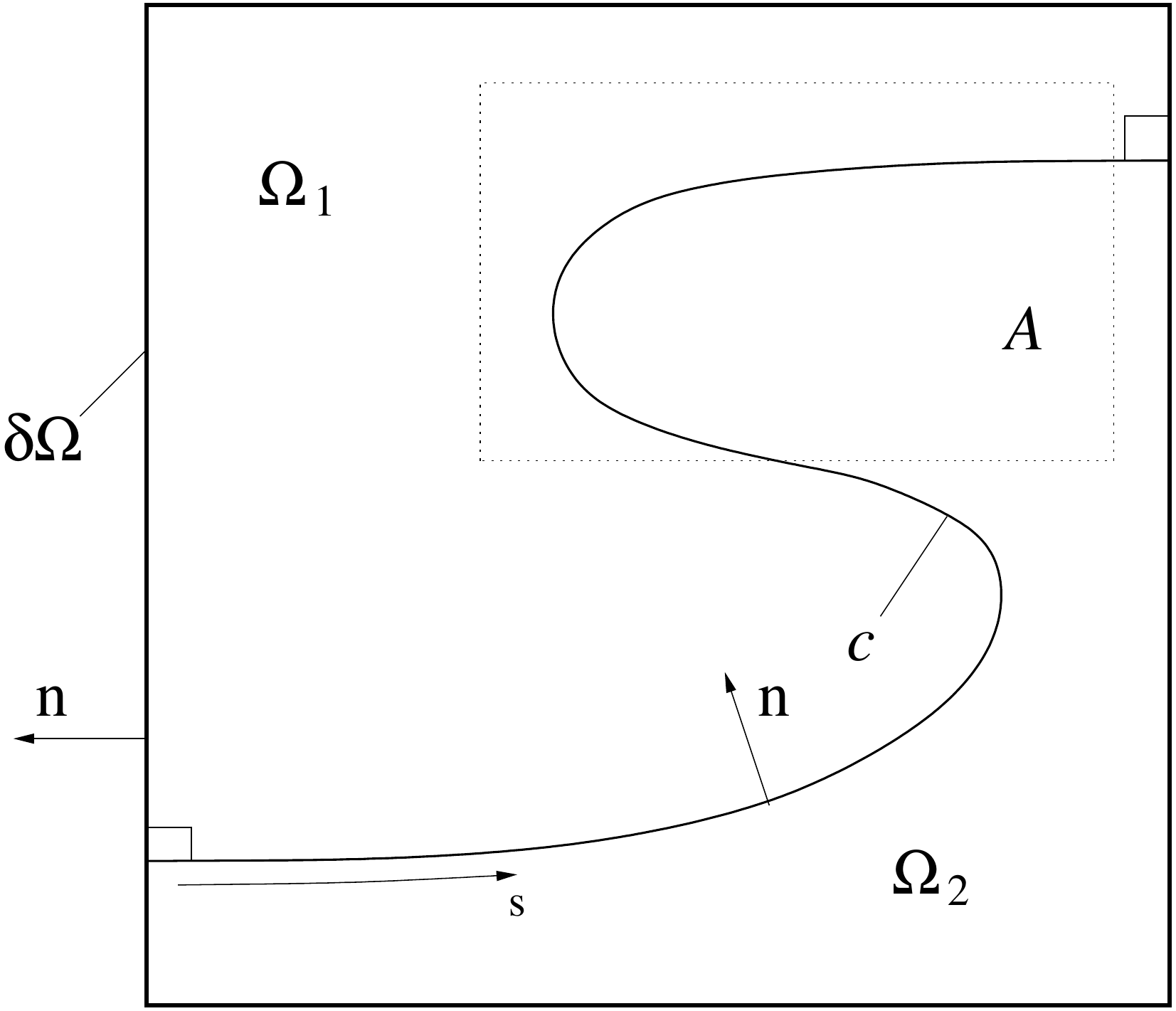}}} 
\caption{Sketch of the domain $\Omega$ with the target region $\A$ for
  {(a) Problem P1 with a closed contour $\C$ and (b) Problem P2 with an
  open contour $\C$ representing the cooling element.}}
\label{fig:domain}
\end{center}
\end{figure}
We restrict our attention to contours which are Lipschitz-continuous
{and will assume that they are parameterized in terms of the
  arc-length coordinate $s \in [0,L]$. Two versions of the problem
  will be considered:
\begin{equation}
\begin{aligned} 
\textrm{P1:} \  \ & \x_{\C}\big|_{s=0} = \x_{\C}\big|_{s=L} \notin \partial\Omega, \\
& u_0 = \textrm{Const}
\end{aligned}\quad
\begin{aligned} 
\textrm{P2:} \  \ & \x_{\C}\big|_{s=0}, \x_{\C}\big|_{s=L} \in \partial\Omega, \ \x_{\C}\big|_{s=0} \neq \x_{\C}\big|_{s=L} \\
& u_0 = u_0(s)
\end{aligned}
\label{eq:P12}
\end{equation}
corresponding, respectively, to a closed contour $\C$ with a constant
reference temperature and to an open contour $\C$ with the reference
temperature $u_0 = u_0(s)$ varying with the arc length.  All
validation tests and a number of optimizations will be performed for
the simper problem P1.  In addition, some optimization problems will
be solved for the more realistic configuration P2 in which $q$ will be
taken to be the distribution of the heat sources in an actual battery
(Figure \ref{fig:battery_pack}b). To fix attention, in Problem P2 we
will assume that $u_0$ increases linearly with the length
corresponding to the coolant liquid heating up as it absorbs heat,
i.e.,
\begin{equation}
u_0(s) = T_a + \frac{T_b - T_a}{L}\, s, \quad s \in [0,L],
\label{eq:u0}
\end{equation}
where $T_a$ and $T_b$ are the prescribed temperatures at the inlet and
outlet.} {Sketches} of the domain $\Omega$ with its different
attributes are shown {for both cases} in Figure \ref{fig:domain}.
{In Problem P2 with an open contour $\C$ the endpoints are
  assumed to attach to the domain boundary at the right angles.}  We
will denote $\Omega_1$ the part of the domain $\Omega$ inside,
{or above}, contour $\C$, and $\Omega_2 \triangleq \Omega
\backslash \overline{\Omega}_1$ its complement, {cf.~Figure
  \ref{fig:domain}} (``$\triangleq$'' means ``equal to by
definition'').  Denoting $u_1 \triangleq u|_{\Omega_1}$ and $u_2
\triangleq u|_{\Omega_2}$ the restrictions of the temperature field to
the {subdomains on the two sides} of contour $\C$, we have the
following mathematical model of the problem
\begin{subequations}
\label{eq:intr_d}
\begin{alignat}{5}
-k\,\Delta u_{1}&=q \quad &&\textrm{in} \ \Omega_{1}, \label{eq:intr_d1}\\
-k\,\Delta u_{2}&=q \quad&&\textrm{in} \ \Omega_{2}, \label{eq:intr_d2}\\
u_{2}&=u_{1} \left(\triangleq u|_{\C} \right) \quad&&\textrm{on} \ \C, \label{eq:intr_d3}\\
k\,\left(\Dpartial{u_{2}}{n}-\Dpartial{u_{1}}{n}\right)&=\gamma\,(u_{1}-u_0) \quad&&\textrm{on} \ \C, \label{eq:intr_d4}\\
k\,\Dpartial{u}{n}&=0 \quad &&\textrm{on}\ \partial\Omega, \label{eq:intr_d5}
\end{alignat}
\end{subequations}
where $\n$ is the unit vector normal to the contour $\C$, or the
boundary $\partial \Omega$, and oriented as shown in Figure
\ref{fig:domain}. {The corresponding unit tangent vector will be
  denoted $\t$.} We add that boundary conditions \eqref{eq:intr_d3}
and \eqref{eq:intr_d4} represent Newton's law of cooling mentioned in
Assumption 1.(\rmnum{4}). {Typically used to model the heat transfer
  in the presence of convection, this law stipulates that the heat
  flux
  $k\left(\Dpartial{u_{2}}{n}-\Dpartial{u_{1}}{n}\right)|_{\x_{\C}}$
  absorbed into the cooling element $\C$ at a given point $\x_{\C} \in
  \C$ is proportional to the difference between the local temperature
  $u|_{\x_{\C}}$ and the reference temperature $u_0$. We remark that
  $u_0$ may be therefore thought of as the temperature of some
  hypothetical coolant liquid circulating in the cooling element,
  although details of this process are neglected in the present model
  (the contour $\C$ has in fact zero thickness).}  Clearly, the
solution $u$ will depend on the shape of the contour, i.e., $u=u(\C)$.
We also remark that, while the differential equations and boundary
conditions in system \eqref{eq:intr_d} are linear in the dependent
variables $u_1$ and $u_2$, problem \eqref{eq:intr_d} is in fact {\em
  geometrically} nonlinear with respect to the shape of the contour
$\C$.  For a discussion of the existence and regularity of solutions
to elliptic boundary-value problems in complicated domains we refer
the reader to monograph \cite{g85}.

The optimization problem, motivated by the industrial applications
discussed in Introduction, is to find an optimal contour $\tilde{\C}$
such that the corresponding solution $\tilde{u} \triangleq
u(\tilde{\C})$ of system \eqref{eq:intr_d} evaluated over $\A$ is as
close as possible to the prescribed target distribution
$\ubar$. Defining the reduced least-squares cost functional as
\begin{equation}
 \J(\C) \triangleq \frac{1}{2}\int_{A}(u-\overline{u})^2\, d\Omega,
\label{eq:J}
\end{equation}
we obtain the following optimization problem
\begin{equation}
\begin{aligned}
& \min_{\C} \J(\C) \\
& \textrm{subject to System} \ \eqref{eq:intr_d}.
\end{aligned}
\label{eq:min1}
\end{equation}
Since in actual applications the length of the contour representing
the cooling element may not be arbitrary, we will also consider a
second optimization problem with the additional constraint on the
contour length, namely,
\begin{equation}
\begin{aligned}
& \min_{\C} \J(\C)  \\
& \textrm{subject to:}\, &&\textrm{System} \ \eqref{eq:intr_d} \\
&\qquad        &&\oint_{\C} \, ds = L_0,
\end{aligned}
\label{eq:min2}
\end{equation}
where $L_0 > 0$ is the prescribed length of the contour $\tC$.
Clearly, problems \eqref{eq:min1} and \eqref{eq:min2} represent
PDE-constrained shape optimization problems. PDE optimization problems
involving shapes of the domains as the control variables require
special treatment \cite{sz92,dz01a,hm03}, and our computational
approach will be based on methods of the shape-differential calculus
recalled in the next Section. Finally, we add that, in principle, in
the statement of optimization problems \eqref{eq:min1} and
\eqref{eq:min2} we should also include the condition $\C \subset
\Omega$ which is equivalent to a suitable set of inequality
constraints.  However, in the interest of simplifying the formulation,
this condition is omitted here, although {as discussed in Section
  \ref{sec:grad} below,} it will be incorporated in the final
computational algorithm.

\section{Gradient-Based Minimization Approach}
\label{sec:grad}

In this Section we review {the} formulation of the optimality
conditions for problems \eqref{eq:min1} and \eqref{eq:min2} and a
gradient-based descent approach for the computational solution of
these problems.  {Since these elements of our approach are rather
  standard, their presentation will be brief.} We consider the
first-order optimality conditions which require the vanishing of a
suitably-defined G\^ateaux (directional) differential evaluated at the
optimal contour $\tC$. We remark that defining such differential and
the related expression for the gradient requires differentiation of
governing system \eqref{eq:intr_d} with respect to the shape of the
domains $\Omega_1$ and $\Omega_2$ on which the PDEs are defined. This
is properly done based on the methods of the {\em shape-differential}
calculus \cite{sz92,dz01a} which rely on a special parametrization of
the domain geometry and provide formulas for shape-differentiation of
general functionals, PDEs and the associated boundary conditions.
Below we briefly {present} this construction and recall the main
results we will need, referring the reader to monographs
\cite{dz01a,hm03} for further details. As a first step, we define the
``velocity'' field $\V \; : \; \Omega \to \RR^2$ which will
parametrize the deformations of the contour $\C$ and of the domains
$\Omega_1$ and $\Omega_2$, so that for every point $\x_\C \in \C$ we
have
\begin{equation}
\x_{\C(\epsilon)} = \x_\C + \epsilon \V,
\label{eq:xt}
\end{equation}
where $0 < \epsilon \ll 1$ is a parameter and $\x_{\C(\epsilon)}$ is
the position of a point on the deformed contour $\C(\epsilon)$.
Relations analogous to \eqref{eq:xt} can also be written for points in
the {deformed} domains {$\Omega_1(\epsilon)$ and
  $\Omega_2(\epsilon)$}. Given a sufficiently regular function $\psi
\; : \; \Omega \to \RR$ and the functionals
{$j_1(\Omega_1(\epsilon)) \triangleq \int_{\Omega_1(\epsilon)}
  \psi(\x;\Omega_1(\epsilon))\, d\Omega$ and $j_2(\C(\epsilon))
  \triangleq \int_{\C(\epsilon)} \psi(\x;\C(\epsilon))\, ds$ defined
  on the perturbed domain and contour}, the corresponding shape
differentials are defined as $j'_1(\Omega_1;\V) \triangleq
\lim_{\epsilon \to 0} \epsilon^{-1} [j_1(\Omega_1(\epsilon)) -
j_1(\Omega_1(0))]$ and $j'_2(\C;\V) \triangleq \lim_{\epsilon \to 0}
\epsilon^{-1} [j_2(\C(\epsilon)) - j_2(\C(0))]$. One of the central
results of the shape-differential calculus is summarized in the
following
\begin{lemma}
The shape differentials of $j_1(\Omega_1)$ and $j_2(\C)$ with
respect to parametrization \eqref{eq:xt} are given by expressions
\begin{subequations}
\label{eq:s}
\begin{align}
j'_1(\Omega_1;\V) &= \int_{\Omega_1}\psi'\,d\Omega + \oint_{\C}\psi\, (\bV\cdot\n)\,ds, \label{eq:s1} \\
j'_2(\C;\V) &= \oint_{\C}\psi'\,d\s + \oint_{\C}\left(\Dpartial{\psi}{n}+\kappa\,\psi\right)(\bV\cdot\n) \, ds, \label{eq:s2}
\end{align}
\end{subequations}
where $\psi'$ is the shape derivative of the integrand function $\psi$
{defined for $\forall \x \in \Omega$ as $\psi'(\x) \triangleq
  \lim_{\epsilon \rightarrow 0} \epsilon^{-1} \left[
    \psi(\x;\Omega_1(\epsilon)) - \psi(\x;\Omega_1(0))\right]$} and
$\kappa$ denotes the signed curvature of the contour $\C$.
\label{thm:lshape}
\end{lemma} 
A detailed proof of Lemma \ref{thm:lshape} can be found, for example,
in \cite{dz01a}. {We remark that, in general, when differentiating
  with respect to the shape of {\em open} contours, expressions
  \eqref{eq:s1} and \eqref{eq:s2} will have additional terms
  proportional to $(\V\cdot\t)$ and localized via Dirac delta
  distributions at the contour endpoints \cite{g93,vp09}.  However, in
  our Problem P2, owing to \eqref{eq:P12} and the assumption that
  contour $\C$ meets the domain boundary $\partial\Omega$ at the right
  angle (cf.~Figure \ref{fig:domain}b), these terms vanish
  identically.  Therefore, for both the closed and open contours only
  the normal component $\zeta \triangleq (\V \cdot \n)|_\C$} of the
perturbation velocity field on the contour $\C$ plays a role in
expressions for shape differentials \eqref{eq:s}. The normal
perturbations {$\zeta = \zeta(s)$, considered} as functions of the
arc-length coordinate, must satisfy certain regularity conditions. It
is sufficient for the perturbation $\zeta$ to belong to the Sobolev
space $H^1(0,L)$ of periodic functions with square-integrable
derivatives on $[0,L]$ (precise definition of the corresponding inner
product will be given in \eqref{eq:ipH1} below).  {We add that contour
  parametrization allows us to recast line integrals, such as
  appearing in \eqref{eq:s1}, \eqref{eq:s2} and below, as definite
  integrals.}

The optimality condition for problem \eqref{eq:min1} is given by
\begin{equation}
\forall_{\zeta \in H^1(0,L)} \quad \J'(\tC;\zeta\n) 
= \int_{\A} (u-\ubar) \, u' \, d\Omega = 0,
\label{eq:opt1}
\end{equation} 
where we note that the subdomain $\A$ is fixed and does not depend on
the perturbation $\zeta$, and $u' = u'(\tC,\zeta\n)$ is the shape
derivative of the solution of governing problem \eqref{eq:intr_d}
evaluated for the optimal contour shape $\tC$. The sensitivity
(perturbation) equation satisfied by $u'$ is obtained by considering a
suitable weak form of system \eqref{eq:intr_d} and
shape-differentiating the resulting integrals using formulas
\eqref{eq:s}, see \cite{hm03},
\begin{subequations}
\label{eq:ptb}
 \begin{alignat}{2}
  k\, \Delta u'_{1}&=0\quad &&\text{in} \ \Omega_{1}, \label{eq:ptb1} \\
  k\, \Delta u'_{2}&=0\quad &&\text{in} \ \Omega_{2}, \label{eq:ptb2} \\
  u'_{2}-u'_{1}&=\left(\Dpartial{u_1}{n}-\Dpartial{u_2}{n}\right)\zeta \quad&&\text{on} \ \C, \label{eq:ptb3} \\
  k\,\left(\Dpartial{u'_{2}}{n}-\Dpartial{u'_{1}}{n}\right)-\gamma\,u'_{1}&=\gamma\,\left[\Dpartial{u_1}{n} + \kappa\,(u_{1}-u_0)\right]\zeta {-\gamma\, u'_0} \quad&&\text{on} \ \C, \label{eq:ptb4} \\
  k\Dpartial{u'}{n}&=0\quad &&\text{on} \ \partial\Omega, \label{eq:ptb5} 
 \end{alignat}
\end{subequations}
where $u'_1 \triangleq u'|_{\Omega_1}$ and $u'_2 \triangleq
u'|_{\Omega_2}$, {and $u'_0$ is the shape-derivative of
  \eqref{eq:u0}
\begin{equation}
u'_0 = u'_0(s;\zeta\n) = \frac{T_b - T_a}{L} \int_0^L \left[ H(s-s') - s / L\right]\kappa\zeta\, ds'
\label{eq:du0}
\end{equation}
in which $H(\cdot)$ is the Heaviside function and whose structure is a
consequence of the dependence of the arc length $s = s(\C)$ on the
contour shape. As a matter of course, in problem P1, $u'_0 \equiv 0$,
cf~\eqref{eq:P12}.}

As regards the second optimization problem \eqref{eq:min2}, we will
incorporate the additional constraint on the length of the contour
$\C$ by defining an augmented cost functional
\begin{equation}
\J_{\alpha}(\C) \triangleq \J(\C) + 
\frac{\alpha}{2} \left(\int_{\C}\,ds - L_0\right)^2,
\label{eq:Ja}
\end{equation}
where $\alpha>0$ is a numerical parameter. The optimality condition
for this second optimization problem is thus
\begin{equation}
\forall_{\zeta \in H^1(0,L)} \quad \J'_{\alpha}(\tC;\zeta\n) 
= \int_{\A} (u-\ubar) u' \, d\Omega + 
\alpha\left(\int_{\tC}\,ds - L_0\right) \oint_{\tC} \kappa\zeta \, ds = 0,
\label{eq:opt2}
\end{equation} 
where we used relationship \eqref{eq:s2} to differentiate the second
term {in \eqref{eq:Ja}. We note that, although it arises from
  rather different mathematical principles, the more systematic
  formulation of the constrained problem using Lagrange multipliers
  would result in an optimality condition quite similar to
  \eqref{eq:opt2}. More precisely, the only difference is that the
  factor $\alpha\left(\int_{\tC}\,ds - L_0\right)$ in \eqref{eq:opt2}
  would be replaced by the Lagrange multiplier $\lambda$. As a result,
  the modification of the descent direction would have the same form
  (but with a different magnitude) in the two cases. On the other
  hand, given the geometric nonlinearity of the constraint
  $\int_{\tC}\,ds = L_0$, the Lagrange multiplier $\lambda$ can be
  rather hard to compute accurately, so for simplicity in this study
  we chose formulation \eqref{eq:Ja}--\eqref{eq:opt2}.}  We emphasize
that optimality conditions \eqref{eq:opt1} and \eqref{eq:opt2} only
characterize {\em local}, rather than global, minimizers and due to
the non-convexity of cost functional \eqref{eq:J}, resulting from the
geometric nonlinearity of system \eqref{eq:intr_d} {and the
  length constraint}, existence of multiple local minima can be
expected.  The (locally) optimal shape $\tC$ can be found
computationally as $\x_{\tC} = \lim_{n \to \infty} \x_{\C^{(n)}}$
using the following gradient-descent algorithm
\begin{equation}
\begin{aligned}
\x_{\C^{(n+1)}} & = \x_{\C^{(n)}} -
\tau_n \bnabla \J\left( \C^{(n)} \right), \qquad n=1,2,\dots, \\
\x_{\C^{(0)}} & = \x_{\C_0},
\end{aligned}
\label{eq:iter}
\end{equation}
where {the points $\x_{\C_0}$ represent the contour $\C_0$ used
  as the ``initial guess'' and} $\tau_n$ is the length of the step
along the descent direction at the $n$-th iteration computed by
solving a line-minimization problem
\begin{equation}
\tau_n=\argmin_{\tau>0}\{\J(\C^{(n)}- \tau \, \bnabla \J( \C^{(n)}) \}.
\label{eq:lmin}
\end{equation}
{There are many different approaches to solving problems of this
  type and in our study we use Brent's iterative method combining the
  golden section search with inverse parabolic interpolation in the
  neighbourhood of the minimum. This approach does not require any
  derivatives with respect to $\tau$ and an efficient implementation
  is discussed in \cite{pftv86}. If $\tau_n$ found by solving problem
  \eqref{eq:lmin} results in the deformed contour $\C^{(n+1)}$
  intersecting the domain boundary $\partial\Omega$, the value of
  $\tau_n$ is suitably reduced to ensure the condition $\C^{(n+1)}
  \subset \Omega$ is always satisfied.}  We add that, while for the
sake of brevity of notation formula \eqref{eq:iter} represents the
steepest descent approach, more advanced optimization methods such as
the Polak-Ribi\'ere version of the nonlinear conjugate gradients
method \cite{nw00} {were used to obtain the results presented in
  Section \ref{sec:optim}. At least for the problems we investigated,
  these approaches were found to systematically outperform the
  steepest descent method.} Clearly, a critical element of
minimization algorithm \eqref{eq:iter} is evaluation at every
iteration of the cost functional gradient $\bnabla \J( \C^{(n)})$. The
Riesz representation theorem \cite{b77} guarantees that it can be
extracted from the G\^ateaux shape differential according to the
formula
\begin{equation}
\J'(\C;\zeta\n) = \Big\langle \bnabla^{H^1} \J(\C), \zeta \Big\rangle_{H^1(0,L)},
\label{eq:rieszH1}
\end{equation}
where 
\begin{equation}
\big\langle z_1, z_2 \big\rangle_{H^1(0,L)}
= \int_0^L z_1 z_2 + \ell^2 \Dpartial{z_1}{s}\Dpartial{z_2}{s}\, ds,
\quad \forall_{z_1, z_2 \in {H^1(0,L)}}
\label{eq:ipH1}
\end{equation}
denotes an inner product in the Sobolev space $H^1(0,L)$ in which $\ell
\in \RR$ is a parameter (which will be shown below to have the meaning of
a length scale). We observe that expressions for the G\^ateaux
differentials $\J'(\C;\zeta\n)$ and $\J'_{\alpha}(\C;\zeta\n)$
appearing in \eqref{eq:opt1} and \eqref{eq:opt2} are not yet in the
form consistent with \eqref{eq:rieszH1}, because the perturbation
$\zeta$ rather than appear as a factor is hidden in boundary
conditions \eqref{eq:ptb3}--\eqref{eq:ptb4} of the sensitivity system
defining $u'$. In order to transform the differential
$\J'(\C;\zeta\n)$ to Riesz form \eqref{eq:rieszH1} we will employ the
{\em adjoint variable} $u^* \; : \; \Omega \to \RR$ which is the
solution of the following {\em adjoint system}
\begin{subequations}
\label{eq:adj}
\begin{alignat}{2}
  -k\, \Delta u^*_{1}&=(u - \overline{u})\chi_{A_1}\quad &&\text{in} \ \Omega_{1}, \label{eq:adj1} \\
  -k\, \Delta u^*_{2}&=(u - \overline{u})\chi_{A_2}\quad &&\text{in} \ \Omega_{2}, \label{eq:adj2} \\
  u^*_{2}-u^*_{1}&=0 \quad&&\text{on} \ \C, \label{eq:adj3} \\
  k\,\left(\Dpartial{u^*_{2}}{n}-\Dpartial{u^*_{1}}{n}\right)&=\gamma\,u^*_{1}\quad&&\text{on} \ \C,  \label{eq:adj4} \\
  \Dpartial{u^*}{n}&=0\quad &&\text{on} \ \partial\Omega, \label{eq:adj5}
\end{alignat}
\end{subequations}
where $u^*_1 \triangleq u^*|_{\Omega_1}$ and $u^*_2 \triangleq
u^*|_{\Omega_2}$, whereas $\chi_{\A_i}$ is the characteristic function
of the region $\A_i \triangleq \Omega_i \bigcap \A$,
$i=1,2$. Following the standard procedure, see e.g.~\cite{g03}, we
obtain
\begin{equation}
\J'(\C;\zeta\n) = \int_0^L \bnabla^{L_2} \J \zeta \, ds, 
\label{eq:dJ}
\end{equation}
where 
\begin{align}
\bnabla^{L_2}\J = & -\gamma\,(u_{1} - u_0) 
\left(\kappa\,u^*_{1} + \Dpartial{u^*_{1}}{n}\right)
 - \gamma\,u^*_{1}\Dpartial{u_{2}}{n} \label{eq:gradJL2} \\
& {-\gamma\, \kappa\, \frac{T_b - T_a}{L} \int_0^L \left[ H(s'-s) - s' / L\right]u^*(s')\, ds'}\qquad \text{on} \ \C.
\nonumber
\end{align}
{The last term in \eqref{eq:gradJL2} stems from the arc length
  dependence of the reference temperature $u_0$ in Problem P2,
  cf.~\eqref{eq:u0}, and vanishes identically in Problem P1.}
Derivation details are presented in \cite{p11,xp11}; in \cite{p11} we
also discuss a symbolic algebra algorithm for automated determination
of the adjoint boundary conditions in PDE optimization problems
characterized by complicated interface conditions such as the problem
considered here. While this is not the gradient we use in the actual
computations, for simplicity in \eqref{eq:dJ}--\eqref{eq:gradJL2} the
gradient $\bnabla^{L_2}\J$ was obtained as the Riesz representer in
the space $L_2(0,L)$ of square-integrable functions. We also add that
the part of G\^ateaux differential \eqref{eq:opt2} associated with the
length constraint is already in the Riesz form, so that the $L_2$
gradient of cost functional $\J_{\alpha}(\C)$ is
\begin{equation}
\bnabla^{L_2}\J_{\alpha} = \bnabla^{L_2}\J + 
\alpha\,\left(\int_{\C}\,ds - L_0\right)\kappa\,\quad \text{on} \ \C. 
\label{eq:gradJaL2}
\end{equation}
{The gradients actually used in minimization algorithm
  \eqref{eq:iter}, namely} the Sobolev gradients $\bnabla^{H^1} \J$
and $\bnabla^{H^1} \J_{\alpha}$, can be obtained from
\eqref{eq:gradJL2} and \eqref{eq:gradJaL2} {by} identifying
\eqref{eq:rieszH1}--\eqref{eq:ipH1} with \eqref{eq:dJ}, and noting the
arbitrariness of the shape perturbations $\zeta \in H^1(0,L)$.
{Then,} after integrating by parts and using the boundary
conditions, we arrive at
\begin{equation}
\begin{alignedat}{2}
& \left(1 - \ell^2\,\Dpartialn{}{s}{2}\right)\,\bnabla^{H^1}\J 
= \bnabla^{L_2}\J &\quad &\textrm{on} \ (0,L), \\
& \textrm{Periodic boundary conditions}  && \textrm{(P1)}, \\
& {\Dpartial{}{s}\,\bnabla^{H^1}\J \Big|_{s=0,L} = 0} && \textrm{(P2)}.
\end{alignedat}
\label{eq:precond}
\end{equation}
Thus, the Sobolev gradient $\bnabla^{H^1}\J$ is obtained {by
  first computing the gradient $\bnabla^{L_2}\J$ from
  \eqref{eq:gradJL2} or \eqref{eq:gradJaL2}, and then by} solving
elliptic boundary-value problem \eqref{eq:precond} defined on the
contour $\C$, a step which is known to be equivalent to low-pass
filtering (smoothing) the $L_2$ gradient with $\ell$ acting as the
cut-off length scale \cite{pbh04}. {In Problem P2 the homogeneous
  Neumann boundary conditions ensure that the Sobolev gradient
  $\bnabla^{H^1}\J$ does not change the angle at which the contour
  $\C$ meets the domain boundary $\partial\Omega$ (which therefore
  always remains $\pi/2$).}  For some other applications of Sobolev
gradients to solution of minimization problems involving PDEs we refer
the reader to {monograph \cite{n10}, articles \cite{kd12,dk10}
  and to articles \cite{vp09,vplg09,b03,sym07} for studies concerned
  specifically with shape optimization}.  The different elements
discussed in the present Section combine into Algorithm \ref{alg1}.
\begin{algorithm}
\begin{algorithmic}
\STATE $n \leftarrow 1$
\STATE $\C^{(0)} \leftarrow $ initial guess $\C_0$
\REPEAT
\STATE solve direct problem \eqref{eq:intr_d}
\STATE solve adjoint problem \eqref{eq:adj}
\STATE evaluate \eqref{eq:gradJL2}--\eqref{eq:gradJaL2} and solve \eqref{eq:precond} to determine $\bnabla^{H^1} \J(\C^{(k)})$
\STATE compute the Polak-Ribi\'ere conjugate direction $\g\left[ \bnabla^{H^1} \J(\C^{(k)}) \right]$
\STATE perform line minimization $\min_{\tau>0}\{\J(\x_\C^{(n)}- \tau \, \g\left[\bnabla \J( \C^{(n)})\right] \}$ 
       to find the step size $\tau_n$, ensuring that $\left(\x_\C^{(n)}- \tau_n \, \g\left[\bnabla \J( \C^{(n)})\right]\right) \notin \partial \Omega$
\STATE obtain $\C^{(n+1)}$ by deforming $\C^{(n)}$ along the conjugate direction
$\g\left[ \bnabla^{H^1} \J(\C^{(n)}) \right]$ with the step size $\tau_n$,
\STATE $n \leftarrow n+1$
\UNTIL $|\,\tau_n|< \varepsilon_{\tau}$ or $|\J(\C^{(n+1)})-\J(\C^{(n)})|<\,
\varepsilon_{\J}|\J(\C^{(n)})|$
\end{algorithmic}
\caption{Iterative minimization algorithm for finding optimal contour shapes $\tC$.
\textbf{Input:} $\varepsilon_{\J}$ and $\varepsilon_{\tau}$ (adjustable tolerances),
 $\C_0$ (initial contour shape) \newline
\textbf{Output:} $\tC$ (optimal contour shape) }
\label{alg1}
\end{algorithm}

\noindent An elegant and accurate numerical solution technique for the
direct and adjoint systems \eqref{eq:intr_d} and \eqref{eq:adj}
{and evaluation of gradient expression \eqref{eq:gradJL2}} is
described in the next Section.

\section{Numerical Implementation}
\label{sec:implement}

In this Section we {present in detail a novel numerical approach
  we devised} to solve the governing and adjoint systems
\eqref{eq:intr_d} and \eqref{eq:adj} at every iteration of Algorithm
\ref{alg1}.  Since these systems have in fact essentially identical
structure, we will focus our discussion on the solution of the first
one. {The methods to tackle Problems P1 and P2 are based on the
  same concept, but differ in regard to some technical details, and to
  fix attention, below we describe the approach applicable to Problem
  P1. Modifications required to solve Problem P2 are summarized
  further below with all details available in \cite{n13}.}  We observe
that both systems \eqref{eq:intr_d} and \eqref{eq:adj} can be regarded
as combinations of two Poisson problems (defined in $\Omega_1$ and in
$\Omega_2$) which are coupled via some complicated (mixed) boundary
conditions on the contour $\C$ separating the two domains. It should
be emphasized that this contour can have an arbitrary, though
non-intersecting, shape. {Given the linearity (with respect to $u_1$ and $u_2$) of equations
  \eqref{eq:intr_d1}--\eqref{eq:intr_d2}, we split problem
  \eqref{eq:intr_d} into two subproblems: a potential problem
  associated with the complex interface boundary condition
  \eqref{eq:intr_d4} and another elliptic problem arising from the
  presence of the source term $q$, which are then coupled using a
  suitable interpolation scheme. Since the solution methods for these
  subproblems are adapted to their analytic structure, we achieve for
  each of them the highest possible (spectral) numerical accuracy.
  While similar techniques have already been used for the solution of
  certain direct problems \cite{ll12}, to the best of our knowledge,
  this direction has not been explored in applications to optimization
  or inverse problems.}

\begin{figure}
\begin{center}
\includegraphics[width=0.5\textwidth]{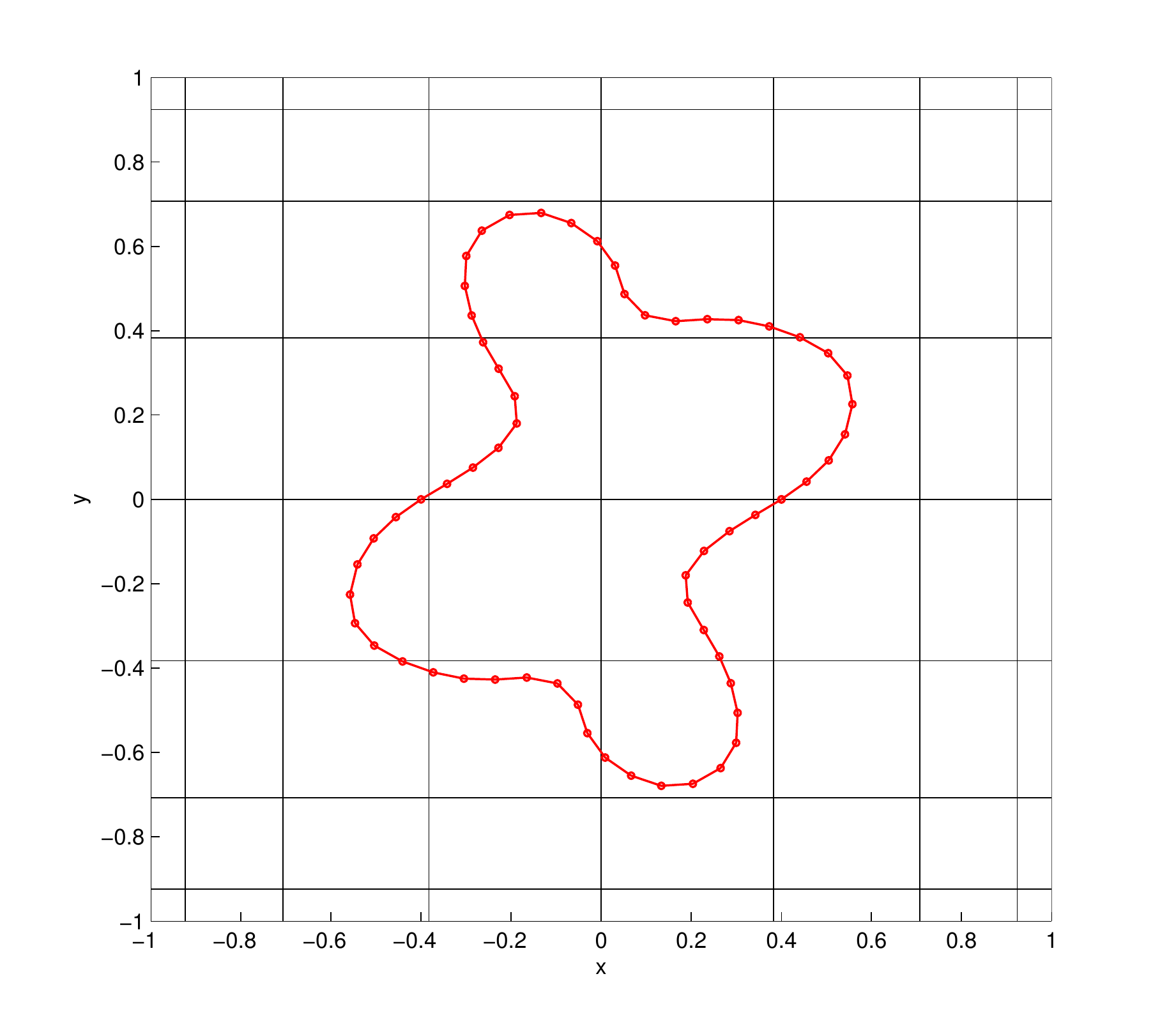}
\caption{Discretization of the domain $\Omega$ and contour $\C$. For clarity, 
the discretization shown is much coarser than used in the actual
computations reported in Section \ref{sec:results}.}
\label{fig:grids}
\end{center}
\end{figure}

As a starting point, we consider the following ansatz for the solution
$u$ of problem \eqref{eq:intr_d}
\begin{equation}
u = u_p + u_h \quad \textrm{in} \ \Omega,
\label{eq:u}
\end{equation}
where the fields $u_p$ and $u_h$ satisfy the following system of PDEs
and boundary conditions equivalent to \eqref{eq:intr_d}
\begin{subequations}
\label{eq:num_d}
 \begin{alignat}{2}
-k\,\Delta u_{p}&=q \qquad \quad && \textrm{in} \ \Omega, \label{eq:num_d_1}\\
\Delta u_{h}&=0 \qquad \quad &&\textrm{in} \ \Omega\setminus\C, \label{eq:num_d_2}\\
u_{h}\big|_1&=u_{h}\big|_2 \qquad \quad &&\textrm{on} \ \C, \label{eq:num_d_3} \\
k\,\left(\Dpartial{u_{h}}{n}\bigg|_2-\Dpartial{u_{h}}{n}\bigg|_1\right)&=\gamma\,(u_p + u_h - u_0) \quad &&\textrm{on} \ \C, \label{eq:num_d_4}\\
\Dpartial{u_{p}}{n}&=-\Dpartial{u_{h}}{n} \qquad \quad &&\textrm{on} \ \partial\Omega. \label{eq:num_d_5}
\end{alignat}
\end{subequations}
We note that the fields $u_p$ and $u_h$ are coupled {only} through 
boundary conditions \eqref{eq:num_d_4} and \eqref{eq:num_d_5}. Since
the field $u_h$ is harmonic in $\Omega \backslash \C$, it admits a
representation in terms of the single-layer potential density $\mu \;
: \; \C \to \RR$
\begin{equation}
\forall_{\x \in \Omega \backslash \C} \quad 
u_{h}(\x)=-\frac{1}{2\pi}\oint_{\C}
\ln\,\left|\x-\x_{\C}\right|\mu(\x_{\C})\,ds.
\label{eq:uh}
\end{equation}
Taking the limit $\x \to \x_\C$ in \eqref{eq:uh}, using boundary
conditions \eqref{eq:num_d_3} and \eqref{eq:num_d_4}, and taking into
account the limiting properties of integrals of type \eqref{eq:uh}
known from the potential theory \cite{h95,k99}, we arrive at a
singular boundary integral equation of Fredholm type II satisfied by
the density $\mu$. Thus, system \eqref{eq:num_d} can be equivalently
rewritten as
\begin{subequations}
\label{eq:upmu}
\begin{alignat}{2}
- k\,\Delta u_p &=q \qquad \quad && \textrm{in} \ \Omega, \label{eq:upmu_1}\\
- \mu(\x)+\frac{\gamma}{2\pi\,k}\oint_{\C}\,\ln\,\left|\x-\x_{\C}\right|\mu(\x_{\C})\,ds &=\frac{\gamma}{k}(u_p-u_0)\, \ && \textrm{on} \ \C, \label{eq:upmu_2}\\
\Dpartial{u_p}{n}&=-\Dpartial{u_h}{n} \qquad \quad && \textrm{on} \ \partial\Omega.\label{eq:upmu_3}
\end{alignat}
\end{subequations}
The new dependent variables are $\{u_p(\x), \, \x \in \Omega; \;
\mu(\x_\C), \x_\C \in \C\}$ and the advantage of this formulation is
that the second variable (potential density) needs to be found on the
contour $\C$ only {and, unlike in original system
  \eqref{eq:intr_d}, there are no differential operators defined on
  the contour $\C$}. For the purpose of discretizing Poisson equation
\eqref{eq:upmu_1} we cover the domain $\Omega$ with a $N\times N$
dyadic Chebyshev grid \cite{t00}, where $N>0$ is the number of grid
points in each direction. Contour $\C$ is represented with $M$ points
equispaced in the arc-length coordinate $s$ ($M$ is taken to be an
even number). These discretizations are shown in Figure
\ref{fig:grids} {(in Problem P2 the discretization of contour
  $\C$ needs to be a bit different, cf.~\cite{n13}).}
We let $u_{p;i,j}^N \triangleq u_p(x_i,y_j)$, $i,j=1,\dots,N$ and
$\mu_l^M \triangleq \mu(s_l)$, $l=1,\dots,M$ denote the discrete nodal
values of the unknowns, where $(x_i,y_j)$ are the coordinates of the
collocation points on the dyadic Chebyshev grid covering $\Omega$,
whereas $s_l$ are the arc-length coordinates of the points
discretizing contour $\C$ {in Problem P1}, i.e., $s_l \triangleq
(l-1) \frac{L}{M}$, $l=1,\dots,M$. We then construct the vectors $\bU$
and $\m$
\begin{subequations}
\label{eq:Umu}
\begin{alignat}{2}
[\bU]_{(i-1)N+j} &= u_{p;i,j}^N, & \qquad & i,j=1,\dots,N, \label{eq:Umu1} \\
[\m]_l &= \mu_l^M, & & l=1,\dots,M,  \label{eq:Umu2}
\end{alignat}
\end{subequations}
and will use the symbol $\bDelta^N$ to denote the discretization of
the Laplace operator $\Delta$ based on the Chebyshev spectral
collocation approach \cite{t00} and corresponding to the Neumann
boundary conditions. Thus, discretization of \eqref{eq:upmu_1} takes
the algebraic form
\begin{equation}
\bDelta^N \bU = \f + \q,
\label{eq:U}
\end{equation}
where $\q \in \RR^{N^2}$ contains the values of the right-hand side
(RHS) function $q$ evaluated at the interior collocation points (and
zeros in the entries corresponding to the boundary nodes) and $\f \in
\RR^{N^2}$ is a vector containing the values of $\Dpartial{u_h}{n}$ at
the boundary nodes, cf.~\eqref{eq:upmu_3}. It can be expressed as
\begin{equation}
\f = \bB \, \m,
\label{eq:f}
\end{equation}
in which $\bB$ is a $N^2 \times M$ matrix operator representing the
discretization via the trapezoidal rule of the relation
\begin{equation}
\Dpartial{u_{h}}{n}\bigg|_{\b_i} = 
- \frac{1}{2\pi}\oint_{S}\frac{(\b_i-\x_{\C})\cdot\,\n}{|\b_i-\x_{\C}|^2}\mu(\x_{\C})\,ds,
\quad i=1,\dots,4N-4,
\label{eq:Dub}
\end{equation}
where $\b_i \in \partial\Omega$ (the rows of $\bB$ corresponding to
the interior grid points are zero). As regards integral equation
\eqref{eq:upmu_2}, we observe that the logarithmic kernel it contains
is in fact singular and, assuming the potential density is a
Lipschitz-continuous function of $s$, the integral is defined as an
improper one. As a standard approach to deal with this issue
\cite{h95,k99}, we rewrite the kernel as
\begin{equation}
\ln|\x_\C(t)-\x_\C(t')| =\frac{1}{2}\ln\left\{\frac{|\x_\C(t)-\x_\C(t')|^2}{4\sin^2\frac{t-t'}{2}}\right\} 
+ \frac{1}{2}\ln\,\left(4\sin^2\frac{t-t'}{2}\right),
\label{eq:kernel}
\end{equation}
where $t,t' \in [0,2\pi]$ are the variables parameterizing contour
$\C$. Therefore, rewriting the line integral in \eqref{eq:upmu_2} as
a definite integral, the boundary integral equation becomes
\begin{equation}
\begin{split}
-\mu(\x_{\C}(t))&+\overset{(\Rmnum{1})}{\overbrace{\frac{\gamma}{2\pi\,k}\int_{0}^{2\pi}
\mu(\x_{\C}(t'))\ln\Bigg|\frac{\x_\C(t)-\x_\C(t')}{2\sin(\frac{t-t'}{2})}\Bigg|\, r(t') \,dt'}}\\
&+\overset{(\Rmnum{2})}{\overbrace{\frac{\gamma}{4\pi\,k}
\int_0^{2\pi}\mu(\x_{\C}(t'))\ln\Bigg[4\sin^2\left(\frac{t-t'}{2}\right)\Bigg]\, r(t') \,dt'}}
=\frac{\gamma}{k}[u_p(\x_\C(t))-u_0], \quad t \in [0,2\pi]
\end{split}
\label{eq:bie}
\end{equation}
where, assuming that the contour parameterization is uniform in the
arc length $s$, we have $r(t) = \big| \frac{d\x(t)}{dt}\big| =
\frac{L}{2\pi}$. We note that integral $(\Rmnum{1})$ has now a regular
kernel (with a removable singularity to be more precise) and can be
evaluated with spectral accuracy in a straightforward manner using the
trapezoidal quadrature. The singularity is now contained in the
improper integral $(\Rmnum{2})$ which can be evaluated analytically as
follows. We approximate the potential density $\mu(t)$ using the
spectrally-accurate trigonometric interpolation \cite{k99}
\begin{equation}
\mu(t) \approx \sum_{j=1}^{M}\mu^M_j\,\L_j(t),
\end{equation}
in which $\L_j(t)$, $j=1,\dots,M$, are the trigonometric cardinal
functions
\begin{equation*}
\L_j(t) \triangleq \frac{1}{M}\sin\left(\frac{M (t-t_j)}{2} \right)\, \cot\left(\frac{t-t_j}{2}\right),
\quad t \in [0,2\pi], \ t\neq\,t_j, \ j=1,\dots,M,
\end{equation*}
where $t_j \triangleq (j-1) \frac{2\pi}{M}$. Defining now
\begin{equation}
R_j^{M}(t) \triangleq -\frac{2}{M}\left\{\sum_{m=1}^{M/2-1}
\frac{1}{m}\cos\left[\,m(t-t_j)\right] + 
\frac{1}{M}\cos\left[\frac{M(t-t_j)}{2}\right]\right\}, \quad j=1,\dots,M,
\label{eq:RjM}
\end{equation}
the improper integral $(\Rmnum{2})$ in \eqref{eq:bie} is approximated as
\begin{equation}
\frac{\gamma}{4\pi\,k}
\int_0^{2\pi}\mu(t')\ln\left[4\sin^2\left(\frac{t-t'}{2}\right)\right]\, r(t') \,dt' 
\approx \frac{\gamma\,L}{4\pi\,k}\sum_{j=1}^{M}\mu^M_jR_j^{M}(t), \ t \in [0,2\pi].
\label{eq:II}
\end{equation}
Therefore, collocating integral equation \eqref{eq:bie} on the grid
points $t_1,\dots,t_M$ yields the following discrete problem
\begin{equation}
\left(\mathbf{I} +\frac{\gamma}{k}\mathbf{K}_1 +\frac{\gamma}{k}\mathbf{K}_2\right)\m 
+ \frac{\gamma}{k}\bP \bU = \frac{\gamma}{k} \, u_0 \mathbf{1},
\label{eq:bie2}
\end{equation}
where $\mathbf{1}$ is a column vector of dimension $M$ with all entries
equal to one and the matrices $\bK_1$ and $\bK_2$ are defined as,
cf.~\eqref{eq:II},
\begin{equation}
[\mathbf{K}_1]_{ij} =  -\frac{L}{2\pi\,M}\ln\Biggl|\frac{\x_{\C}(t_i)-\x_{\C}(t_j)}{2\sin\left(\frac{t_i-t_j}{2}\right)}\Biggr|, \quad
[\mathbf{K}_2]_{ij} = \frac{L}{4\pi}\,R_j^M(t_i), \qquad i,j=1,\dots,M,
\label{eq:KK}
\end{equation}
whereas $\bP$ is an $M \times N^2$ matrix representing interpolation
of the field $u_p$ from the Chebyshev grid onto the points
$\{\x_{\C}(t_1),\dots,\x_{\C}(t_M)\}$ discretizing the contour
$\C$. Thus, combining \eqref{eq:U} and \eqref{eq:bie2}, the final
discrete form of system \eqref{eq:upmu} is
\begin{equation}
\begin{bmatrix}
-\bDelta^{N}\hfil & \mathbf{B}\\
\frac{\gamma}{k}\mathbf{P}\hfil&\mathbf{I} +\frac{\gamma}{k}\mathbf{K}_1 +\frac{\gamma}{k}\mathbf{K}_2
\end{bmatrix}
\begin{bmatrix}
\bU \\
\m
\end{bmatrix}
= \frac{1}{k} \, \begin{bmatrix}
\q  \\
\gamma \, u_0 \, \mathbf{1} 
\end{bmatrix}.
\label{eq:upmu2}
\end{equation}
The accuracy of approximation represented by system \eqref{eq:upmu2}
is ultimately determined by the accuracy of the interpolation operator
$\bP$, and in principle can be spectral, although for reasons of the
numerical stability we have used spline interpolation in the present
study. System \eqref{eq:upmu2} is readily solved using standard
methods of numerical linear algebra, and we refer the reader to thesis
\cite{xp11} for numerical validation and tests of
accuracy. Discretization of adjoint system \eqref{eq:adj} leads to a
discrete problem with the same matrix as in \eqref{eq:upmu2}, but with
a different right-hand side. The $L_2$ gradient $\bnabla^{L_2}\J$ is
obtained from the solution $\left[ (\bU^*)^T \ (\m^*)^T \right]^T$ of
the discrete adjoint problem using relation \eqref{eq:gradJL2}, where
the different terms are computed using boundary conditions
\eqref{eq:adj3}--\eqref{eq:adj4} and the following identities, known
from the potential theory \cite{h95,k99},
\begin{subequations}
\label{eq:jj}
 \begin{alignat}{2}
  \Dpartial{u^*_{h}}{n}\bigg|_1 - \Dpartial{u^*_{h}}{n}\bigg|_2 &= \mu^*,  \label{eq:jj1} \\
   \frac{1}{2}\left[\Dpartial{u^*_{h}}{n}\bigg|_1 + \Dpartial{u^*_{h}}{n}\bigg|_2\right]&=
-\frac{1}{2\pi}\oint_{\C}\frac{\big\langle \n(\x_\C),\x_\C -\x'\big\rangle}{|\x_\C - \x'|^2} \mu^*(\x')\,ds, \label{eq:jj2}
 \end{alignat}
\end{subequations}
valid for all points $\x_\C \in \C$, where $\langle\cdot,\cdot\rangle$
denotes the inner product in $\RR^2$, $u^* = u^*_p + u^*_h$ and
$\mu^*$ is the single-layer potential density associated with $u^*_h$.
We add that the kernel of the integral on the RHS in \eqref{eq:jj2} is
in fact bounded, as we have $\forall_{\x' \in \C}$ $\lim_{\x' \to
  \x_{\C}} \frac{\langle \n(\x_\C),\x_\C -\x'\rangle}{|\x_\C - \x'|^2}
= \frac{\kappa(\x_\C)}{2}$ \cite{h95}.  Accuracy of the cost
functional gradients computed in this way is assessed in the next
Section. {Finally, we remark that after each step of gradient
  algorithm \eqref{eq:iter}, the points $\x_{\C}(t_i)$, $i=1,\dots,M$,
  are no longer distributed uniformly in the arc length $s$. In order
  to retain the spectral accuracy of the solution of equation
  \eqref{eq:bie}, at every iteration we therefore construct, using
  spectral Fourier interpolation \cite{t00}, a new set of collocation
  points $\{\x_{\C}(t_1),\dots,\x_{\C}(t_M)\}$ which are equispaced in
  the arc-length coordinate. The main modification required to adapt
  the method described above to Problem P2 concerns the solution of
  boundary-integral equation \eqref{eq:upmu_2}. Since the integration
  domain is no longer periodic, identity \eqref{eq:kernel} must be
  replaced with a different one and contour $\C$ must be discretized
  using a different set of points \cite{n13}. Moreover, the
  trapezoidal quadratures need to be replaced with the
  Clenshaw-Curtis quadratures whereas the trigonometric interpolation
  with a suitable polynomial technique.}

\section{Computational Results}
\label{sec:results}

In this Section we first perform tests to {thoroughly} validate
{the computational algorithm introduced in Section
  \ref{sec:implement} for the evaluation of cost functional gradients
  $\bnabla\J$. We do this here for Problem P1 and refer the reader to
  \cite{n13} for the corresponding validation tests for Problem P2.
  Next, we apply this method in the framework of} Algorithm \ref{alg1}
to perform shape optimization in a number of test {cases
  concerning Problems P1 and P2.} Throughout this Section we take the
domain to be $\Omega = [-1,1] \times [-1,1]$.

\subsection{Validation of Gradients}
\label{sec:validate}

\begin{figure}
\centering
\mbox{
\subfigure[]{
\includegraphics[scale=0.25]{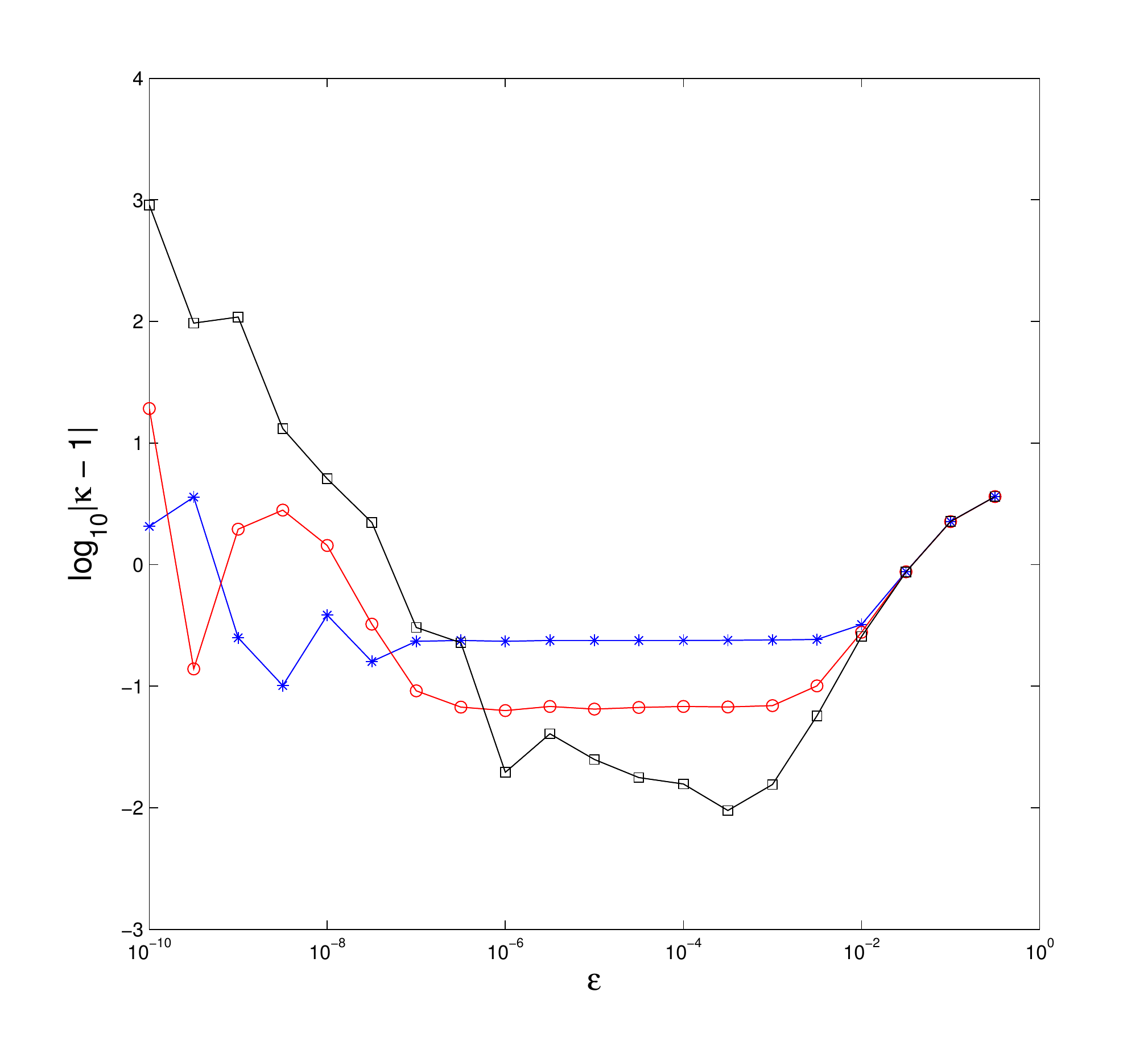}}
\subfigure[]{
\includegraphics[scale=0.25]{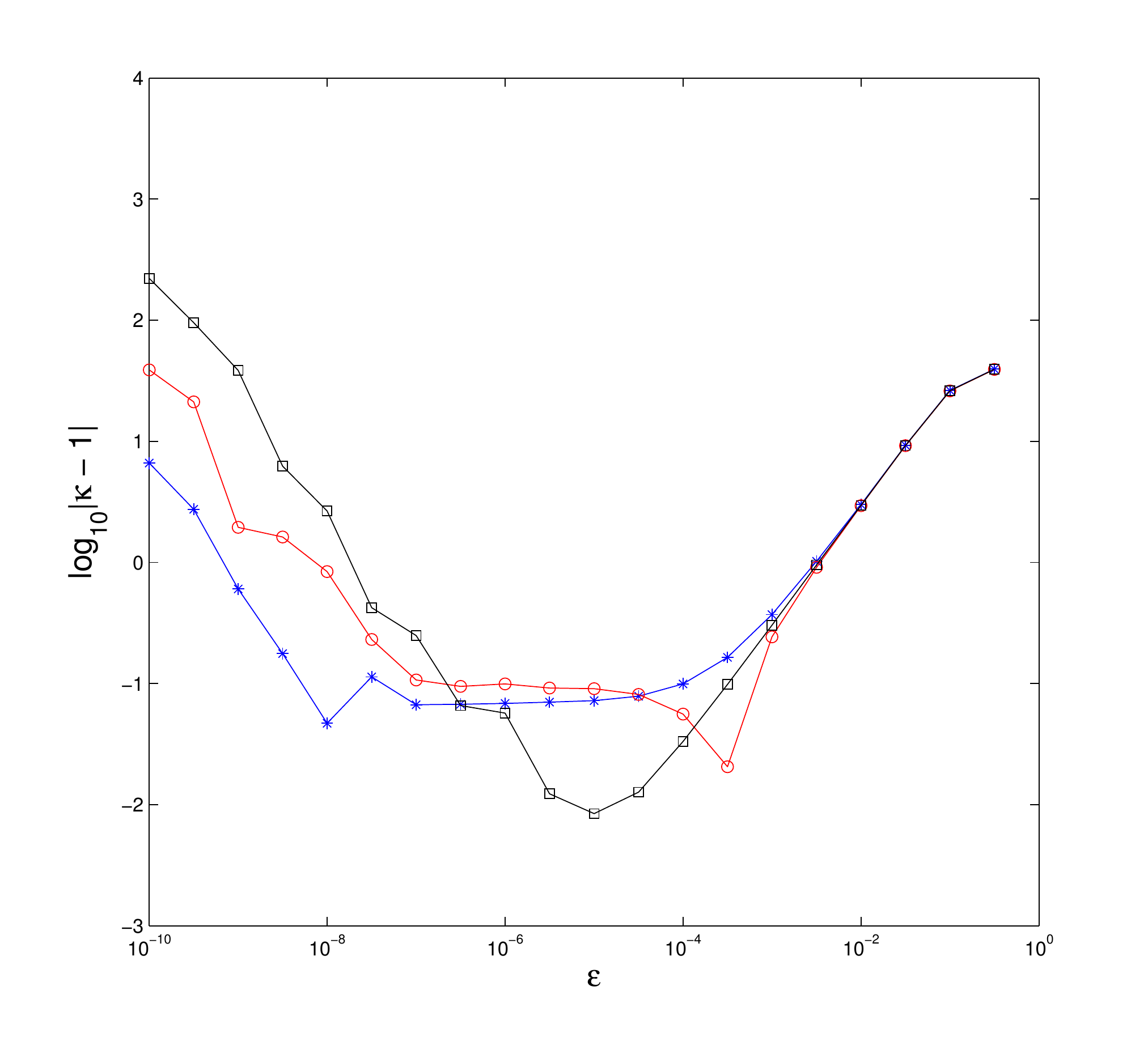}}}
\mbox{
\subfigure[]{
\includegraphics[scale=0.25]{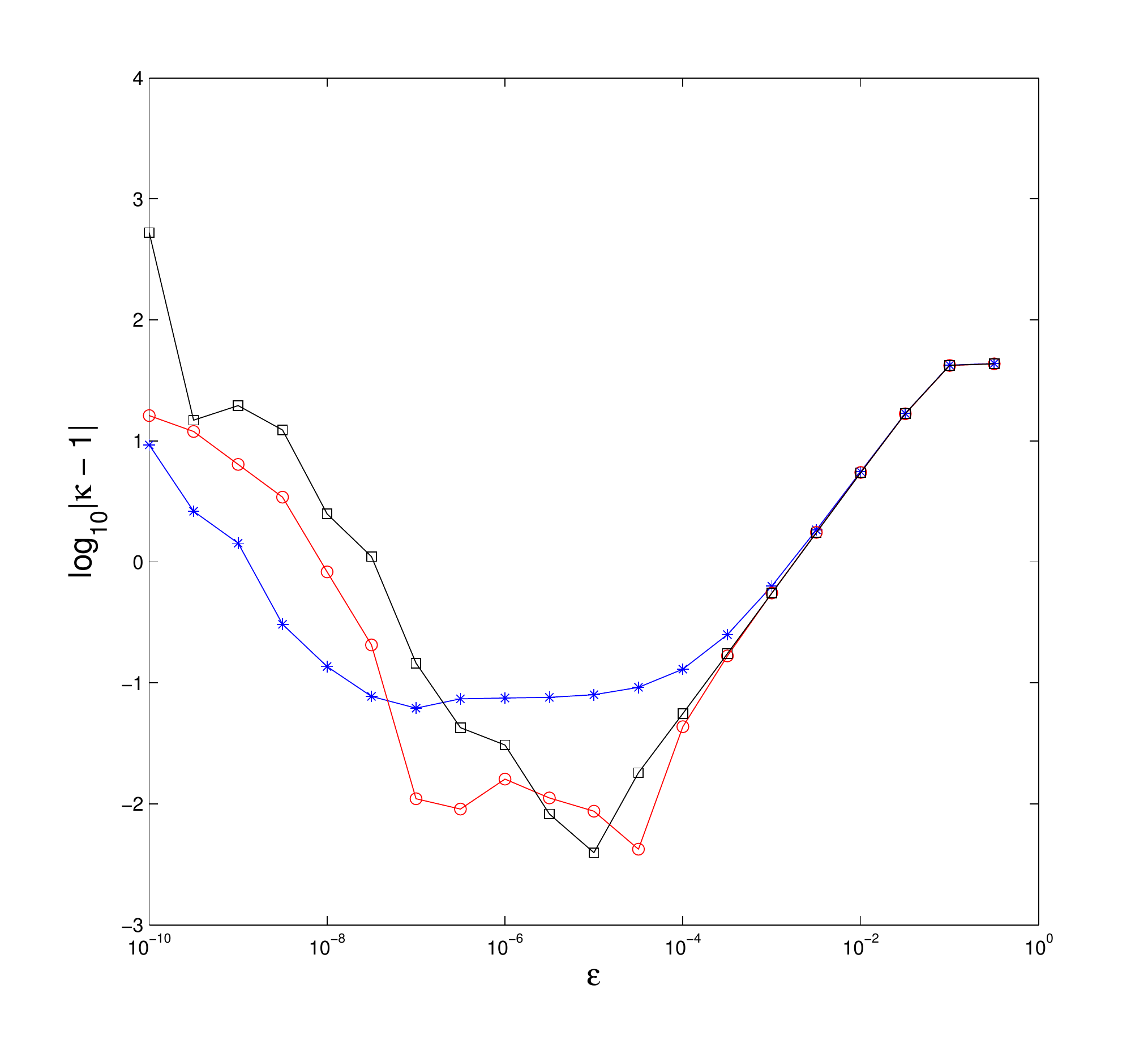}}
\subfigure[]{
\includegraphics[scale=0.2425]{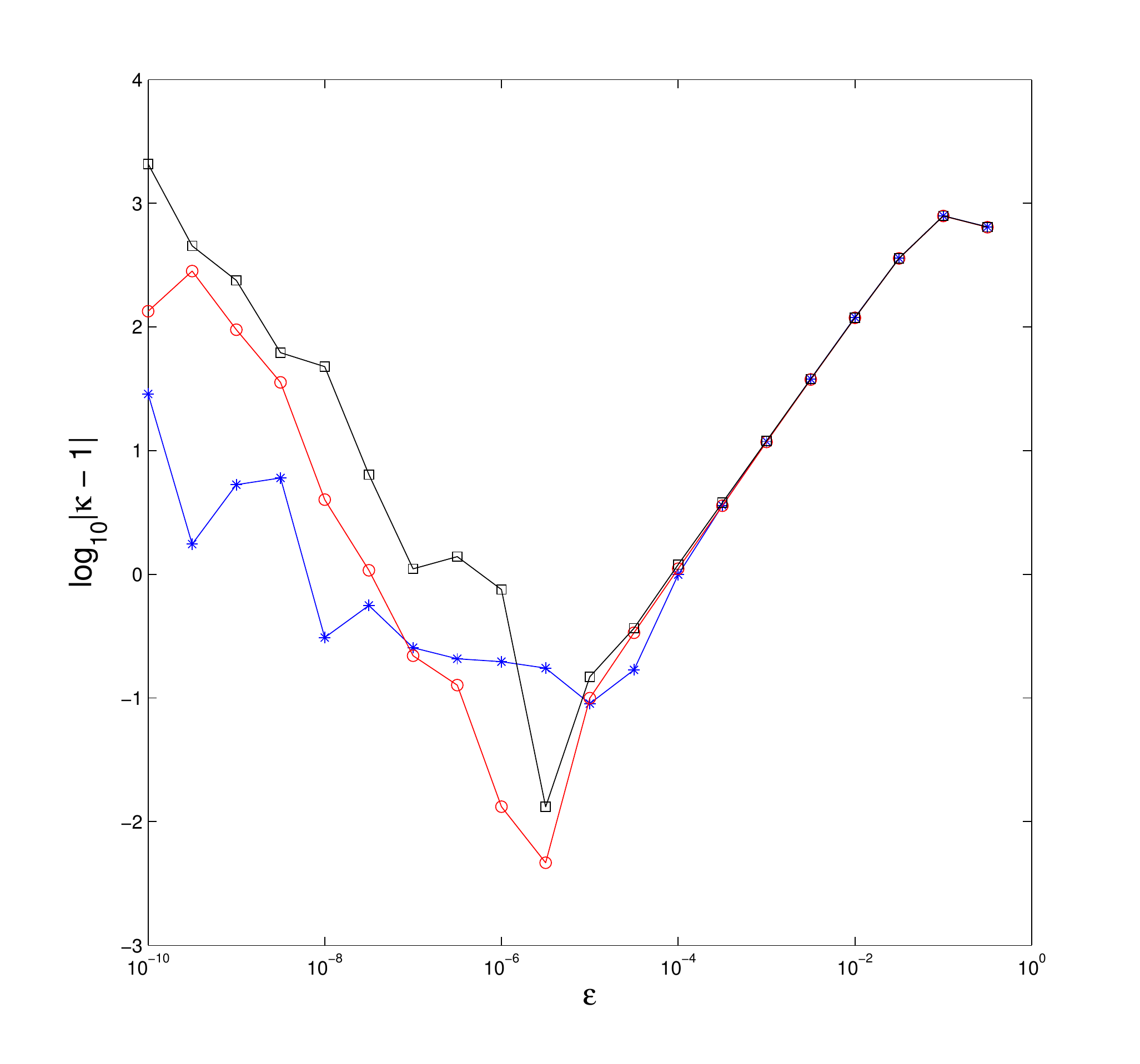}}}
\caption{TEST \#1 (Table \ref{tab:kappa}): Dependence of
  $\log_{10}|\kappa(\epsilon) - 1|$ on the step size $\epsilon$ in
  \eqref{eq:kappa} for different perturbations (a) $\zeta_1$, (b)
  $\zeta_2$, (c) $\zeta_3$ and (d) $\zeta_4$, cf.~Equation
  \eqref{eq:zeta}, and different resolutions (asterisks) $N=50,\,
  M=50$, (circles) $N=100, \, M=100$, and (squares) $N=80, \, M=300$.}
\label{fig:kappa1}
\mbox{
\subfigure[]{
\includegraphics[scale=0.27]{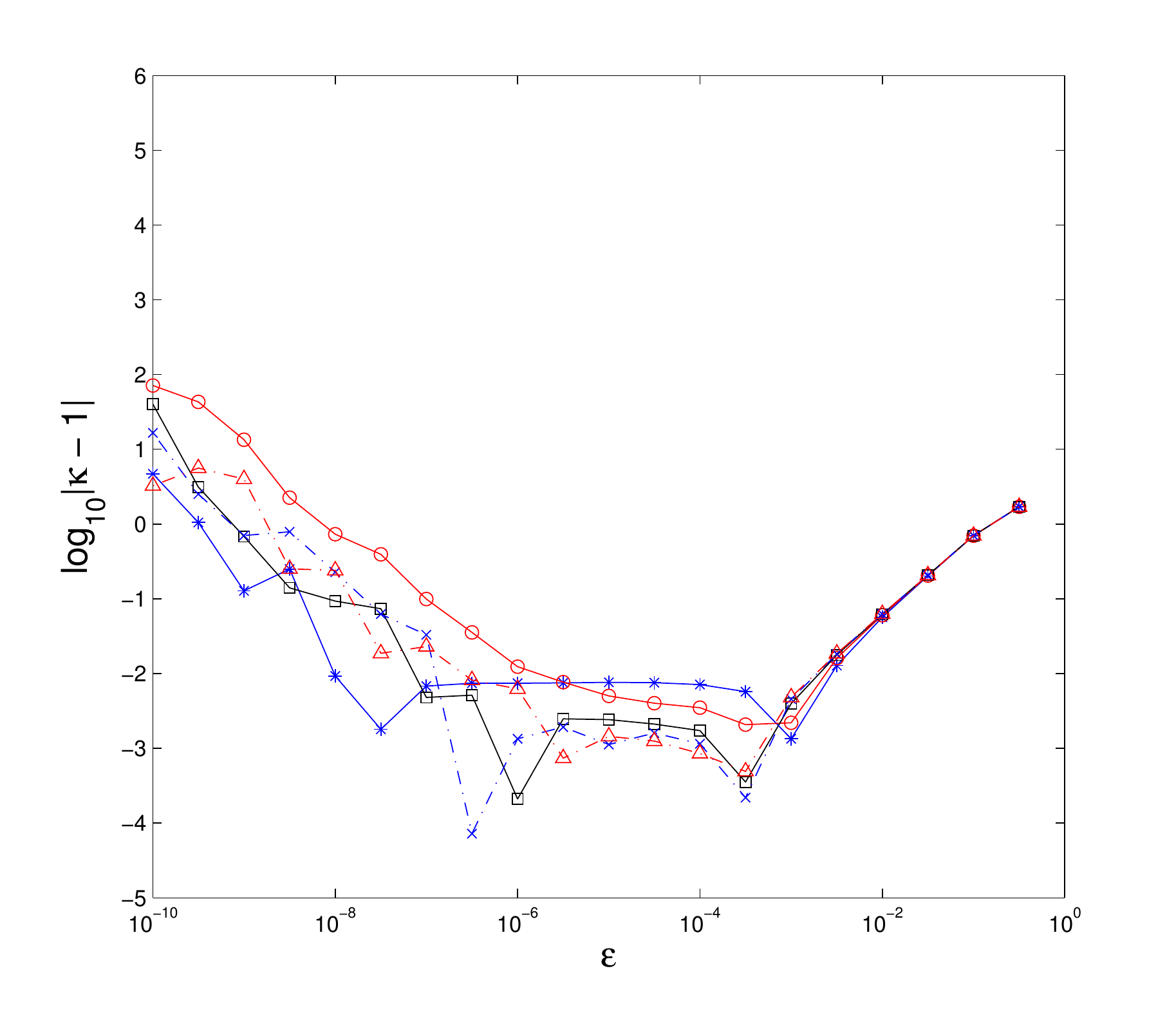}}
\subfigure[]{
\includegraphics[scale=0.27]{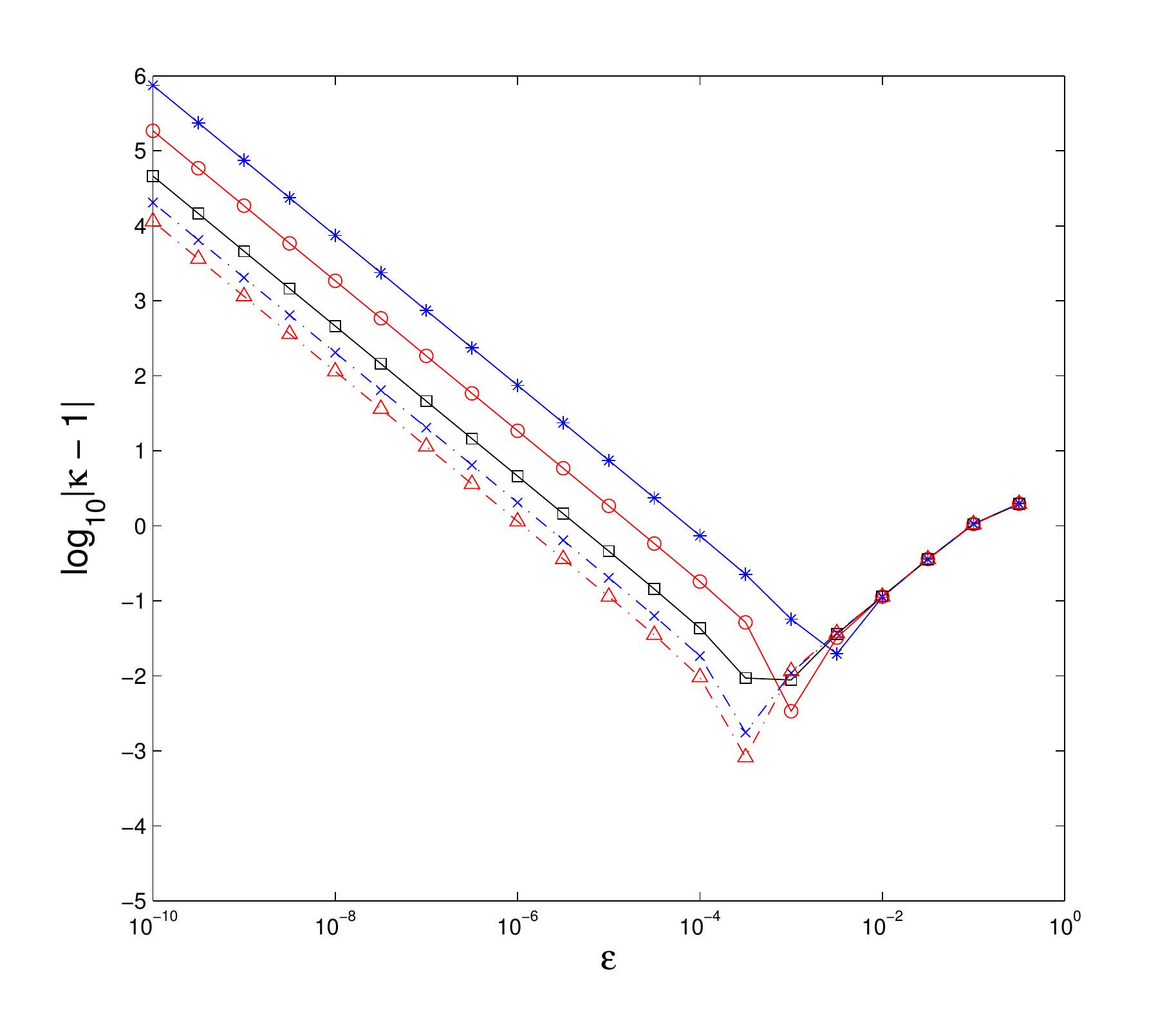}}}
\mbox{
\subfigure[]{
\includegraphics[scale=0.235]{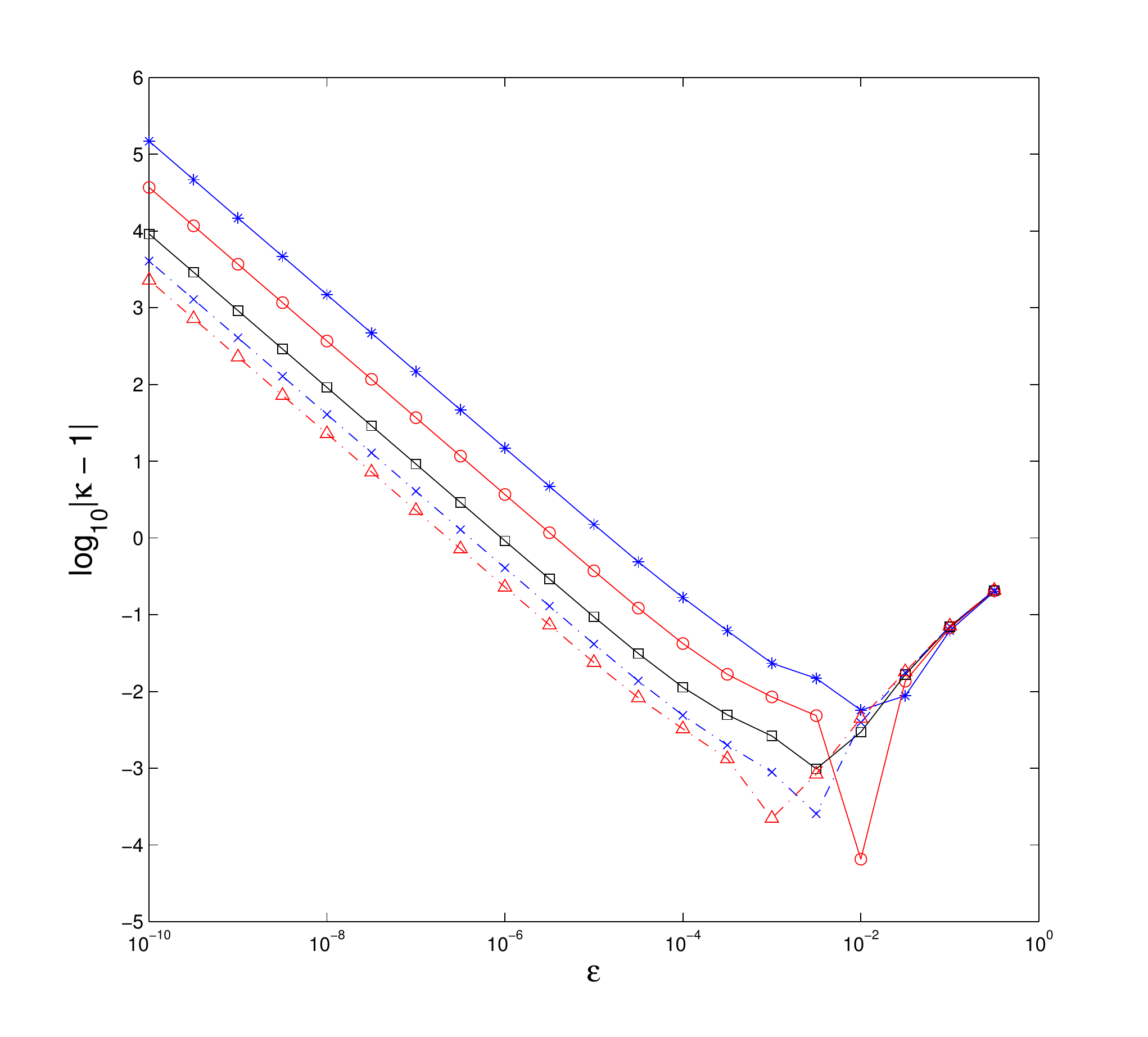}}
\subfigure[]{
\includegraphics[scale=0.275]{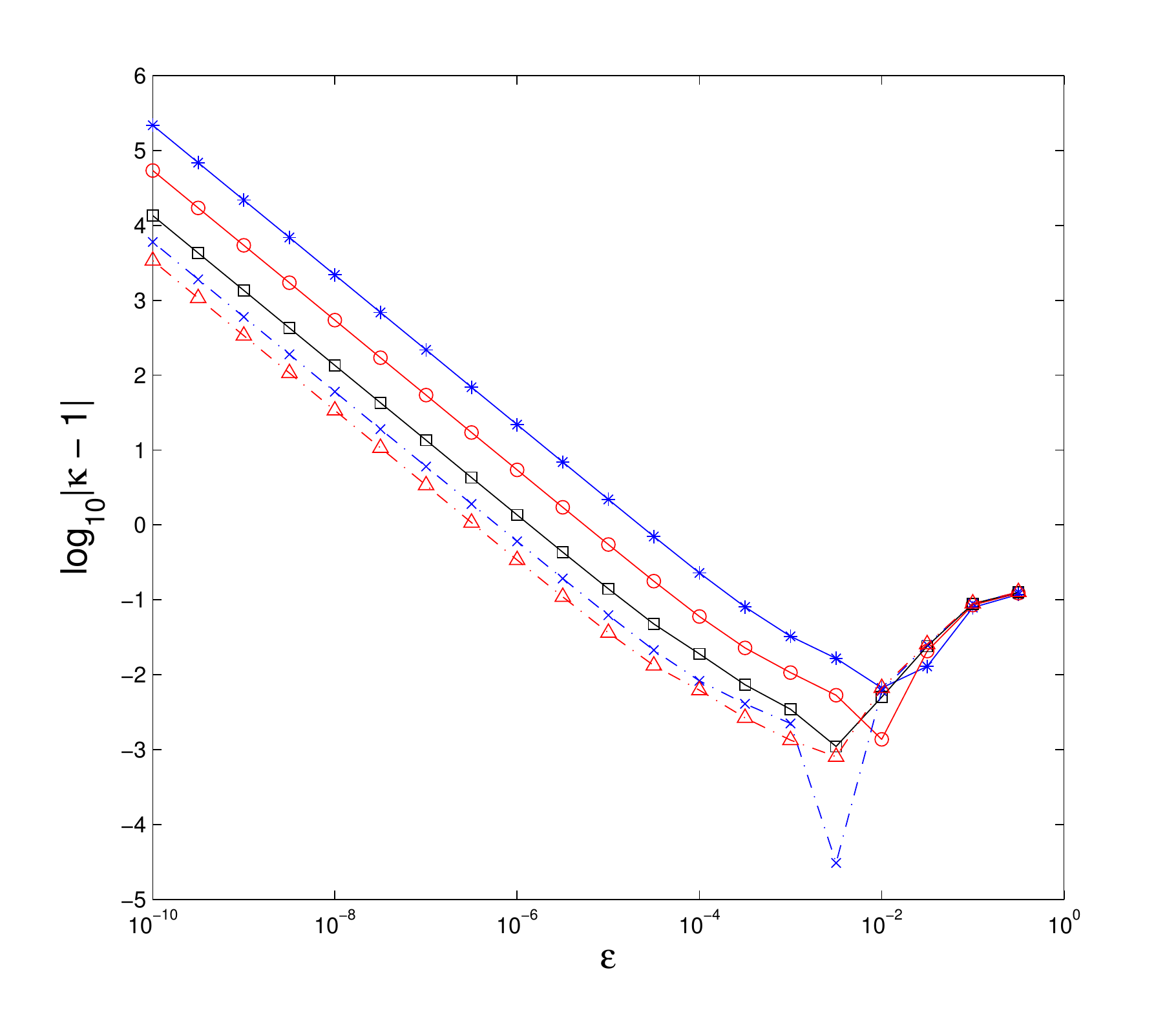}}}
\caption{TEST \#2 (Table \ref{tab:kappa}): Dependence of
  $\log_{10}|\kappa(\epsilon) - 1|$ on the step size $\epsilon$ in
  \eqref{eq:kappa} for different contours (a) $\C_2$, (b) $\C_3$, (c)
  $\C_4$ and (d) $\C_5$, cf.~Table \ref{tab:C}, and different
  resolutions (asterisks) $N=50,\, M=50$, (circles) $N=80, \, M=100$,
  (squares) $N=80, \, M=200$, (crosses) $N=80, \, M=300$ and
  (triangles) $N=80, \, M=400$.}
\label{fig:kappa2}
\end{figure}

A standard computational test employed to ascertain the accuracy of
the cost functional gradients in PDE optimization problems is to
calculate the G\^ateaux differential $\J'(\C;\zeta\n)$ for a given
contour $\C$ and its perturbations $\zeta$ in two different ways:
using an approximate finite-difference formula and Riesz identity
\eqref{eq:rieszH1} \cite{hnl02}. Thus, the ratio of these two
expressions, denoted
\begin{equation}
\kappa(\eps) \triangleq \frac{\J(\C(\epsilon)) - \J(\C(0))}{\epsilon \, 
\big\langle \bnabla^{L_2}\J(\C(0), \zeta \big\rangle_{L_2(0,L)}},
\label{eq:kappa}
\end{equation}
should be approximately equal to unity for a range of values of
$\epsilon$. Plotting $|\kappa(\epsilon) - 1|$ using the logarithmic
scale allows one to see the number of significant digits of accuracy
captured in the computation. We remark that, since different Riesz
representations ($L_2$ vs.~$H^1$) give the same differential
$\J'(\C;\zeta\n)$, for simplicity in \eqref{eq:kappa} we can use the
$L_2$ inner product together with the corresponding gradient. To focus
attention, we present our validation results for the functional
$\J(\C)$, i.e., without the length constraint, as the gradient of the
latter part does not involve the adjoint variable $u^*$. We analyze
two sets of results: one in which we fix the contour $\C$ and consider
different perturbations $\zeta$ and vice versa. For every pair of the
contour and the perturbation we study the effect of different
resolutions $N$ and $M$.  Details concerning the two test cases are
collected in Table \ref{tab:kappa}, where {the different contours
  are specified in Table \ref{tab:C}, whereas} the perturbations
tested are given by
\begin{equation}
\zeta_j(t) = \sin(j\, t), \quad t = [0,2\pi], \  j=1,2,3,4.
\label{eq:zeta}
\end{equation}
In both validation tests we assume that $\A = \Omega$ and use the
following distribution of the heat sources and the target temperature
profile
\begin{align}
q(x,y) &= 50-15x^2-15(y-0.5)^2, \label{eq:q} \\
\ubar(x,y) &= 15+\sin(4x-1)\cos(4y-1), \label{eq:ubar}
\end{align}
where $-1 \le x,y \le 1$. The results of TEST \#1 and TEST \#2 are
shown in Figures \ref{fig:kappa1} and \ref{fig:kappa2}, respectively.
In both cases we note that $\kappa(\epsilon)$ is fairly close to the
unity for values of $\epsilon$ spanning several orders of magnitude.
The quantity $\kappa(\epsilon)$ deviates from the unity for very small
values of $\epsilon$ which is due to the subtractive cancellation
(round-off) errors, and for large values of $\epsilon$ which is due
to the truncation errors, both of which are well-known effects
\cite{pl08}. Since we use the ``differentiate-then-discretize''
formulation, one should not expect $|\kappa(\epsilon) - 1|$ to be at
the level of the machine precision, although this quantity 
approaches zero as the resolution is refined. We also tested cases in
which $\A \neq \Omega$ and the length constraint was included
obtaining similar results as in Figures \ref{fig:kappa1} and
\ref{fig:kappa2}. Having thus validated the cost functional gradients,
we now move on to discuss solution of the actual optimization
problems.
\begin{table}
\begin{center}
  \caption{Settings for the validation tests of the cost functional gradient $\bnabla \J$ {(Problem P1)}.
    The contours and perturbations used are defined in Table \ref{tab:C} and equation \eqref{eq:zeta}, respectively.} 
  \begin{tabular}{| c | c | c | c | c |}
    \hline
    TEST & Contour $\C$ & Perturbations $\zeta$ &Resolution $(N,~M) $ & \Bmp{2.0cm} \vspace*{0.1cm}  Target \\ Domain \vspace*{0.1cm} $\A$ \Emp \\ \hline
    $\#1$ & $\C_1$ & $\zeta_1$, $\zeta_2$, $\zeta_3$, $\zeta_4$ & \Bmp{3.0cm} \vspace*{0.1cm} $(50,50)$, $(100,100)$, $(80,300)$ \vspace*{0.1cm}  \Emp &$ \Omega$ \\ \hline
    $\#2$ & $\C_2$, $\C_3$, $\C_4$, $\C_5$ & $\zeta_1$ & \Bmp{3.0cm} \vspace*{0.1cm} $(50,50)$, $(80,100)$, $(80,200)$, $(80,300)$, $(80,400)$ \vspace*{0.1cm} \Emp & $\Omega$ \\ \hline
  \end{tabular}
  \label{tab:kappa}
\vspace*{0.4cm}
\caption{Definitions of contours $\C_1,\dots,\C_7$ used in the different cases studied in Section \ref{sec:results}.} 
  \begin{tabular}{| c | c | c |}
    \hline
    \Bmp{1.5cm} \vspace*{0.1cm} \centering Contour \vspace*{0.1cm} \Emp & \Bmp{8.0cm} \vspace*{0.1cm} \centering Parametrization $(0 \le t \le 2\pi)$ \vspace*{0.1cm} \Emp & 
    \Bmp{1.5cm} \vspace*{0.1cm} \centering Plot \vspace*{0.1cm} \Emp  \\ \hline
    $\C_1$  & \Bmp{8.0cm} \vspace*{0.1cm} \centering $x(t) = 0.4 \cos(t) + 0.1$,  $y(t) = 0.4 \sin(t) - 0.1$ \vspace*{0.1cm} \Emp & 
    \Bmp{1.5cm} \vspace*{-0.15cm} \centering \includegraphics[scale=0.1]{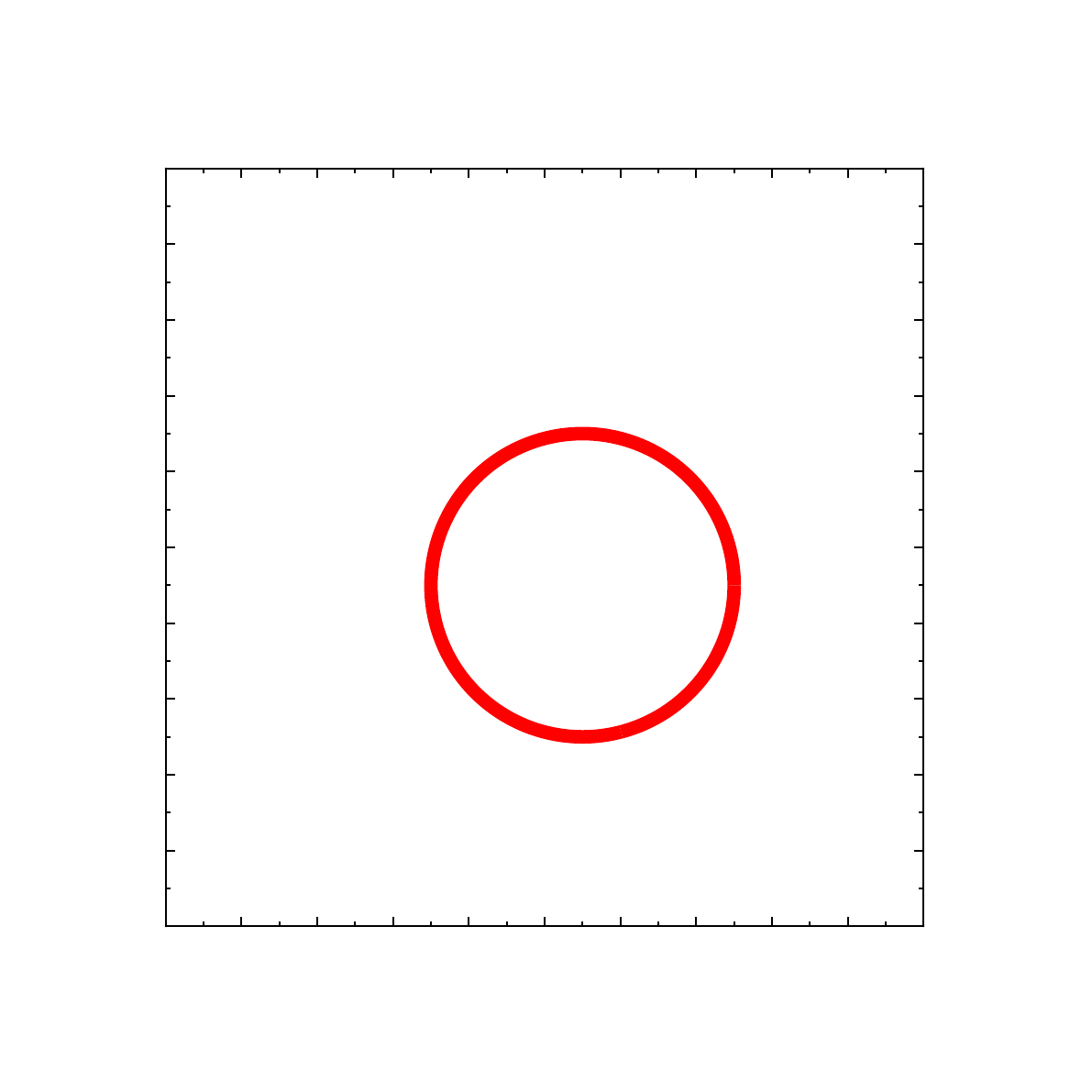} \\ \vspace*{-0.12cm} \Emp \\ \hline
    $\C_2$  & \Bmp{8.0cm} \vspace*{0.1cm} \centering $x(t) = 0.2 \cos(t) + 0.4$,  $y(t) = 0.2 \sin(t) + 0.4$ \vspace*{0.1cm} \Emp & 
    \Bmp{1.5cm} \vspace*{-0.15cm} \centering \includegraphics[scale=0.1]{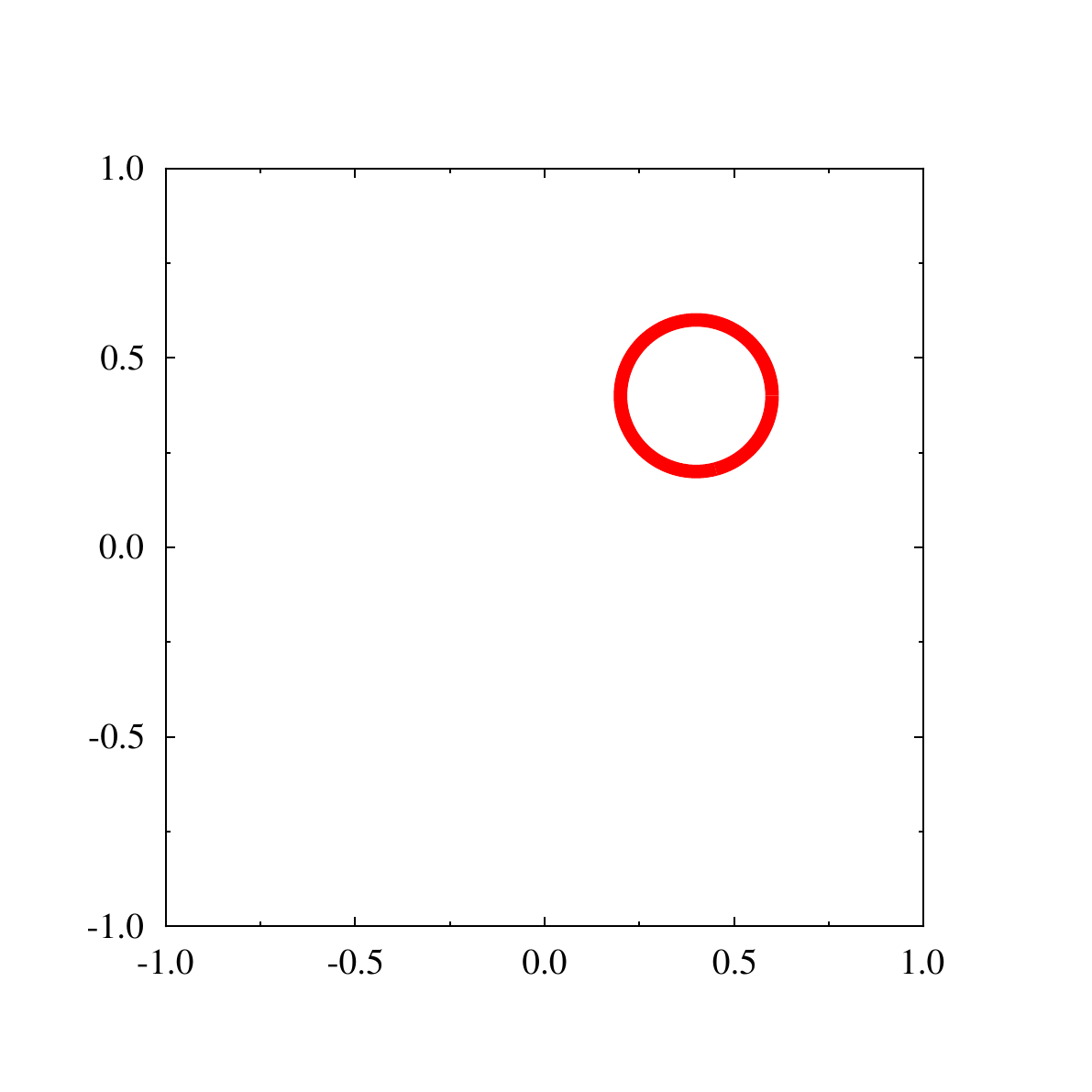} \\ \vspace*{-0.12cm} \Emp \\ \hline
    $\C_3$  & \Bmp{8.0cm} \vspace*{0.1cm} \centering $x(t) = 0.3 \cos(t)$,  $y(t) = 0.2 \sin(t)$ \vspace*{0.1cm} \Emp & 
    \Bmp{1.5cm} \vspace*{-0.15cm} \centering \includegraphics[scale=0.1]{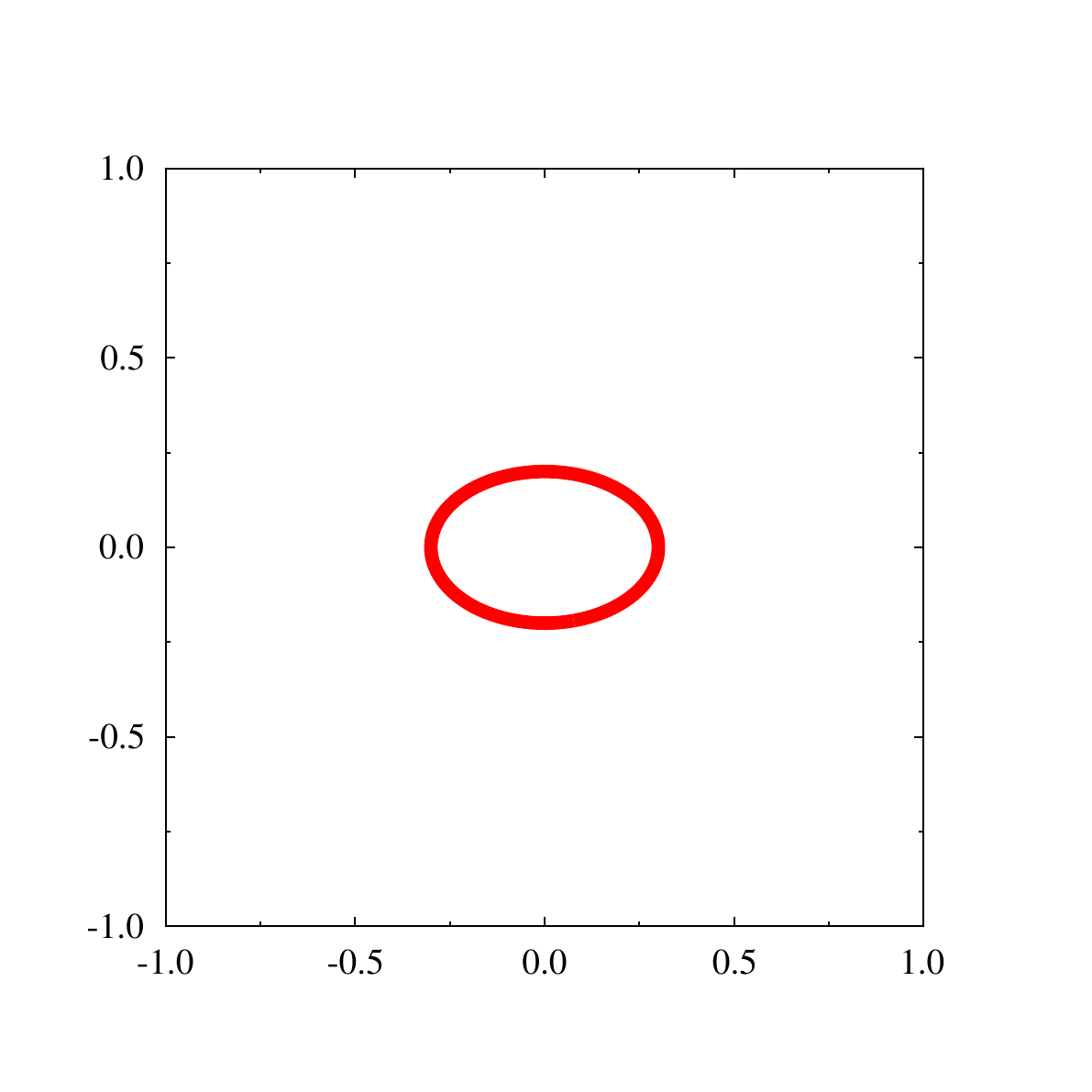} \\ \vspace*{-0.12cm} \Emp \\ \hline
    $\C_4$  & \Bmp{8.0cm} \vspace*{0.1cm} \centering $x(t) = 0.4(1+0.1\cos(3t))\cos(t) + 0.1$,  \\ $y(t) = 0.4(1+0.1\cos(3t))\sin(t) + 0.1$ \vspace*{0.1cm} \Emp & 
    \Bmp{1.5cm} \vspace*{-0.15cm} \centering \includegraphics[scale=0.1]{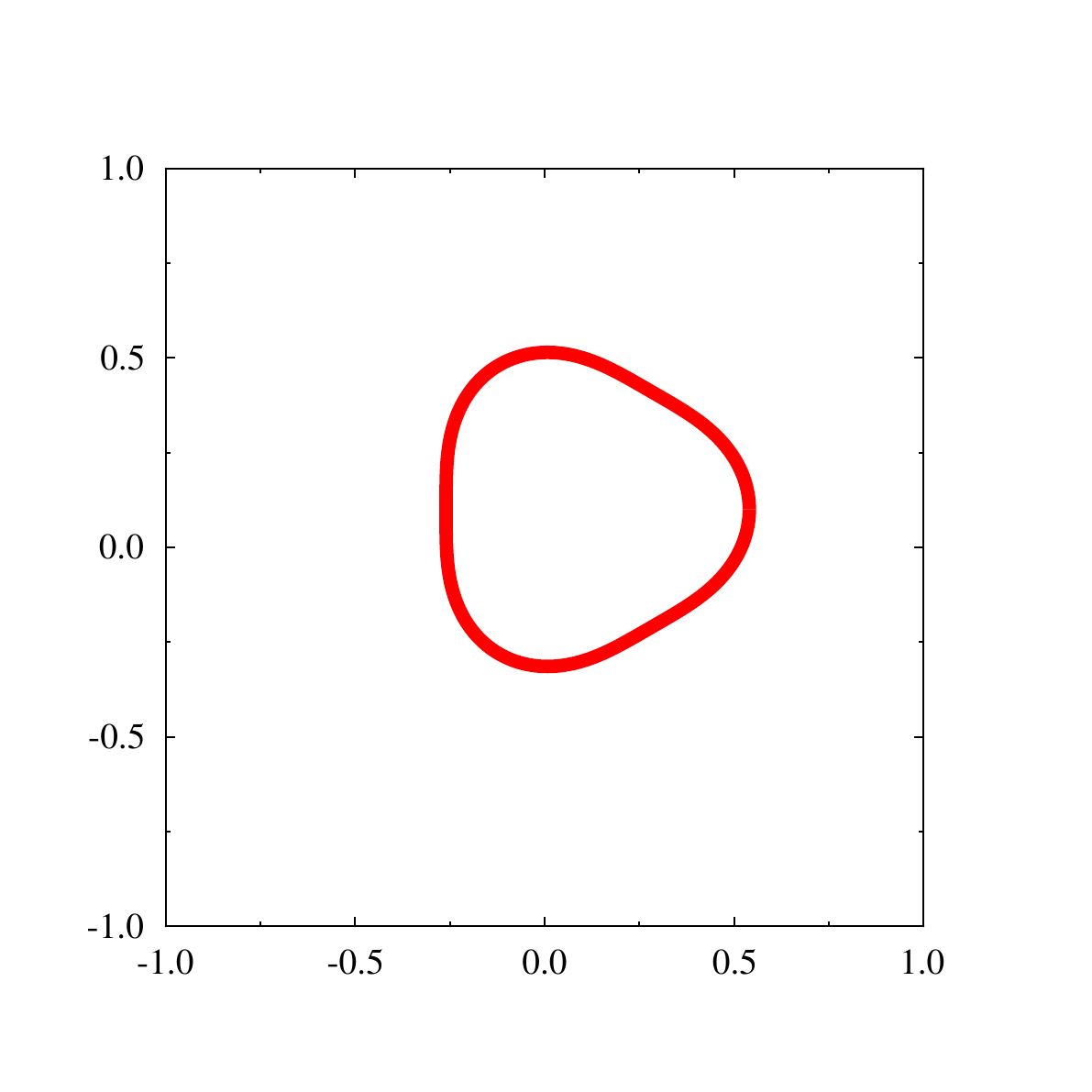} \\ \vspace*{-0.12cm} \Emp \\ \hline
    $\C_5$  & \Bmp{8.0cm} \vspace*{0.1cm} \centering $x(t) = 0.4(1+0.1\cos(4t))\cos(t)+0.1$,  \\ $y(t) = 0.4(1+0.1\cos(4t))\sin(t)+0.1$ \vspace*{0.1cm} \Emp & 
    \Bmp{1.5cm} \vspace*{-0.15cm} \centering \includegraphics[scale=0.1]{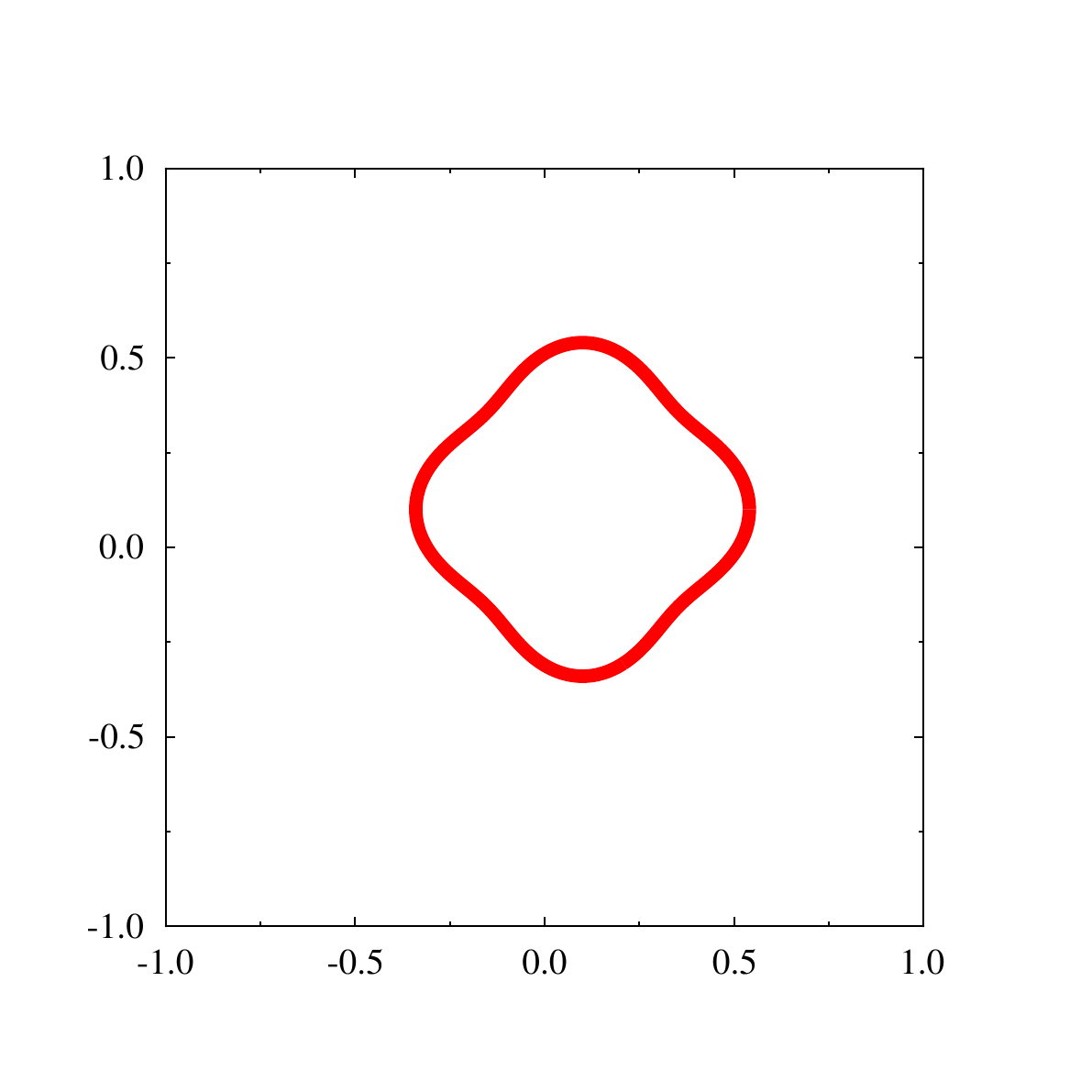} \\ \vspace*{-0.12cm} \Emp \\ \hline
    $\C_6$  & \Bmp{8.0cm} \vspace*{0.1cm} \centering $x(t) = \frac{3}{2\pi} \cos(t)-0.4$,  $y(t) = \frac{3}{2\pi} \sin(t)+0.3$ \vspace*{0.1cm} \Emp & 
    \Bmp{1.5cm} \vspace*{-0.15cm} \centering \includegraphics[scale=0.1]{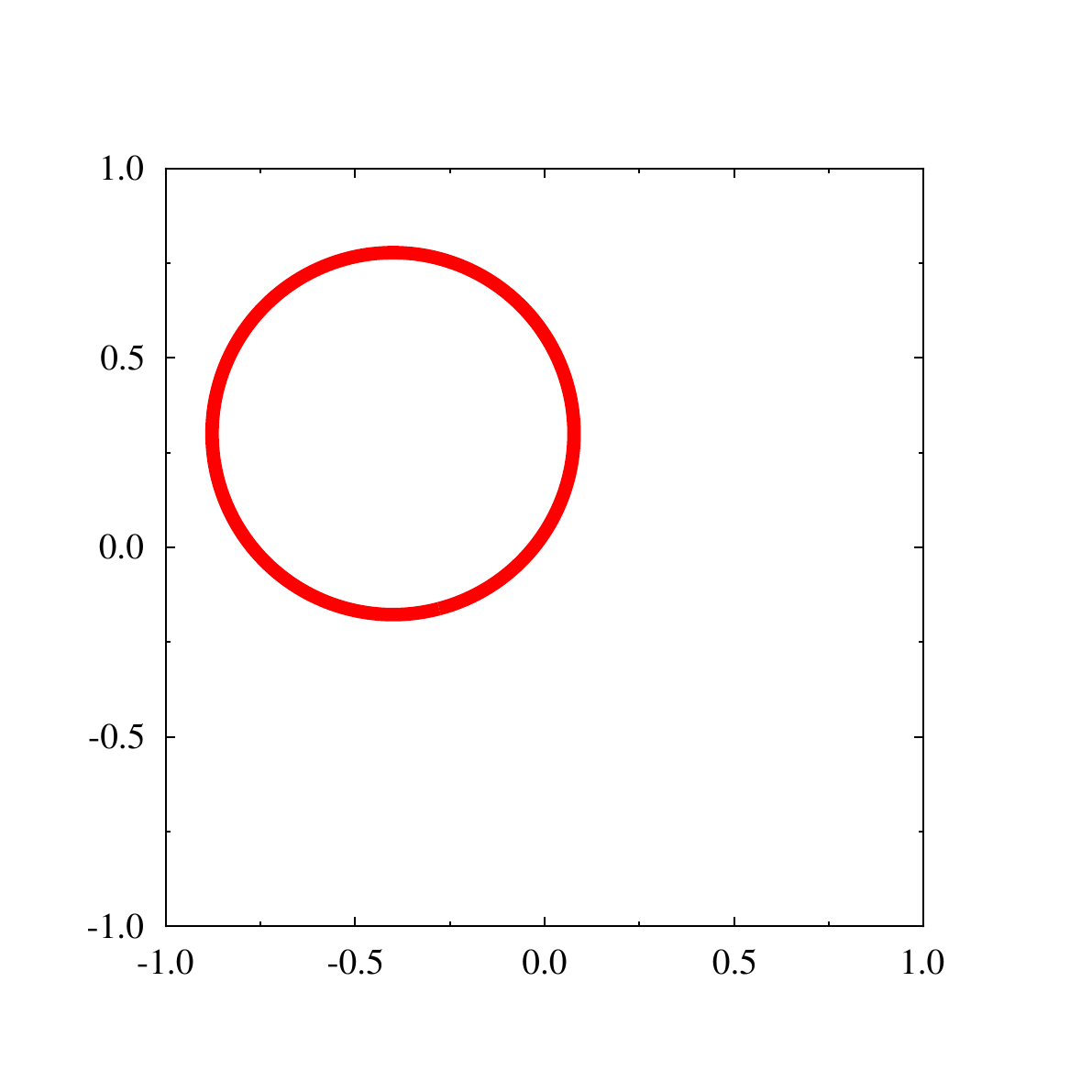} \\ \vspace*{-0.12cm} \Emp \\ \hline
    $\C_7$  &  \Bmp{8.0cm} \vspace*{0.1cm} \centering $x(t) = \frac{t - \pi}{\pi}$,  $y(t) = 0.78$ \vspace*{0.1cm} \Emp  & 
    \Bmp{1.5cm} \vspace*{-0.15cm} \centering \includegraphics[scale=0.1]{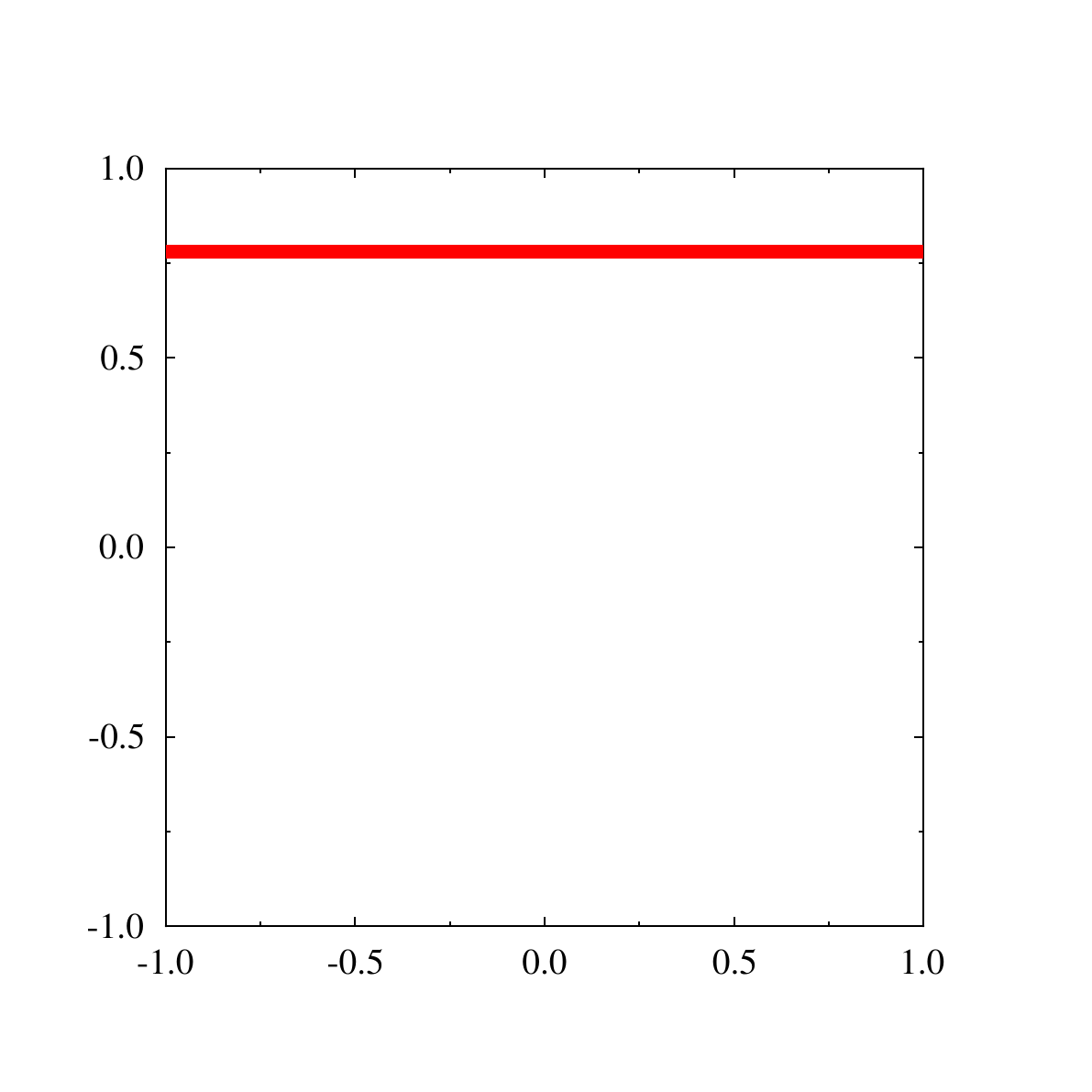} \\ \vspace*{-0.12cm} \Emp \\ \hline
\end{tabular}
  \label{tab:C}
\end{center}
\end{table}

\subsection{Solution of Optimization Problems}
\label{sec:optim}

\begin{figure}
\centering
\mbox{
\subfigure[]{
\includegraphics[scale=0.45]{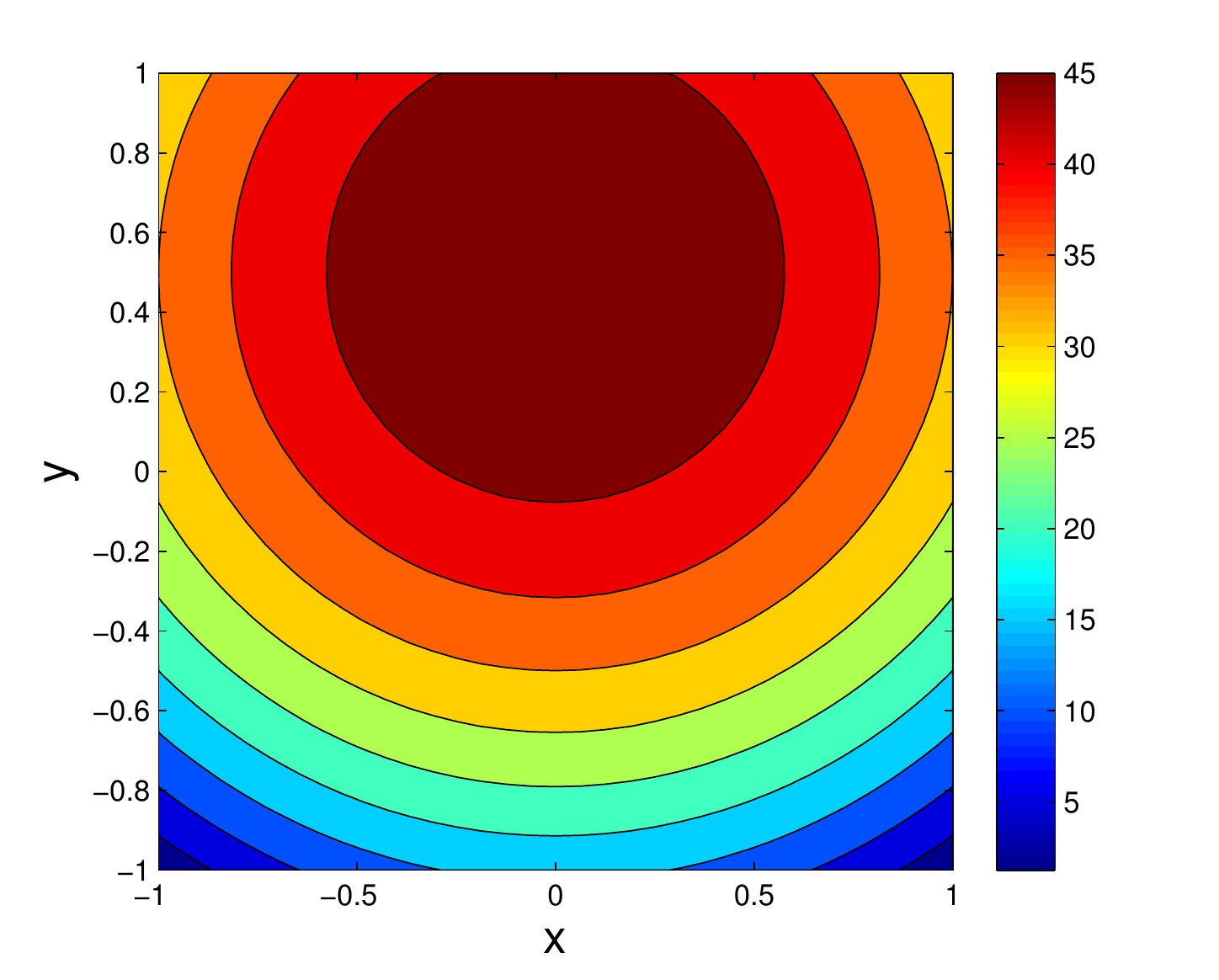}}
\subfigure[]{
\includegraphics[scale=0.45]{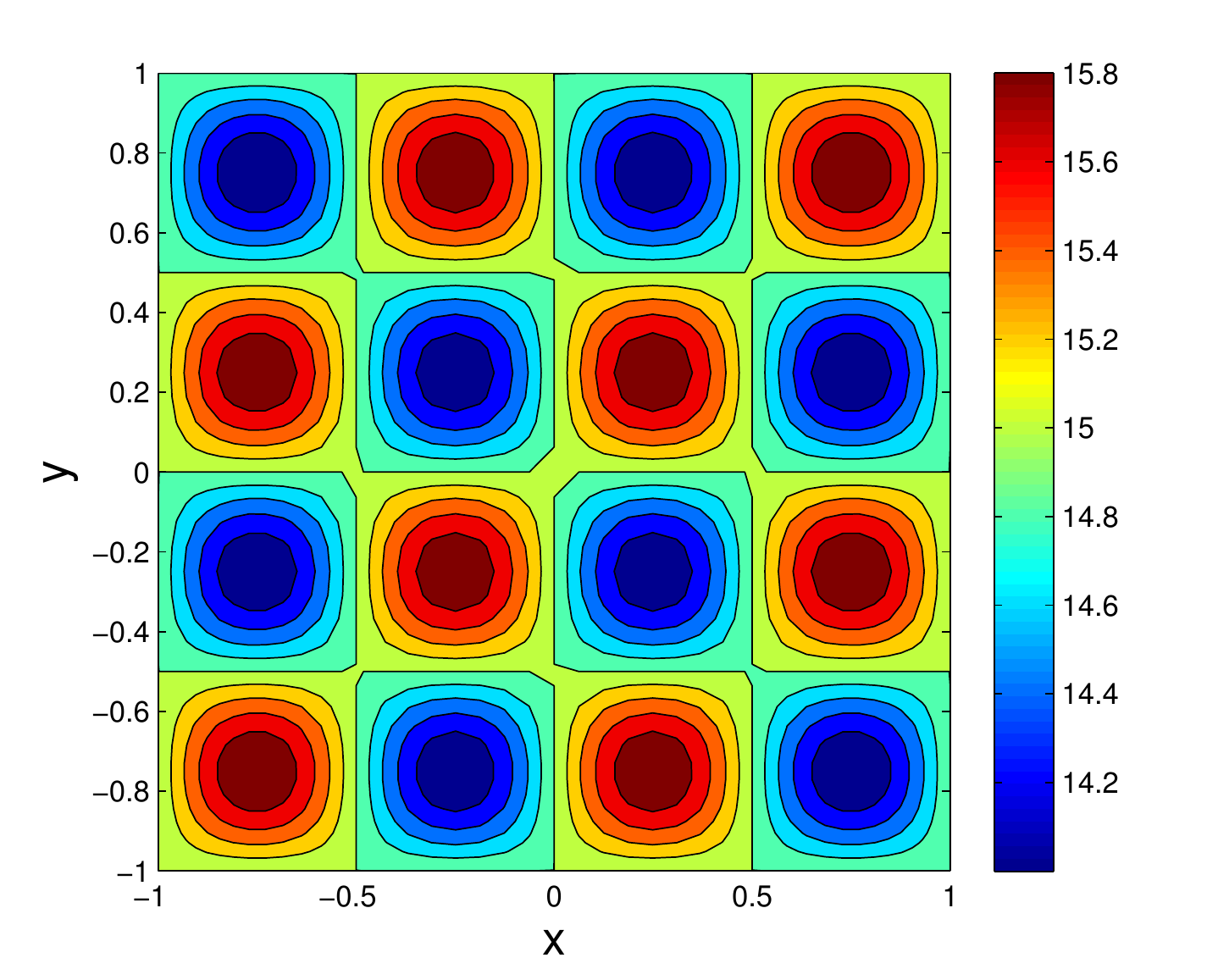}}}
\caption{(a) Distribution of heat sources $q$, cf.~\eqref{eq:q2}, 
and (b) target temperature field $\ubar$, cf.~\eqref{eq:ubar2a}, used
in CASE \#1 and \#2.}
\label{fig:qub}
\end{figure}

We will study in detail solution of the following three optimization
problems with and without the length constraint, as indicated below:
in CASE \#1 {for Problem P1} we examine the convergence of
Algorithm \ref{alg1} without the length constraint for several
different initial guesses $\C^{(0)}$ for the contour and using $\A =
\Omega$, in CASE \#2 {for Problem P1} we consider {a
  configuration in which $\A \neq \Omega$} and also study the effect
of the length constraint, and in CASE \#3 {for Problem P2} we
investigate {a system in which} the heat source distribution $q$
corresponds to the temperature field in an actual battery cell, also
{in the presence} of the length constraint. Parameters
characterizing the three cases are collected in Table \ref{tab:cases}.
As concerns the heat source distribution $q$, in CASES \#1 and \#2 it
is given by the following expression (Figure \ref{fig:qub}a)
\begin{equation}
q(x,y)=50-15x^2-15\left(y-\frac{1}{2}\right)^2\qquad\,(x,y)\in\, \ \Omega,
\label{eq:q2}
\end{equation}
whereas in CASE \#3 it is obtained (by applying the Laplace operator
and suitable smoothing) to the temperature distribution determined
experimentally in an actual battery cell \cite{ksk09}, see Figure
\ref{fig:case3}a. In the different cases the target temperature field
is given by the following expressions (see also Figure \ref{fig:qub}b)
\begin{subequations}
\label{eq:ubar2}
\begin{alignat}{3}
&\textrm{CASE}\ \#1,2: \quad  & \ubar(x,y) & = 15 +\sin(2\pi\,x+\pi)\cos\left(2\pi\,y+\frac{\pi}{2}\right), \quad && (x,y)\in\, \ \A, \label{eq:ubar2a} \\
&\textrm{CASE}\ \#3: & \ubar(x,y) & = 30 && (x,y)\in\, \ \Omega. \label{eq:ubar2b}
\end{alignat}
\end{subequations}
In CASE \#1 and \#2 the distribution of heat sources \eqref{eq:q2} and
the target temperature field \eqref{eq:ubar2a} have been chosen to
{test the algorithm in the situation when the source field varies
  slowly, whereas the target field exhibits a significant
  variability,} cf.~Figures \ref{fig:qub}a and \ref{fig:qub}b.
{On the other hand,} in CASE \#3 the constant target temperature
field \eqref{eq:ubar2b} represents a typical engineering objective.
{The specific values assumed by the fields $q$ and $\overline{u}$
  do not have a physical significance and were selected to make the
  optimization problem sufficiently challenging.} The tolerances in
Algorithm \ref{alg1} are set to $\varepsilon_{\J}=10^{-3}$ and
$\varepsilon_{\tau}=10^{-8}$.
\begin{table}\footnotesize
\renewcommand{\arraystretch}{2}\addtolength{\tabcolsep}{-1pt}
\begin{center}
\caption{Parameters used in the solution of the three optimization problems in Section \ref{sec:optim}. {The contours used as the initial guesses $\C^{(0)}$ are defined in Table \ref{tab:C}.}}
     \begin{tabular}{ | c | c | c | c | c | c | c | c | c |}
      \hline
      CASE & $q$ & $\ubar$ & $(N,~M)$ & $\alpha$ & $L_0$ & $\ell$ & $\C^{(0)}$ & $\A$ \\ \hline
      \#1 (P1) & Eq.~\eqref{eq:q2} & Eq.~\eqref{eq:ubar2a} & (50,100) & 0 & --- & 0.1, 0.3 
           & \Bmp{1.25cm} \vspace*{0.05cm} \centering $\C_2$, $\C_3$, $\C_4$, $\C_5$ \vspace*{0.05cm} \Emp & $\Omega$ \\ \hline
      \#2 (P1) & Eq.~\eqref{eq:q2} & Eq.~\eqref{eq:ubar2a} & (50,100) 
           & $10^2$ & 3.0 & 0.1 & $\C_6$ &  $[-0.5,1]\times[-0.5,1]$ \\ \hline
      \#3 (P2) & Fig.~\ref{fig:case3}a & Eq.~\eqref{eq:ubar2b} & (50,100) 
           & \Bmp{1.0cm} \vspace*{0.2cm} \centering $1, 10$, \\ $10^2, 10^3$ \Emp & 6.0 & 0.1,0.25 & $\C_7$ & $\Omega$ \\ \hline 
    \end{tabular}
\label{tab:cases}
\end{center}
\end{table}

\begin{figure}
\centering
\mbox{
\subfigure[Initial contour shapes $\C_2$, $\C_3$, $\C_4$ and $\C_5$ 
(for clarity they are represented using fewer points than used in the actual computations)]{
\includegraphics[scale=0.34]{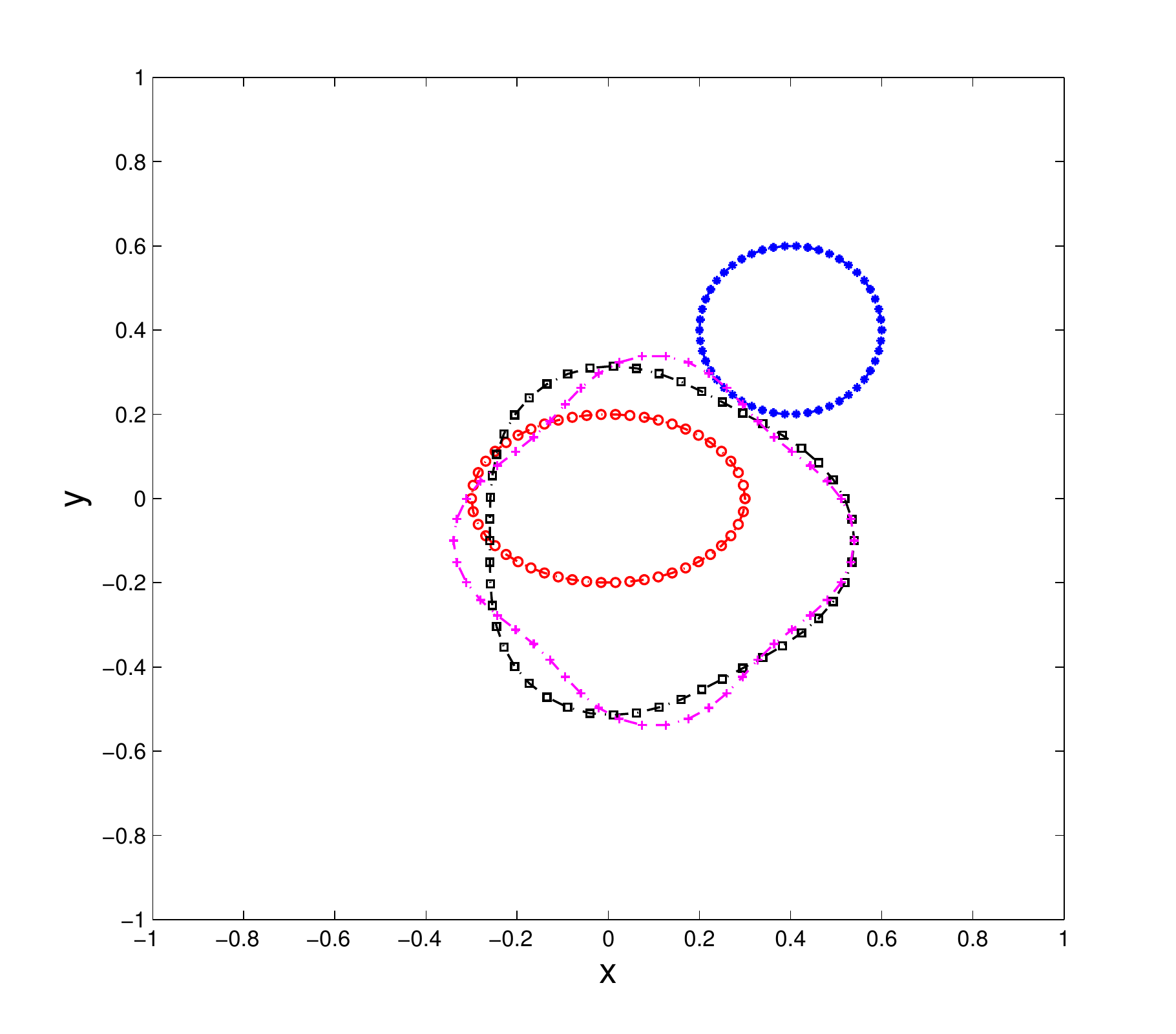}}\qquad
\hspace*{-0.3cm}
\subfigure[Optimal contour shapes $\tC$ corresponding to the 
initial contours shown in Figure (a)]{
\includegraphics[scale=0.335]{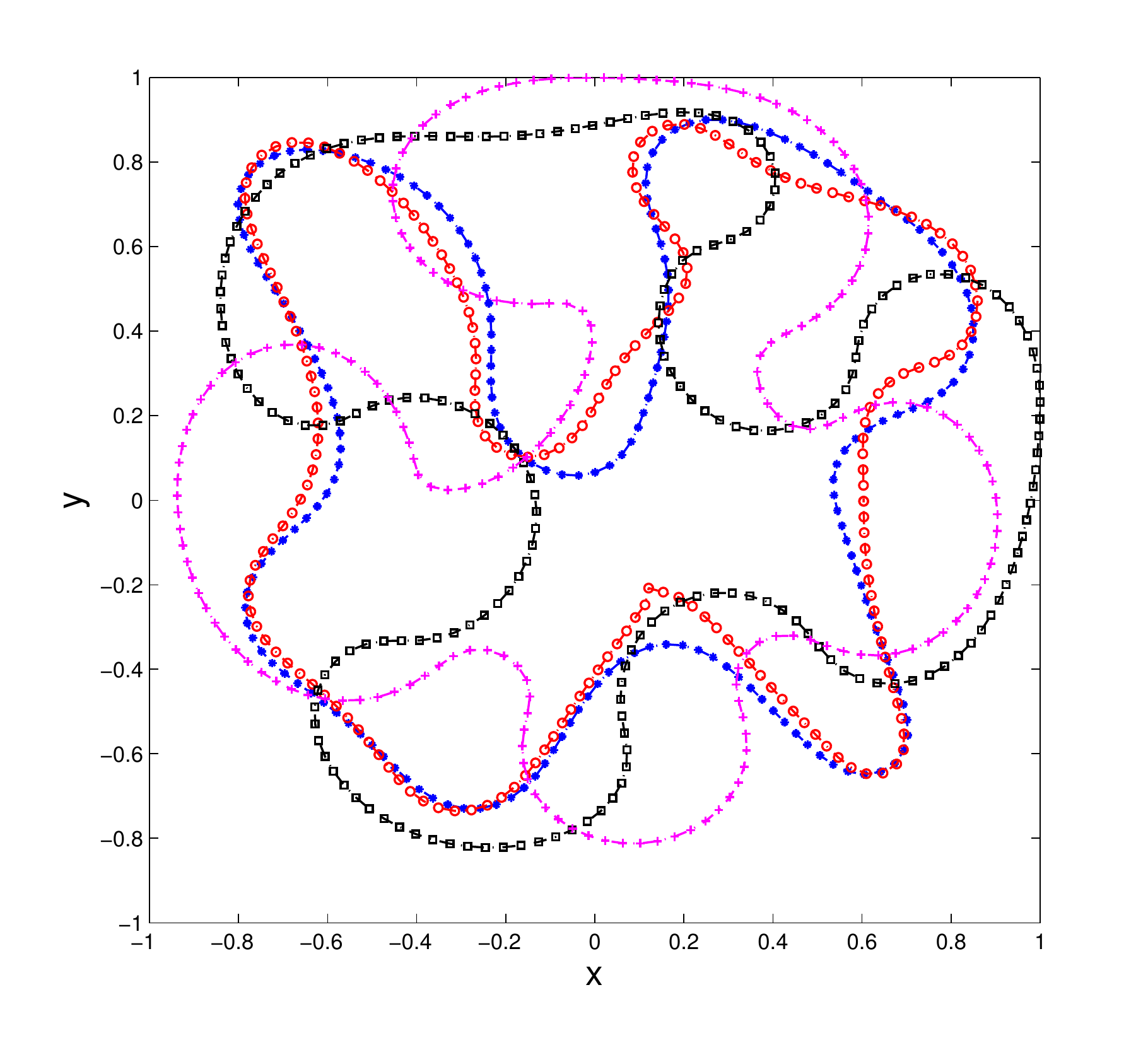}}}
\mbox{
\subfigure[Evolution of functionals $\J(\C^{(n)})$ for iterations 
starting with the initial contours shown in Figure (a)]{
\includegraphics[scale=0.33]{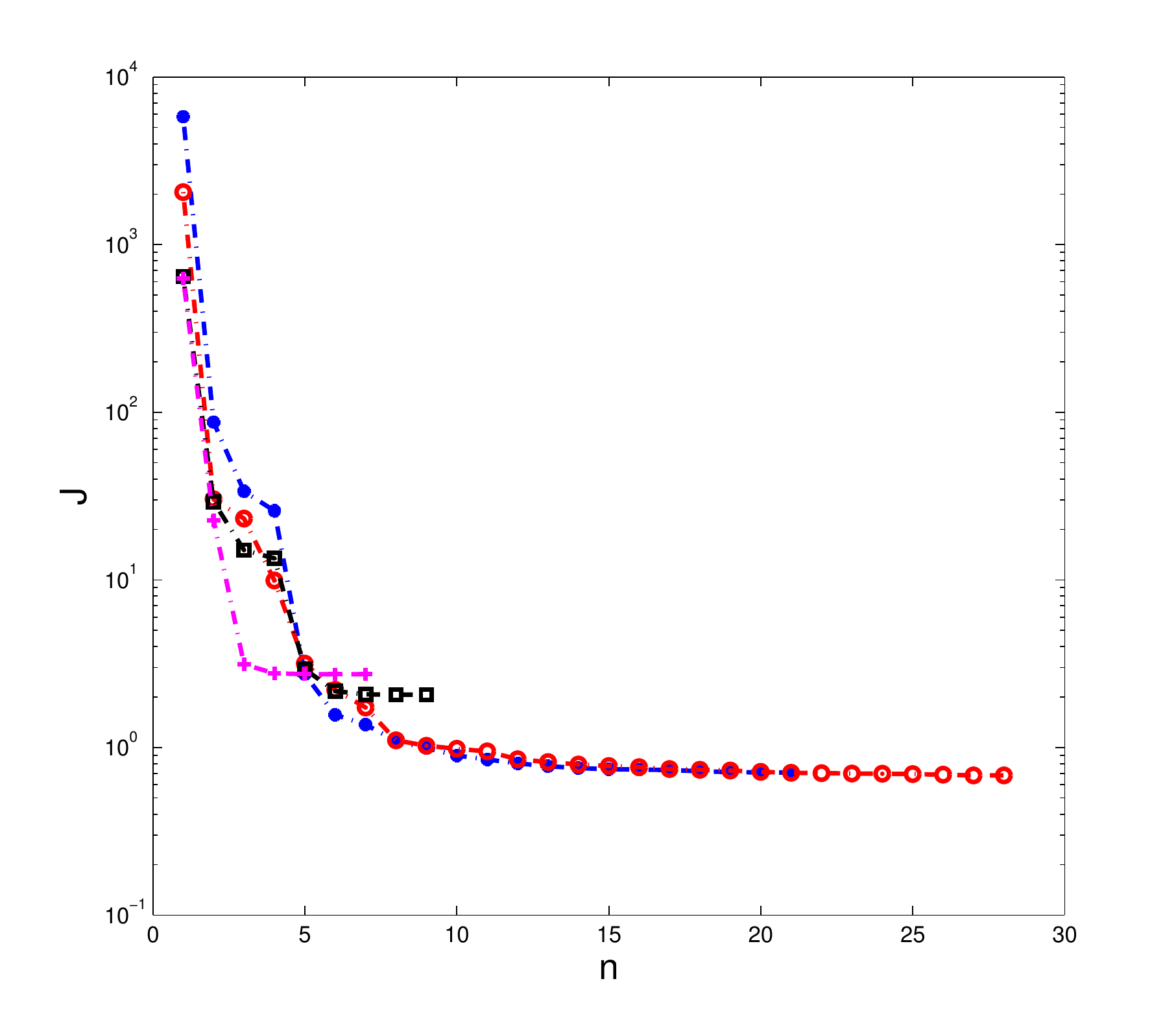}}\qquad
\subfigure[Evolution of contours $\C^{(n)}$ with iterations 
in the problem with initial shape $\C_2$ (dashed line)]{
\includegraphics[scale=0.32]{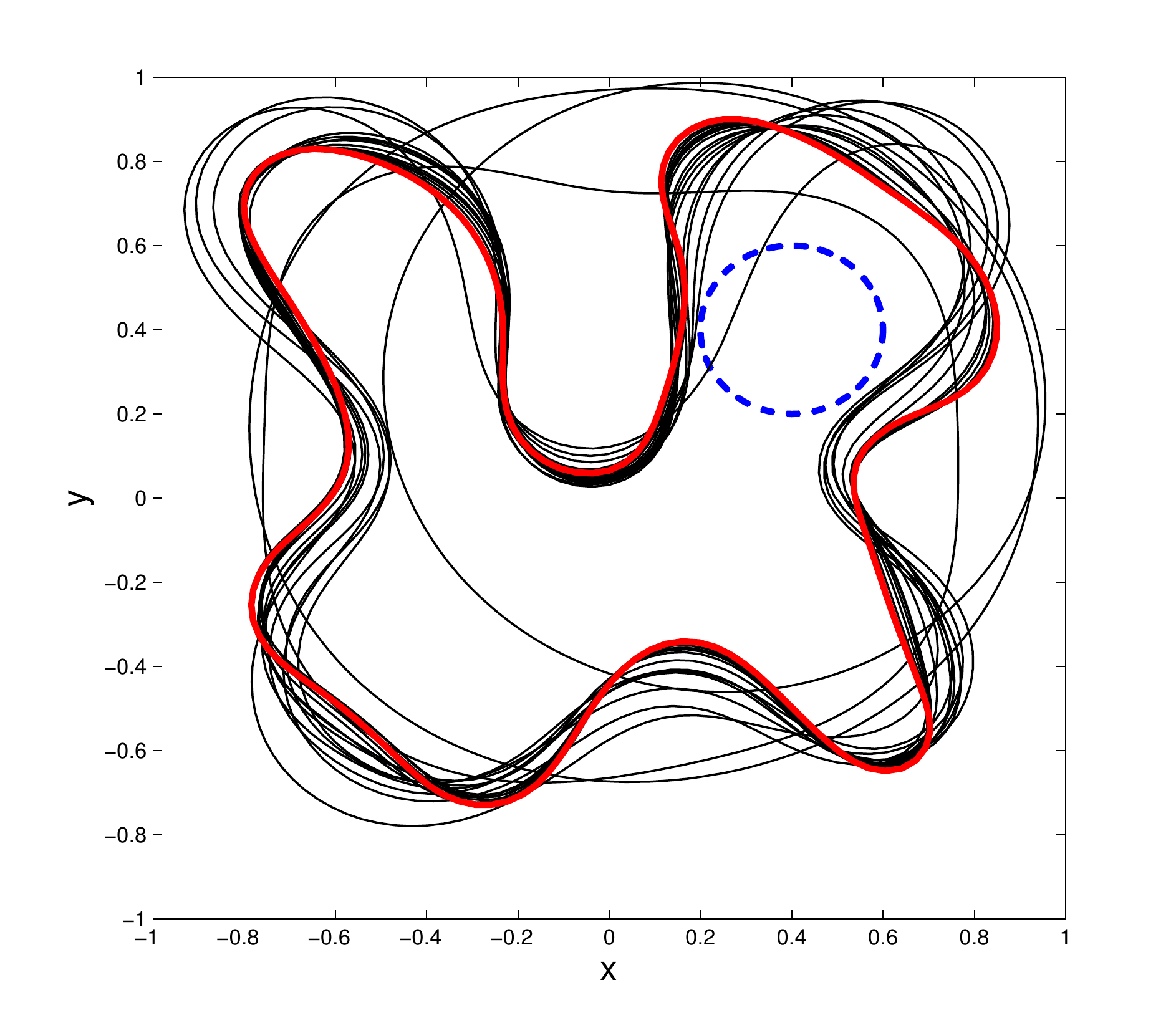}}}
\mbox{
\subfigure[Initial temperature field $u(\C_2)$]{
\includegraphics[scale=0.34]{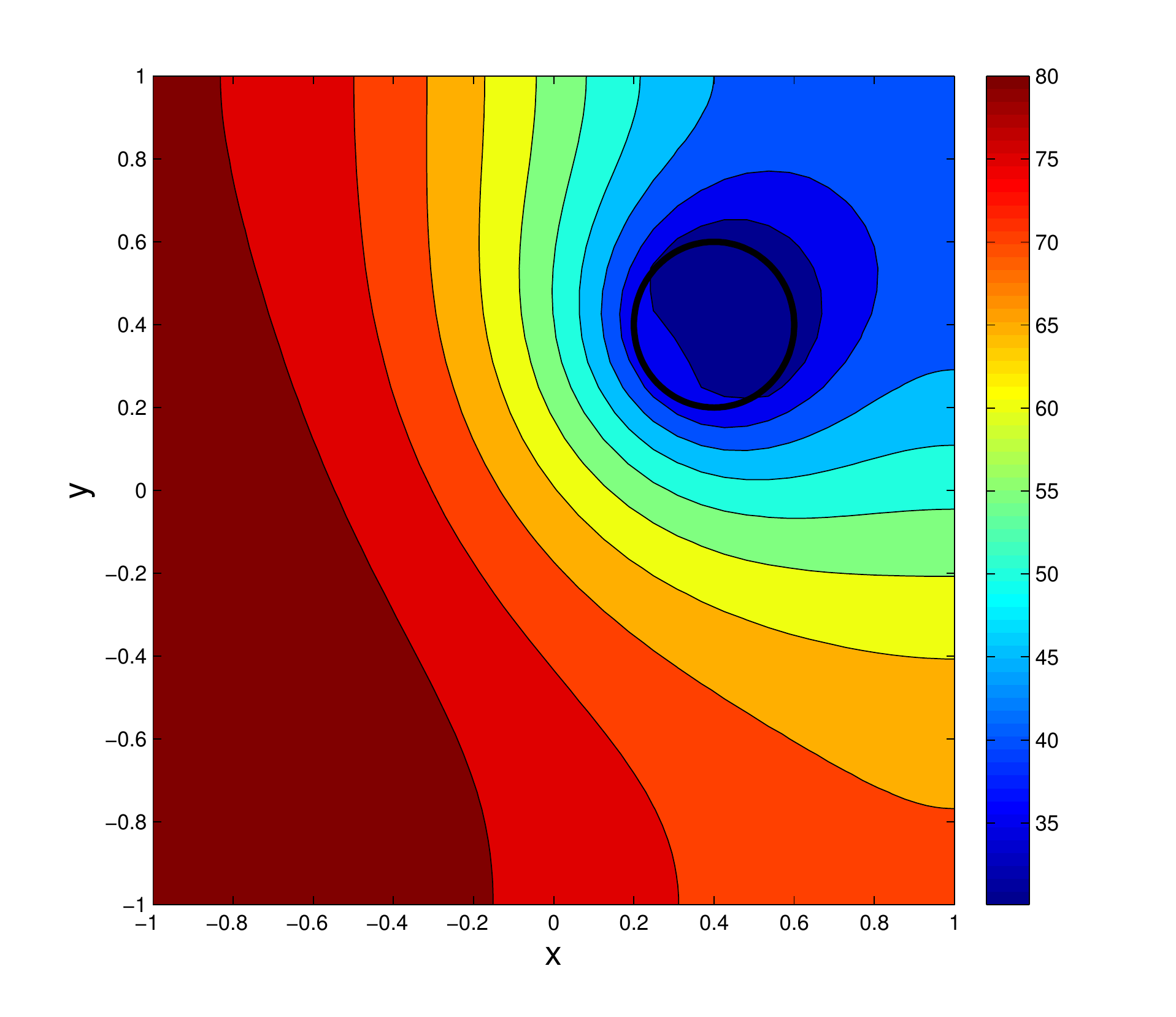}}\quad
\subfigure[Optimal temperature field $u(\tC)$ in the case with initial contour $\C_2$]{
\includegraphics[scale=0.35]{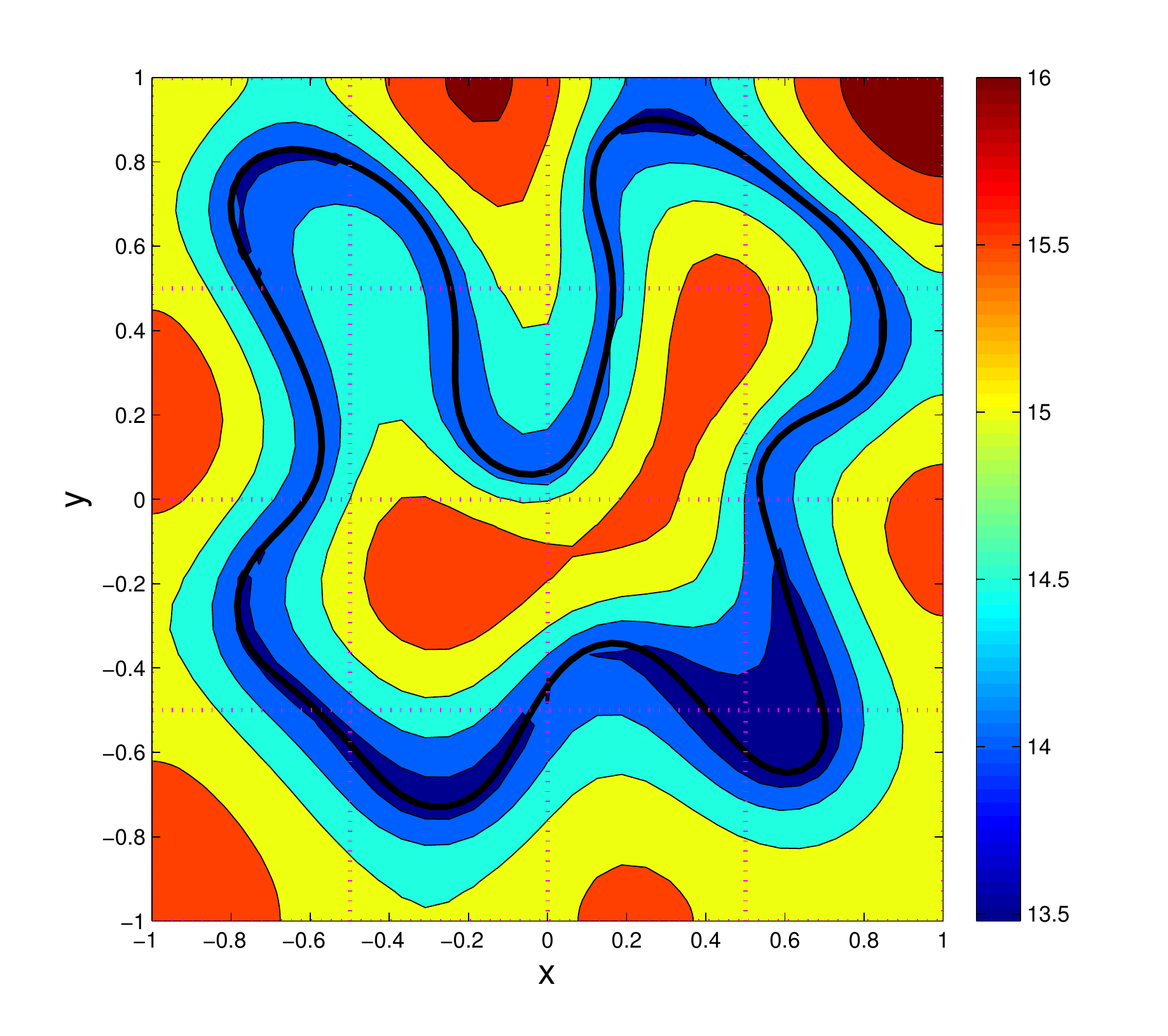}}}
\caption{Results illustrating solution of the optimization problems in
  CASE \#1, cf.~Table \ref{tab:cases}, using different initial
  contours (asterisks) $\C_2$, (circles) $\C_3$, (squares) $\C_4$,
  (pluses) $\C_5$, cf.~Table \ref{tab:C}.}
\label{fig:case1}
\end{figure}

\begin{figure}
\setcounter{subfigure}{6}
\addtocounter{figure}{-1}
\centering
\mbox{
\subfigure[Initial temperature field $u(\C_3)$]{
\includegraphics[scale=0.33]{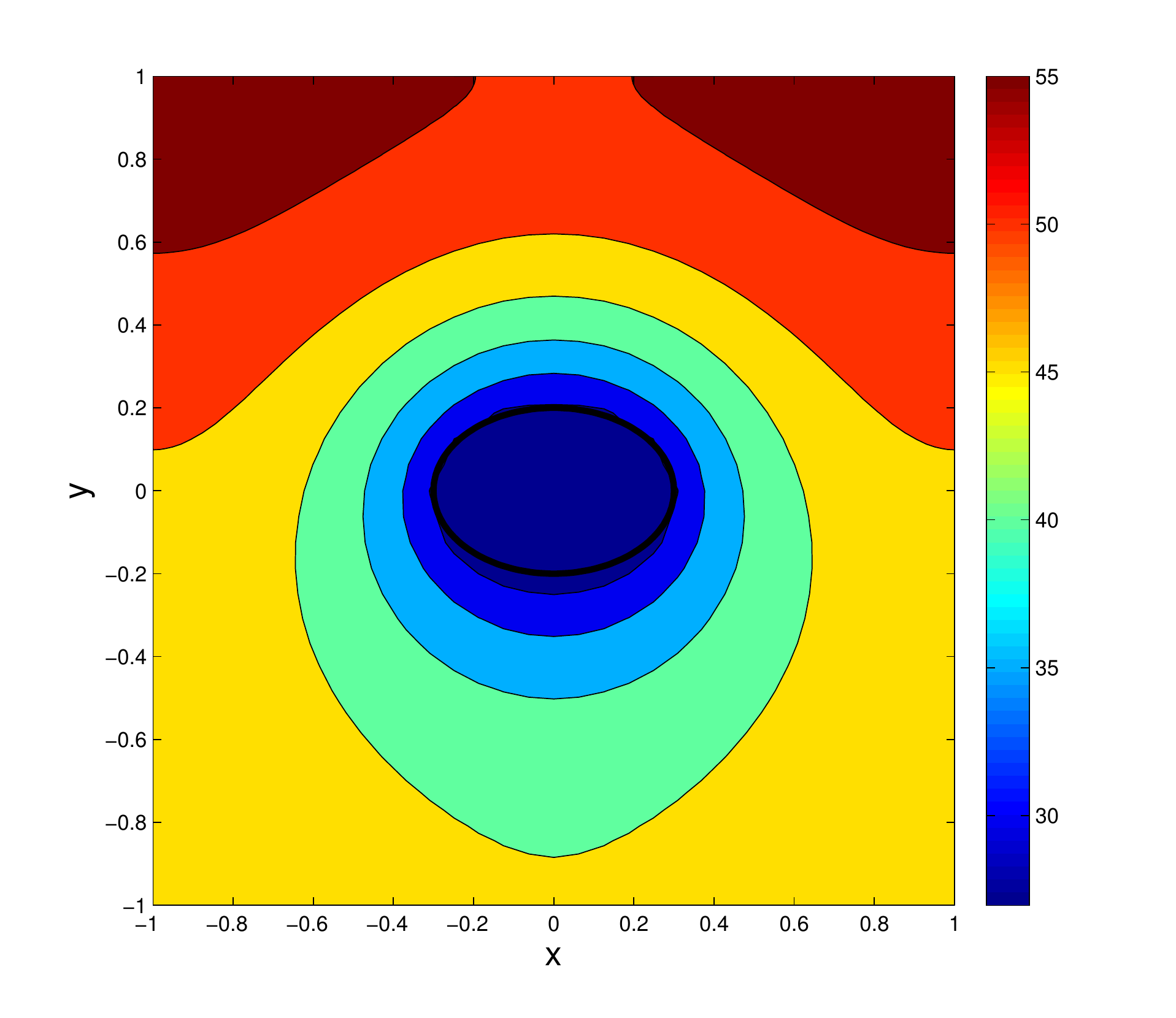}}\qquad
\subfigure[Optimal temperature field $u(\tC)$ in the case with initial contour $\C_3$]{
\includegraphics[scale=0.34]{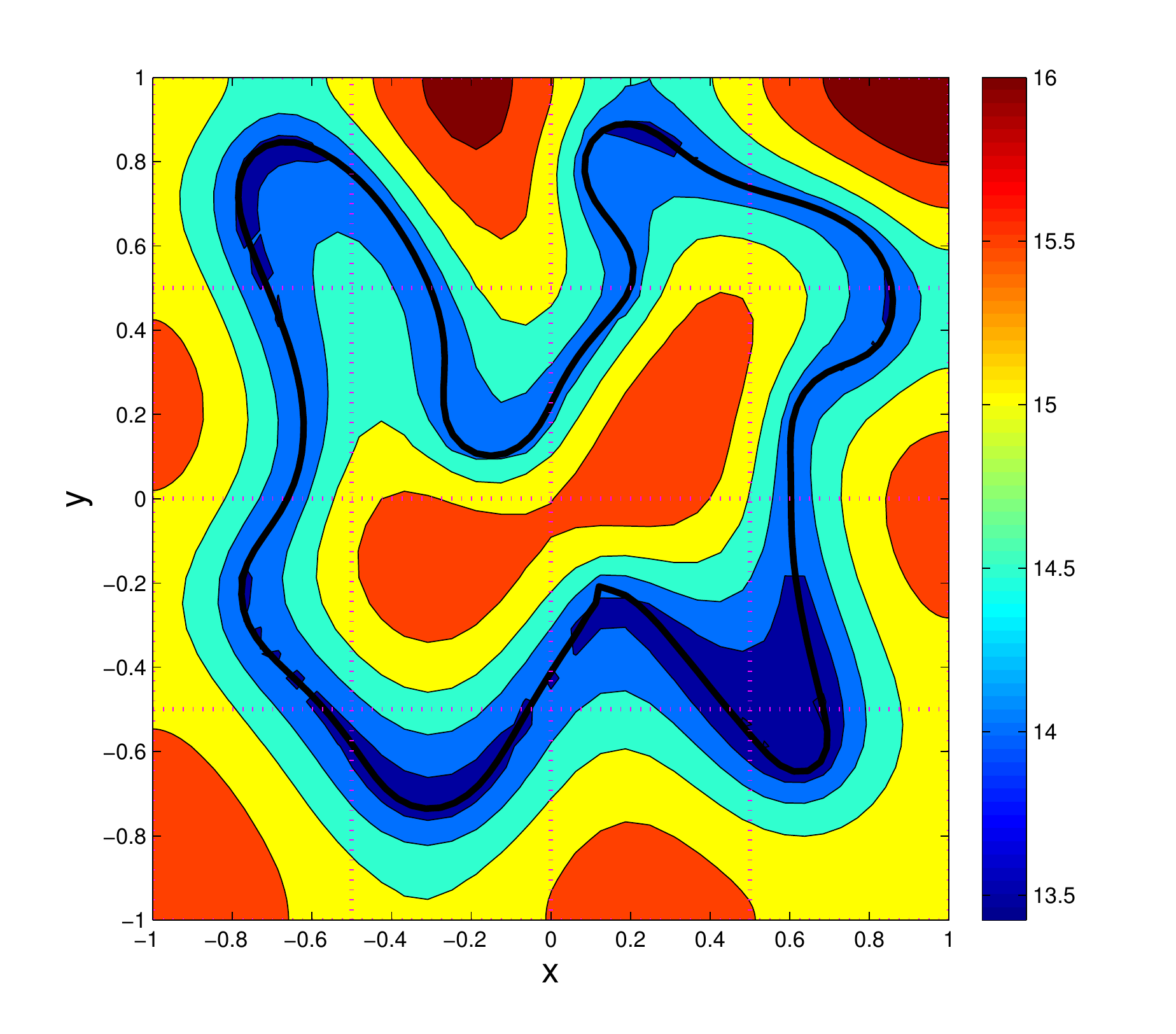}}}
\mbox{
\subfigure[Initial temperature field $u(\C_4)$]{
\includegraphics[scale=0.33]{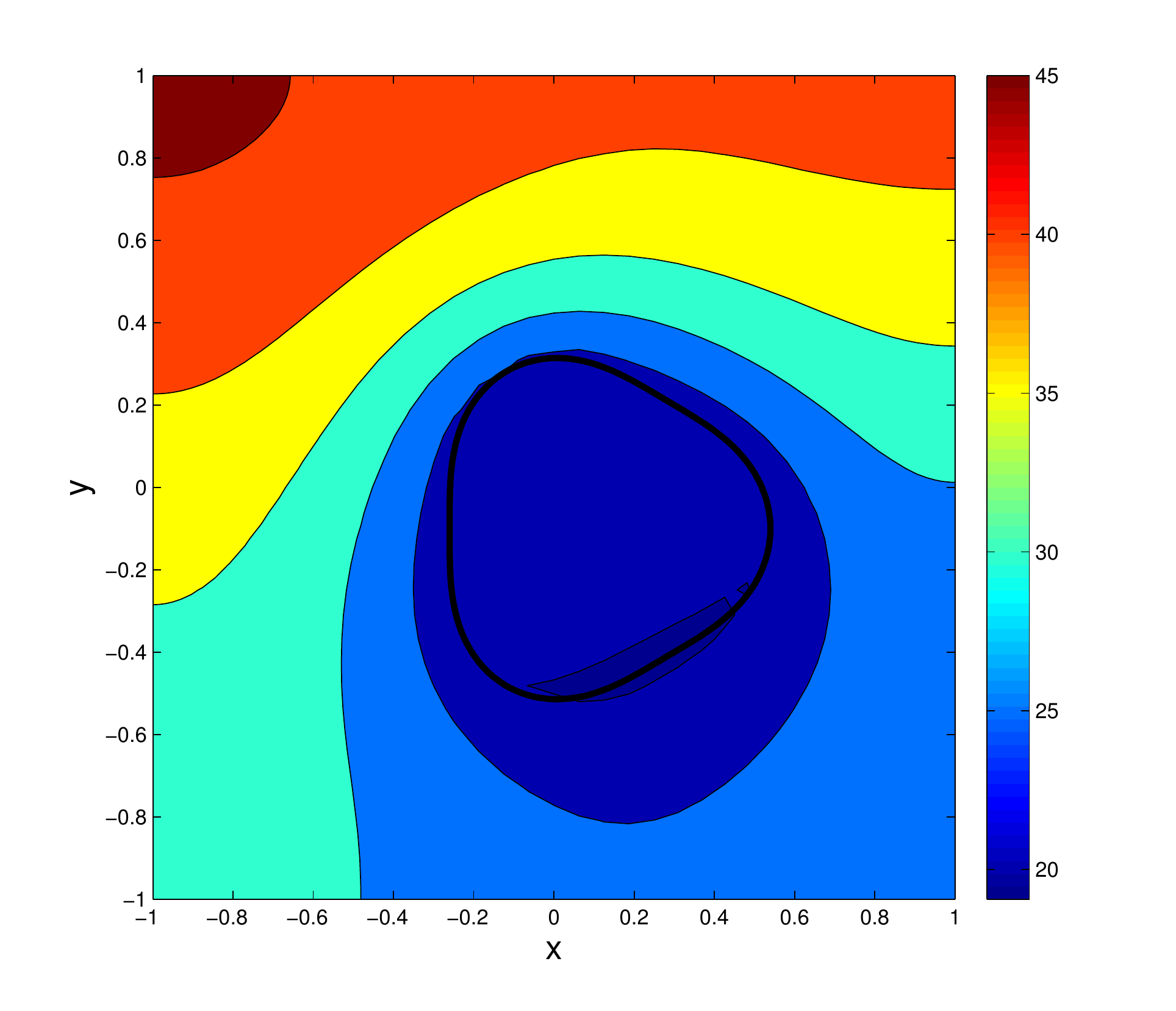}}\qquad
\subfigure[Optimal temperature field $u(\tC)$ in the case with initial contour $\C_4$]{
\includegraphics[scale=0.34]{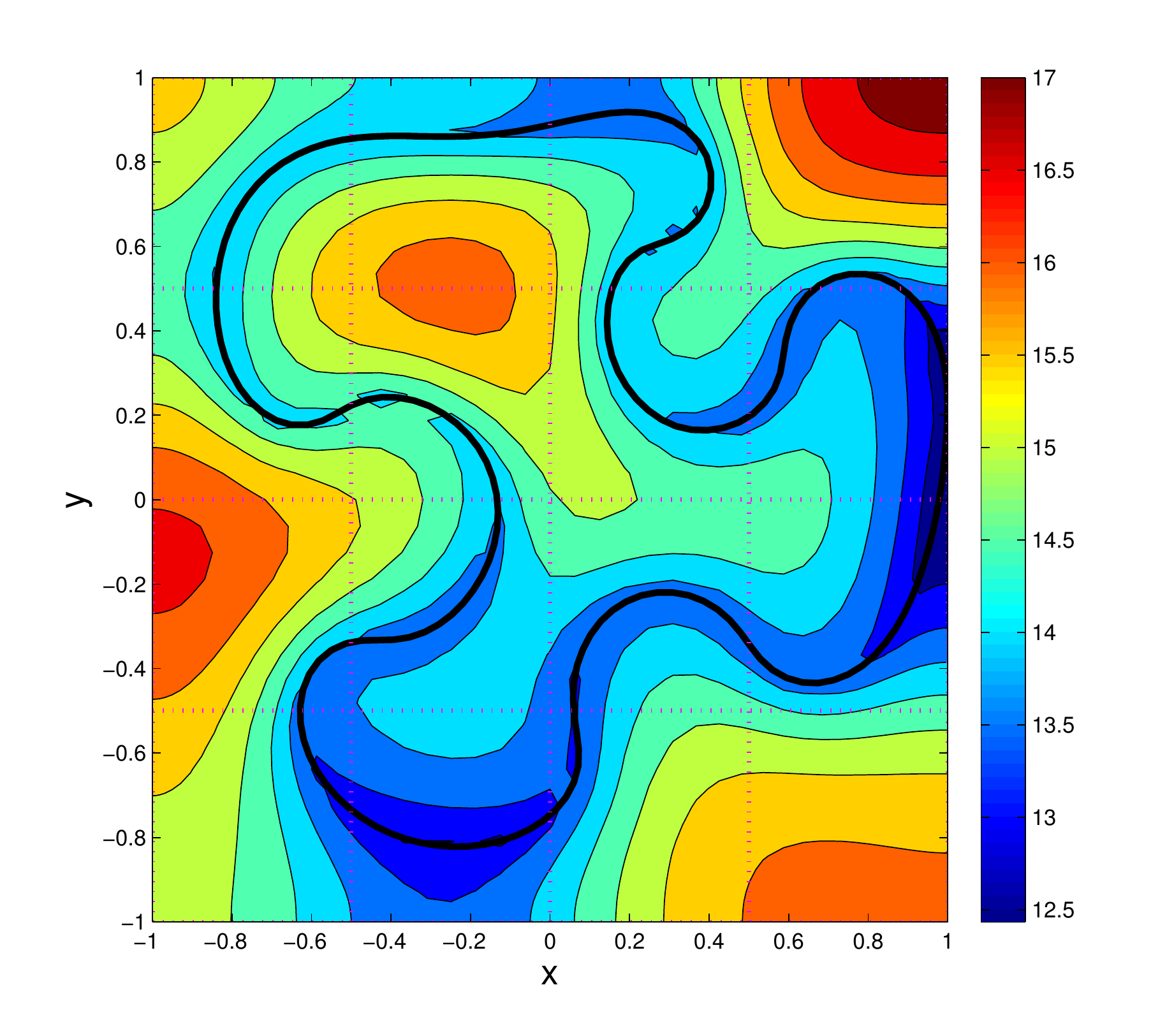}}}
\mbox{
\subfigure[Initial temperature field $u(\C_5)$]{
\includegraphics[scale=0.33]{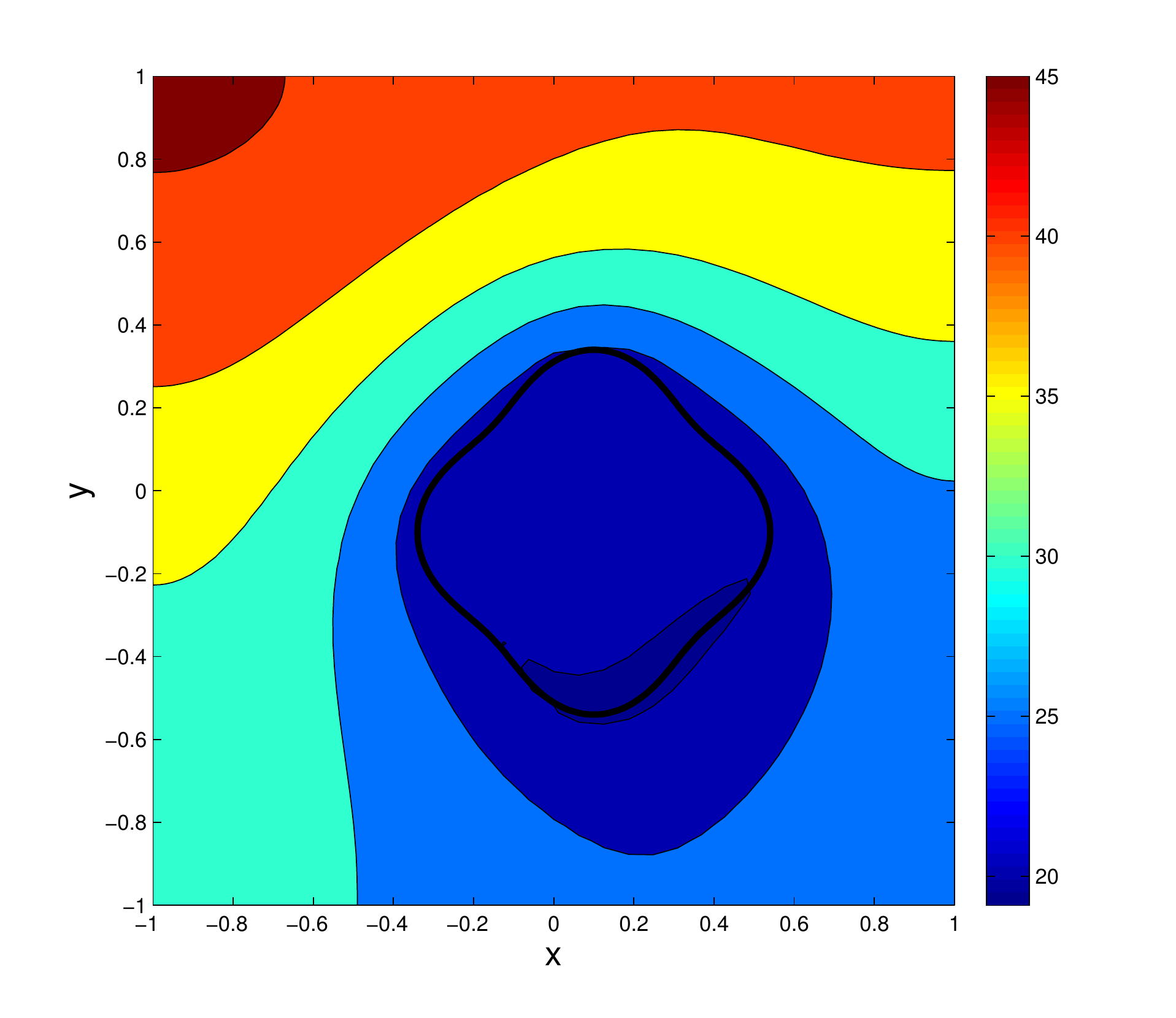}}\qquad
\subfigure[Optimal temperature field $u(\tC)$ in the case with initial contour $\C_5$]{
\includegraphics[scale=0.34]{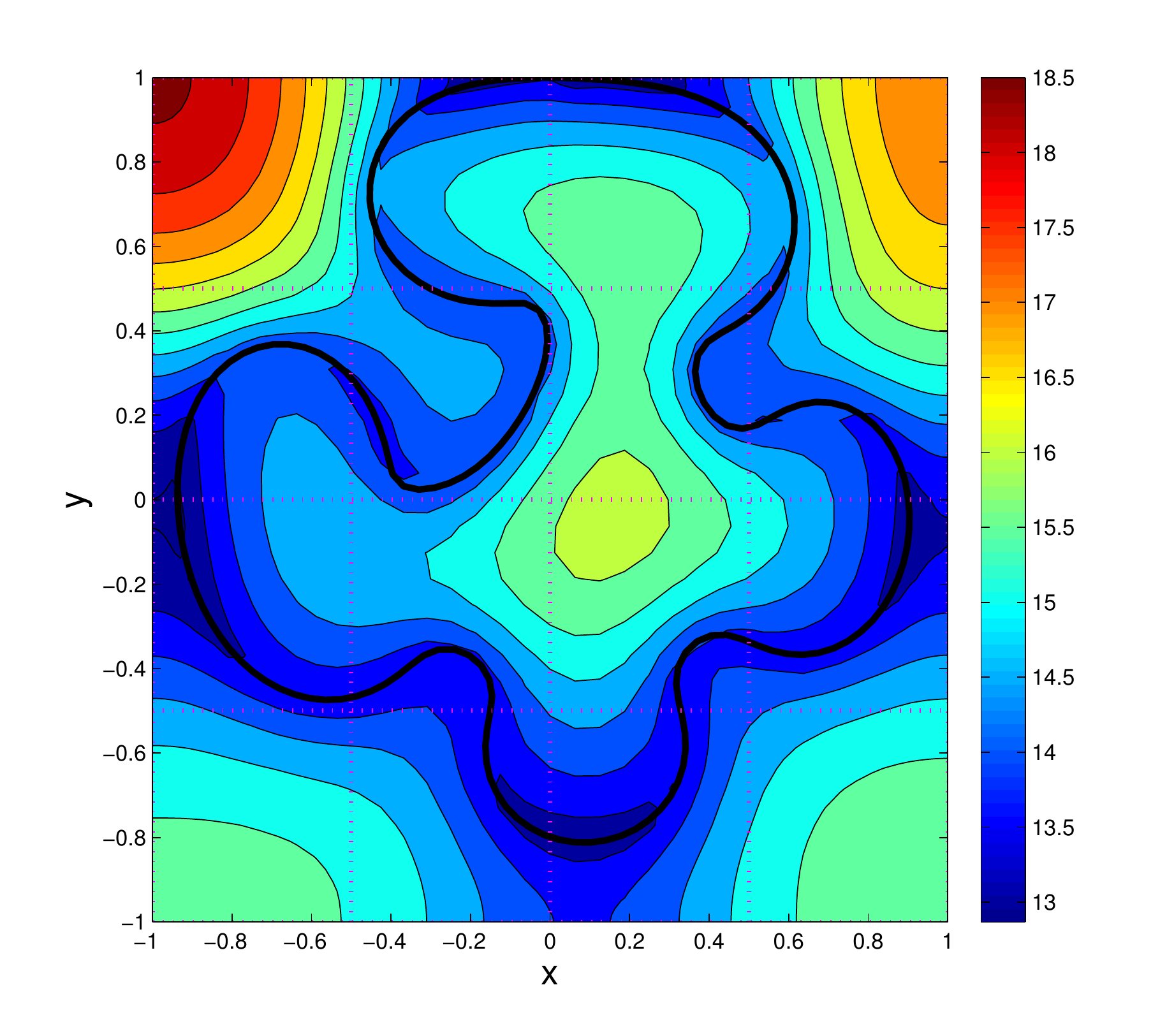}}}
\caption{(continued) Results illustrating solution of the optimization
  problems in CASE \#1, cf.~Table \ref{tab:cases}. The grid shown in
  the Figures in the right column corresponds to the cellular pattern
  of the target field $\ubar$, cf.~Figure \ref{fig:qub}(b).}
\end{figure}

The results characterizing the performance of Algorithm \ref{alg1} in
CASE \#1 are collected in Figure \ref{fig:case1}. First of all, we
note that, depending on the choice of the initial guess $\C^{(0)}$ for
the contour in \eqref{eq:iter}, cf.~Figure \ref{fig:case1}a, the
iterations in fact converge to quite distinct locally optimal shapes,
cf.~Figure \ref{fig:case1}b, providing evidence for the existence of
multiple local minima in the optimization problem, as discussed in
Section \ref{sec:grad}. We also note that the decrease of cost
functional $\J(\C^{(n)})$ with iterations $n$ is quite different in
these different cases, cf.~Figure \ref{fig:case1}c. In Figure
\ref{fig:case1}d we show the intermediate shapes found at the
consecutive iterations of the algorithm starting from the initial
guess $\C_2$ for which the largest decrease was obtained in the cost
functional. We observe that ``simpler'' initial guesses (i.e., a
circle or an ellipse) tend to lead to ``better'' local minimizers.
However, the final temperature distributions $u(\tC)$ obtained from
the different initial guesses all capture features of the cellular
pattern characterizing the target distribution $\ubar$, see Figures
\ref{fig:case1}f,h,j,l vs. Figure \ref{fig:qub}b.

\begin{figure}
\centering
\mbox{
\subfigure[Evolution of cost functional $\J_{\alpha}(\C^{(n)})$ with iterations]{
\includegraphics[scale=0.32]{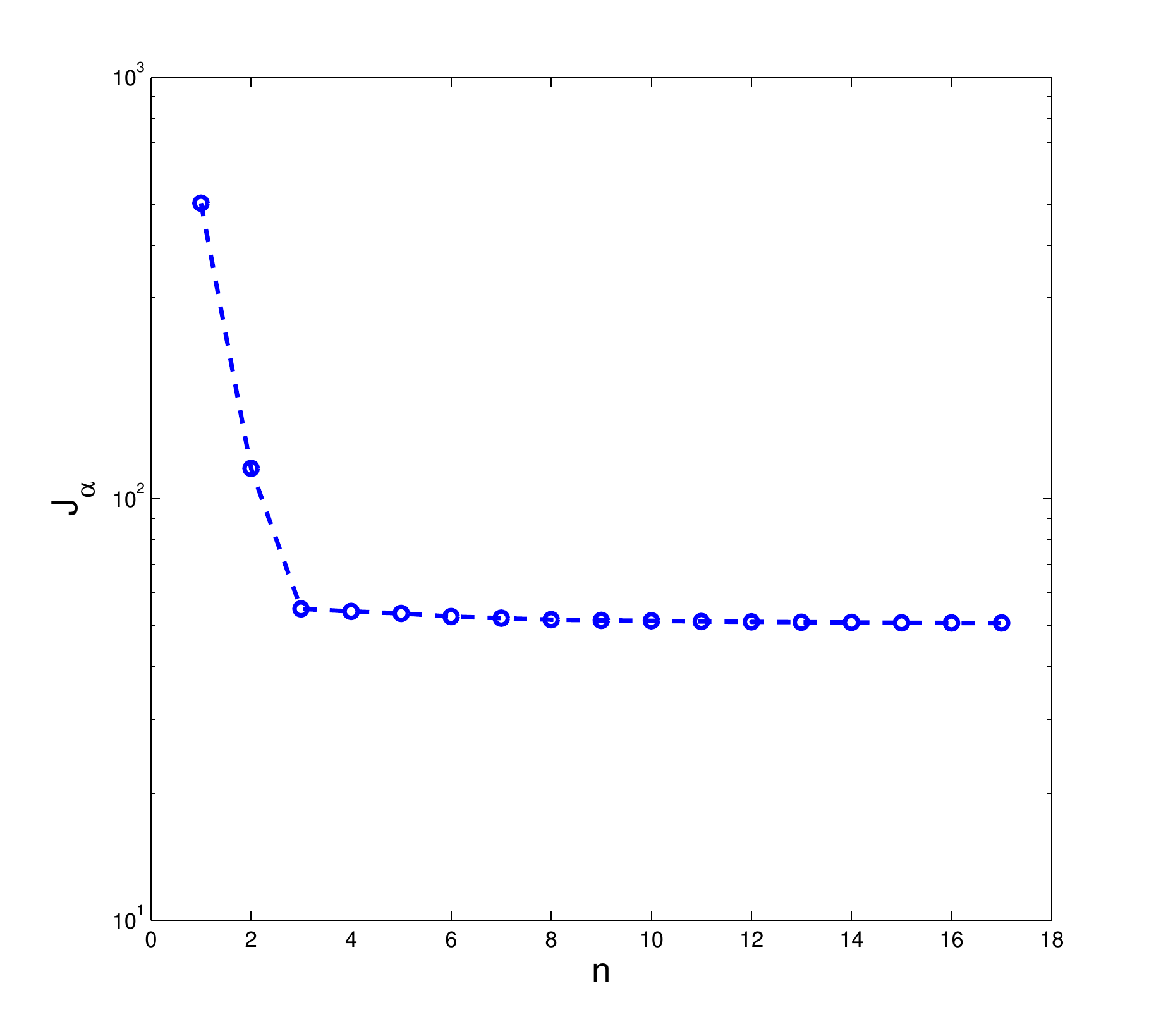}}\qquad
\subfigure[Evolution of contour length $L(\C^{(n)})$ with iterations]{
\includegraphics[scale=0.33]{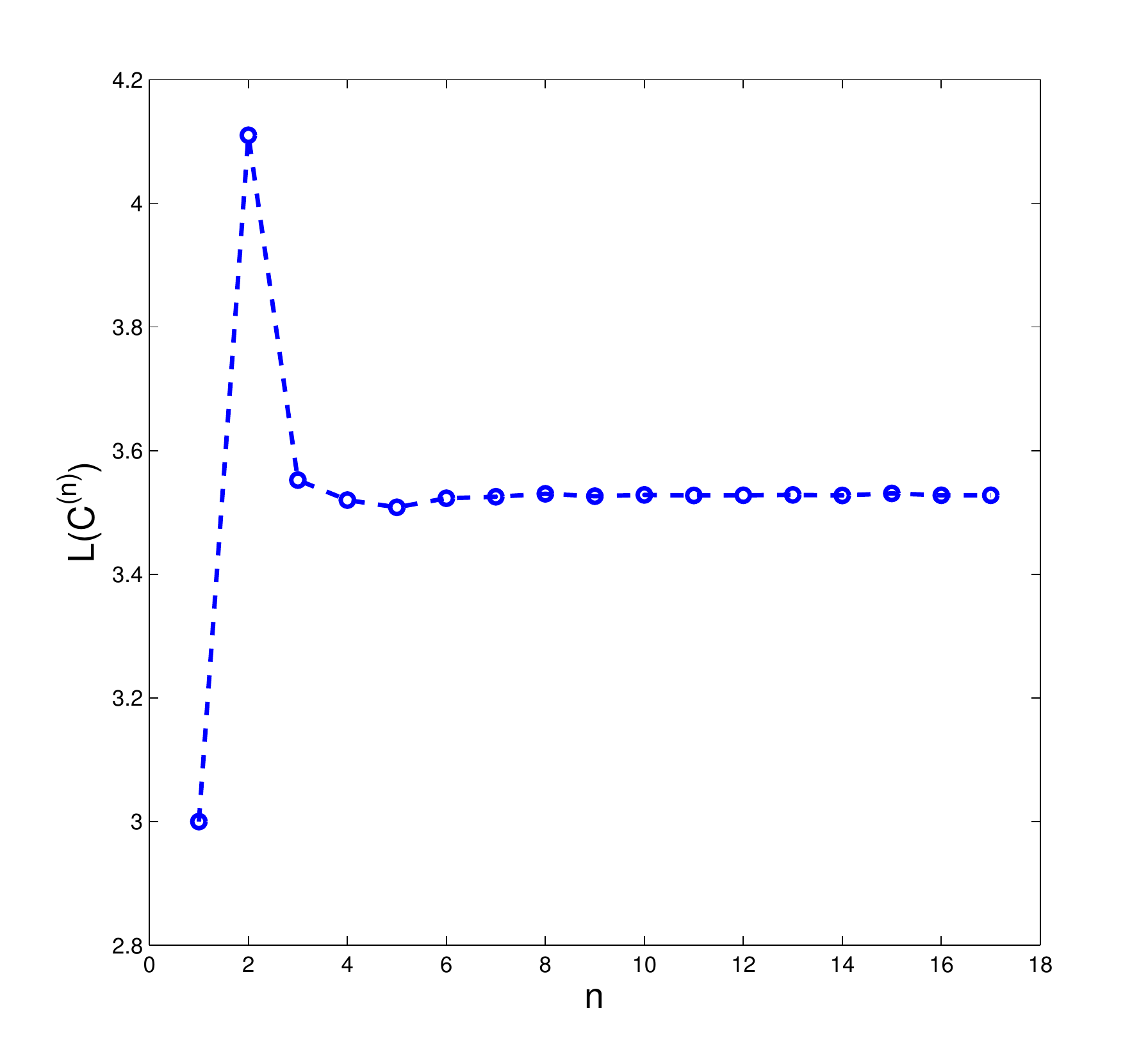}}}
\mbox{
\subfigure[Evolution of contours $\C^{(n)}$ with iterations; the 
initial contour $\C_7$ is marked with thin dashed line and the optimal
shape $\tC$ appears in red]{
\includegraphics[scale=0.32]{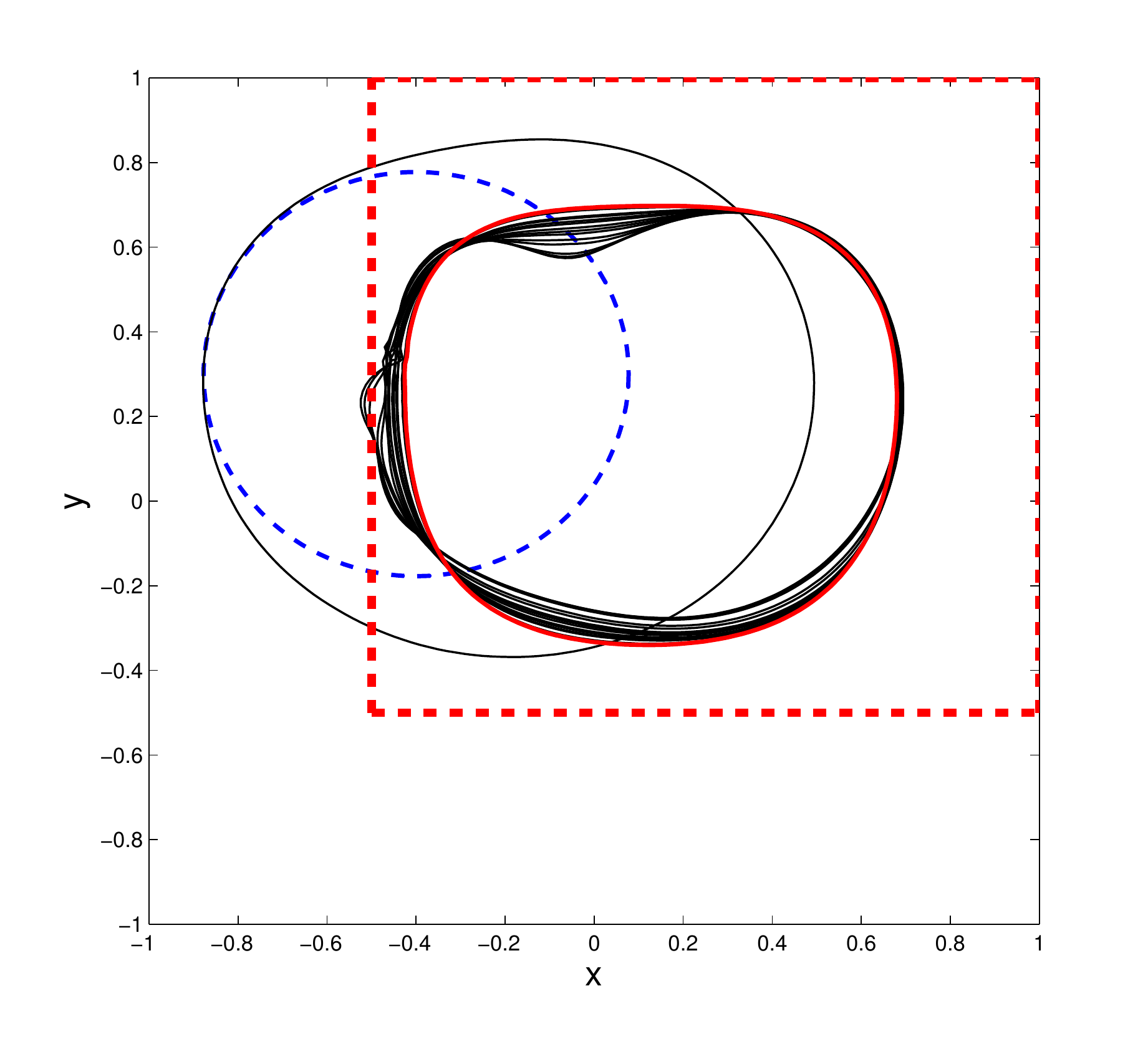}}\qquad
\subfigure[Optimal temperature distribution $u(\tC)$]{\hspace*{-0.15cm}
\includegraphics[scale=0.335]{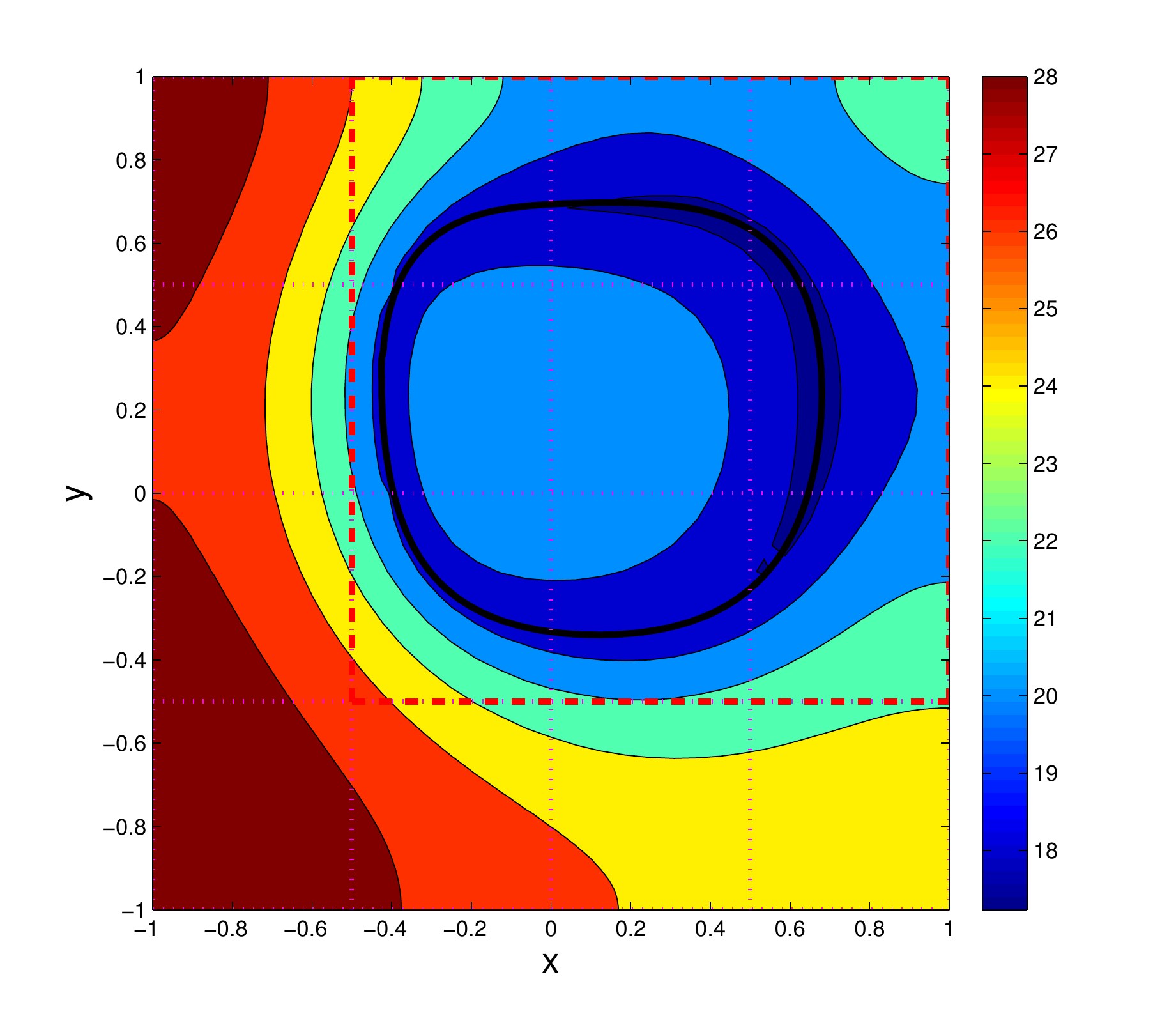}}}
\caption{Results illustrating solution of the optimization problem in
  CASE \#2, cf.~Table \ref{tab:cases}. The rectangles marked with
  thick dashed lines in Figures (c) and (d) indicate the subregion
  $\A$ where the target temperature $\ubar$ is specified.}
\label{fig:case2}
\end{figure}

The data illustrating the performance of Algorithm \ref{alg1}
in CASE \#2 is shown in Figure \ref{fig:case2}. Since in this case we
include the length constraint with a rather large value of the penalty
parameter ($\alpha = 10^2$), the optimal contours are not allowed to
deform much (Figure \ref{fig:case2}c). However, we remark that the
algorithm is able to ``shift'' the contour so that the optimal shape
$\tC$ is enclosed within the target domain $\A$ in which the
temperature field $\ubar$ is defined (Figure \ref{fig:case2}d).

\begin{figure}
\centering
\mbox{
\subfigure[Distribution of heat sources $q$ \cite{ksk09}]{
\includegraphics[scale=0.475]{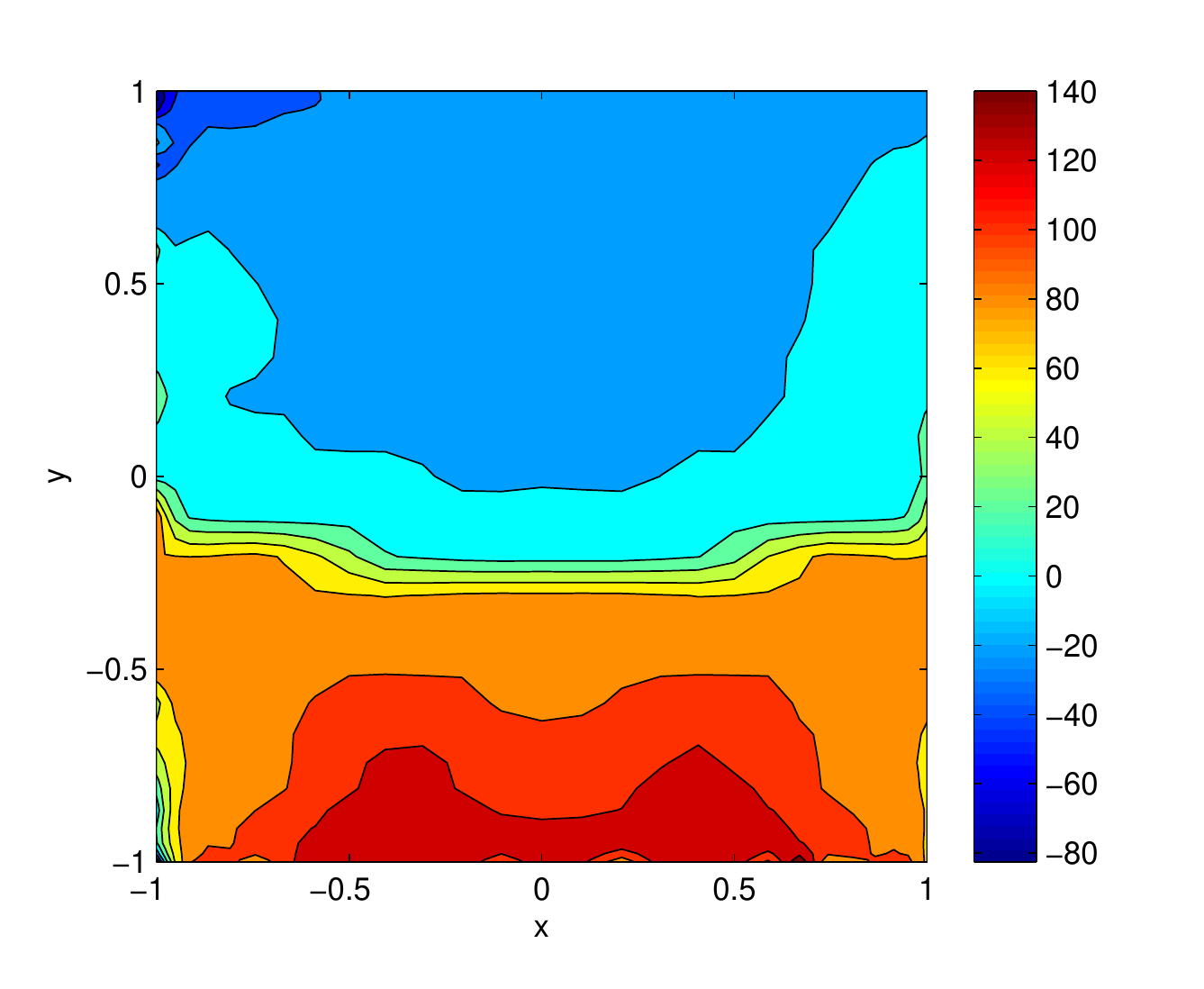}}
\subfigure[Initial temperature field $u(\C_7)$]{
\includegraphics[scale=0.475]{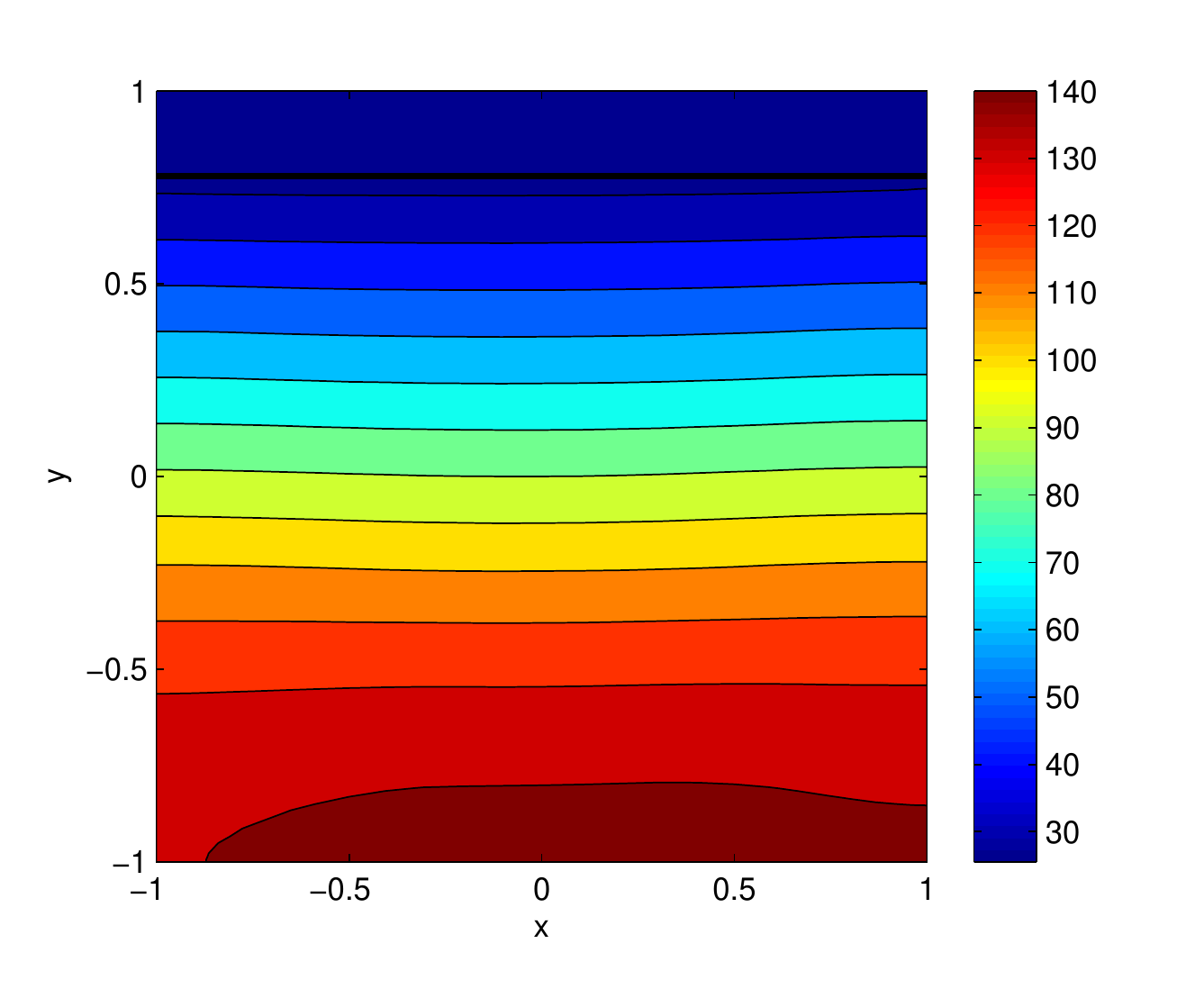}}}
\mbox{
\subfigure[Evolution of cost functional $\J(\C^{(n)})$ with iterations]{
\includegraphics[scale=0.55]{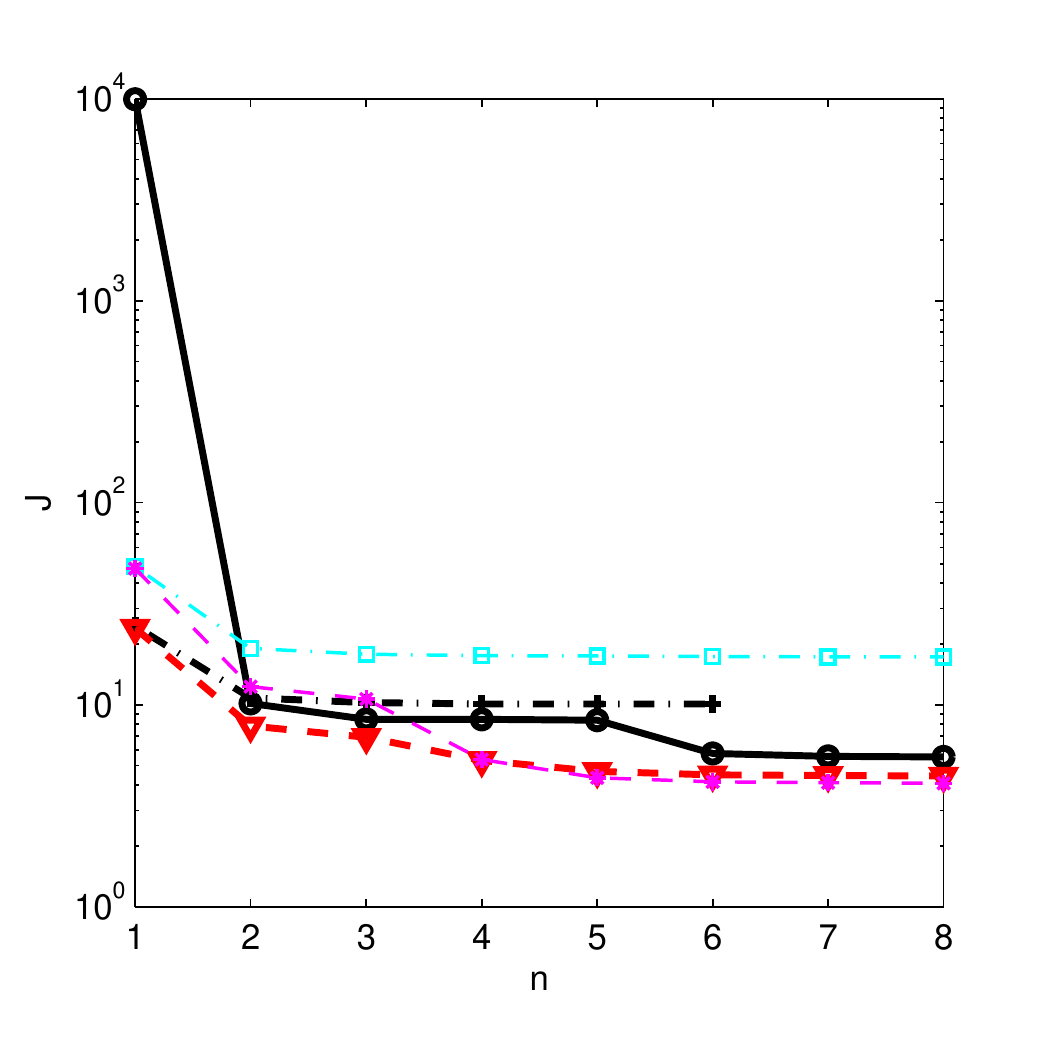}}\quad
\subfigure[Evolution of contour length $L(\C^{(n)})$ with iterations]{
\includegraphics[scale=0.55]{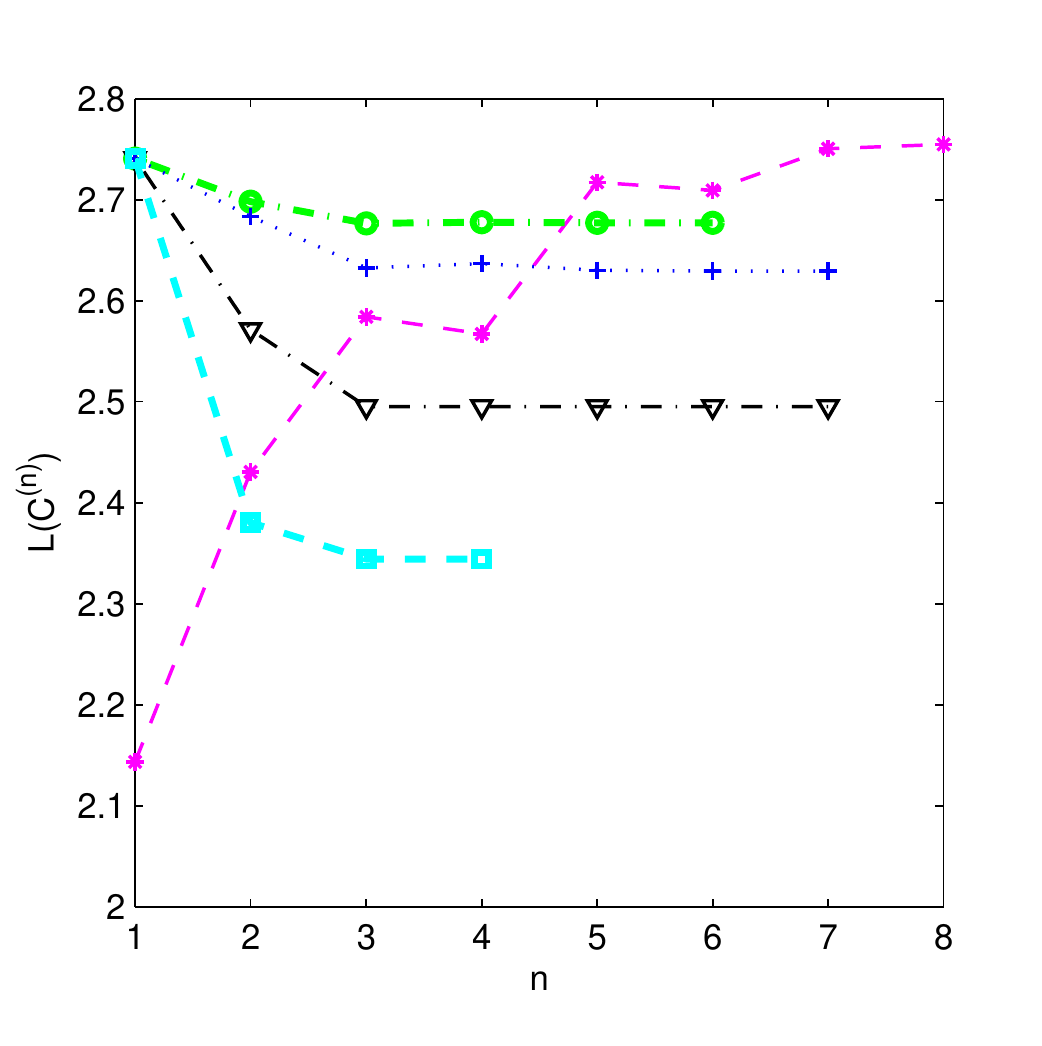}}}
\mbox{
\subfigure[Optimal contours $\tC$ obtained without the length constraint]{
\includegraphics[scale=0.55]{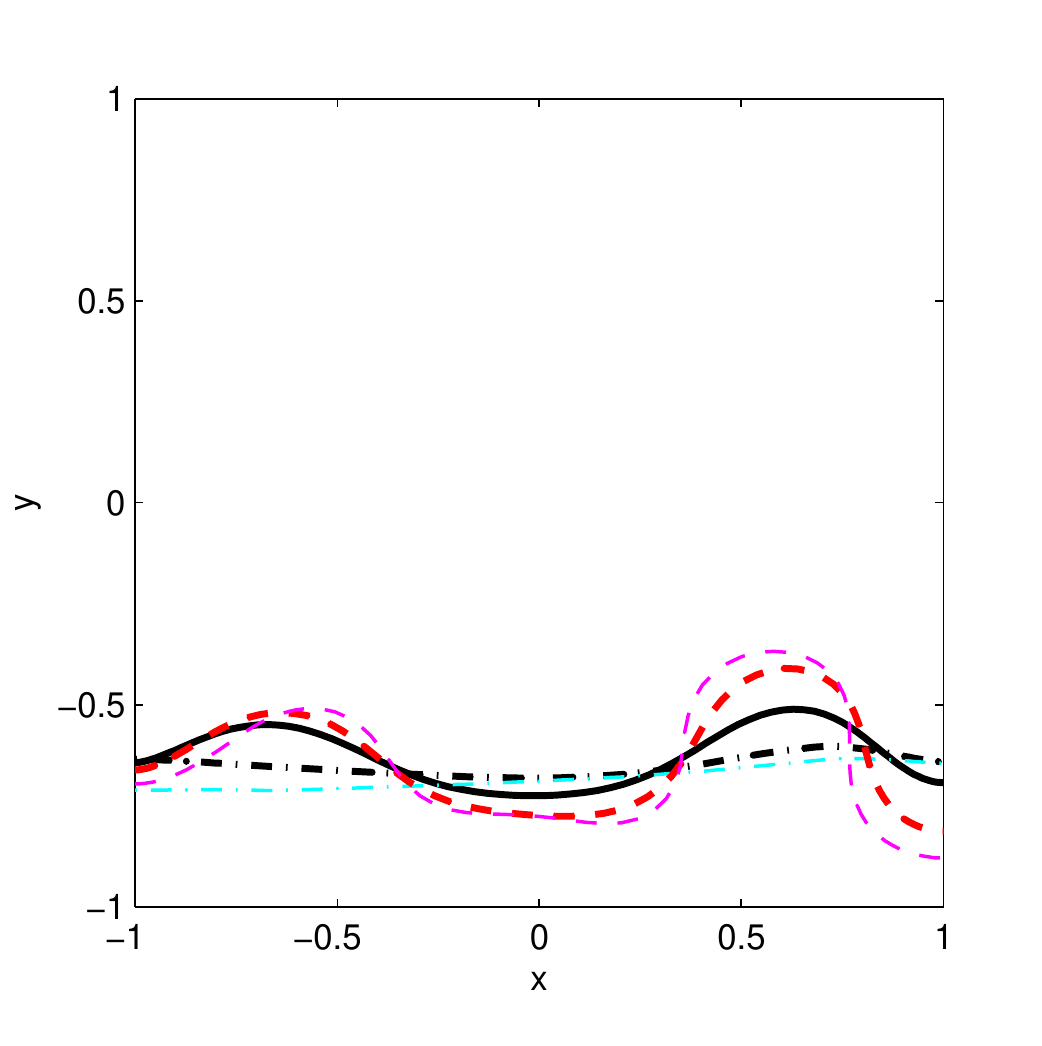}}\quad
\subfigure[Optimal contours $\tC$ obtained with the length constraint]{
\includegraphics[scale=0.55]{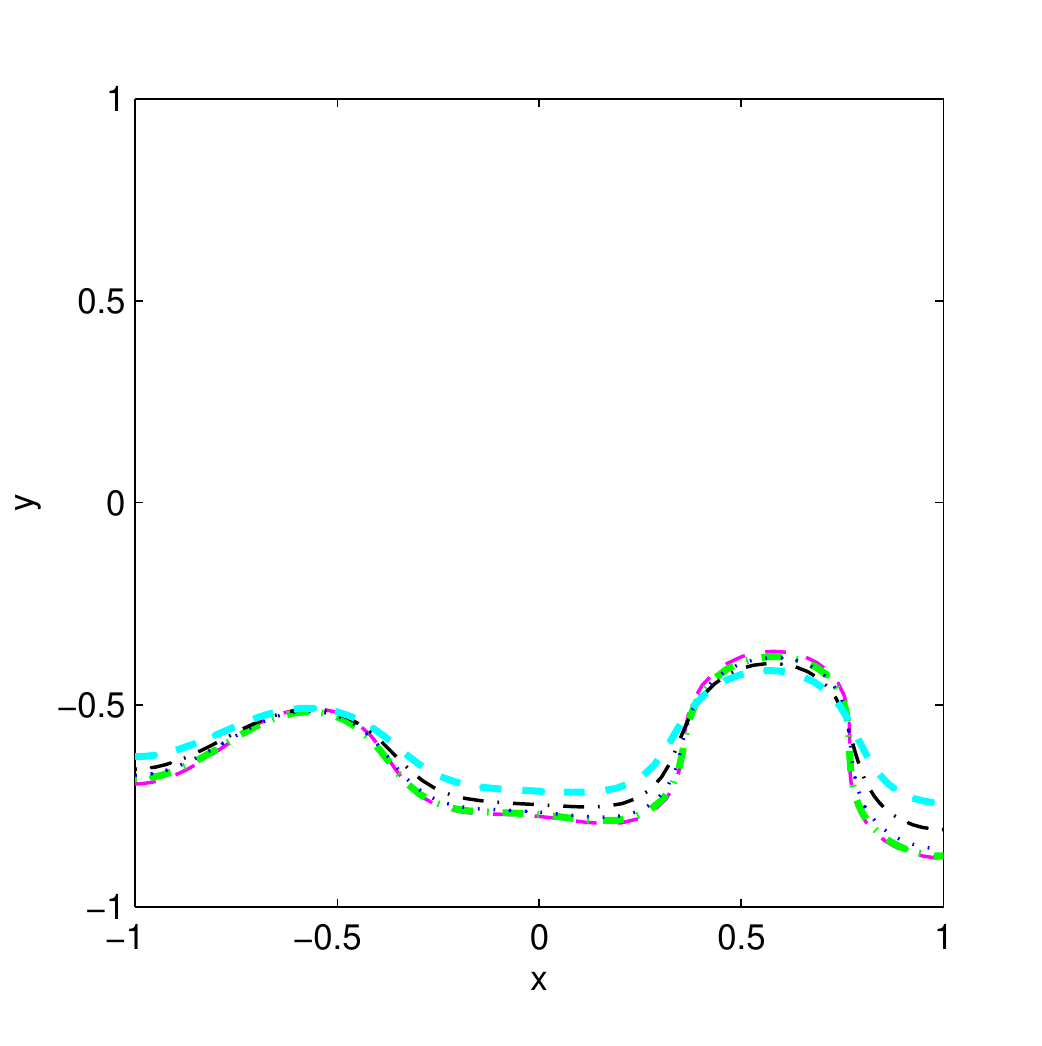}}}
\caption{Results illustrating solution of the optimization problems in
  CASE \#3, cf.~Table \ref{tab:cases}; {in Figures (c) and (e):
    (solid) $u_0 = 10 = \textrm{Const}$, (thick dashed) $T_a=10$ and
    $T_b=16$, (thin dashed) $T_a=10$ and $T_b=19$, (thick dash-dotted)
    $T_a=4$ and $T_b=10$, (thin dash-dotted) $T_a=1$ and $T_b=10$; in
    Figures (d) and (f): (thin dashed) $\alpha=0$, (thick dash-dotted)
    $\alpha=1$, (dotted) $\alpha=10$, (thin dash-dotted)
    $\alpha=10^2$, (thick dashed) $\alpha=10^3$.}}
\label{fig:case3}
\end{figure}

\begin{figure}
\setcounter{subfigure}{6}
\addtocounter{figure}{-1}
\centering
\mbox{
\subfigure[Optimal temperature distribution $u(\tC)$
obtained with $u_0 = 10 = \textrm{Const}$ and $\alpha=0$]{
\includegraphics[scale=0.45]{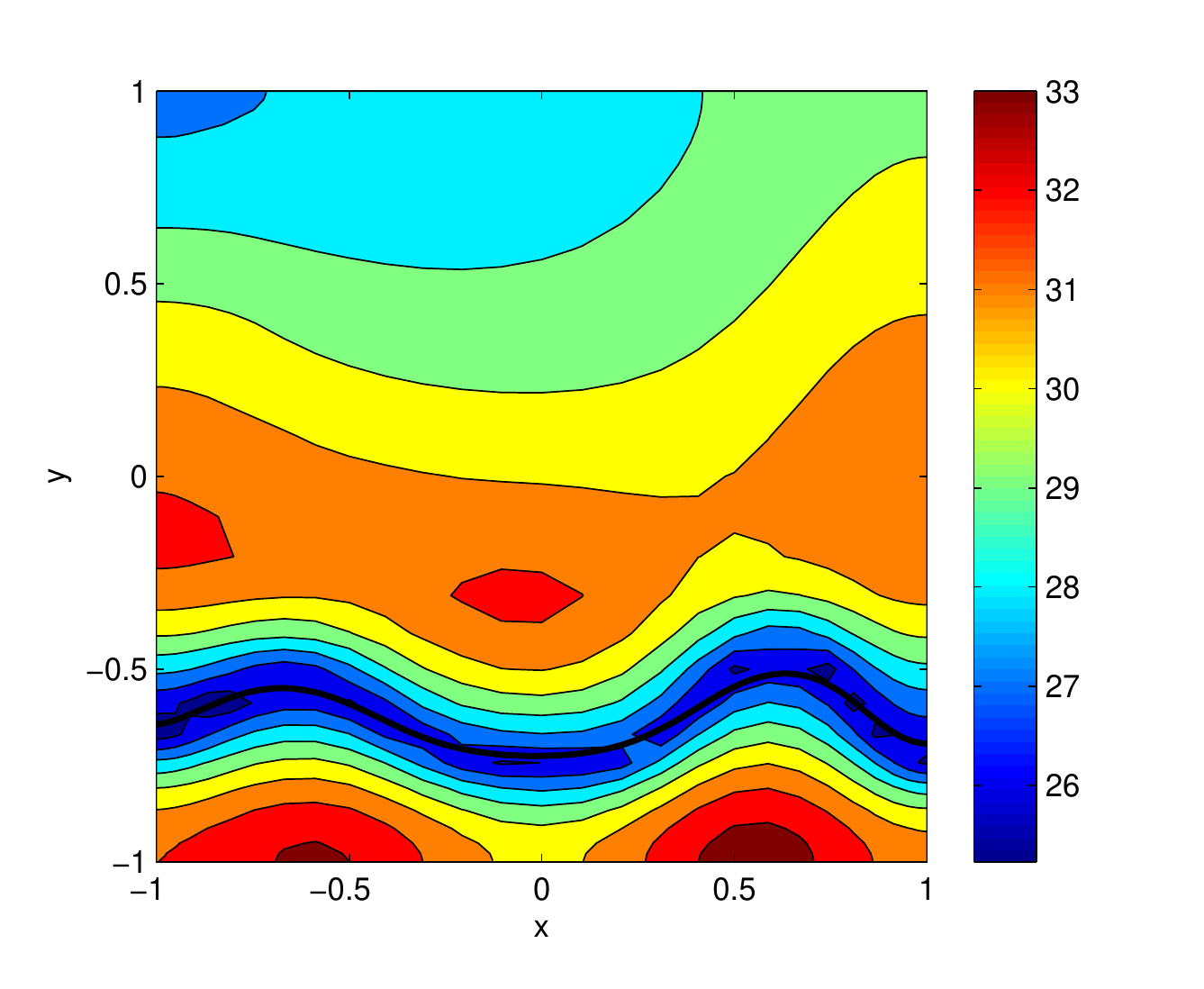}}\quad
\subfigure[Optimal temperature distribution $u(\tC)$
obtained with $T_a=1$, $T_b=10$ and $\alpha=0$]{
\includegraphics[scale=0.45]{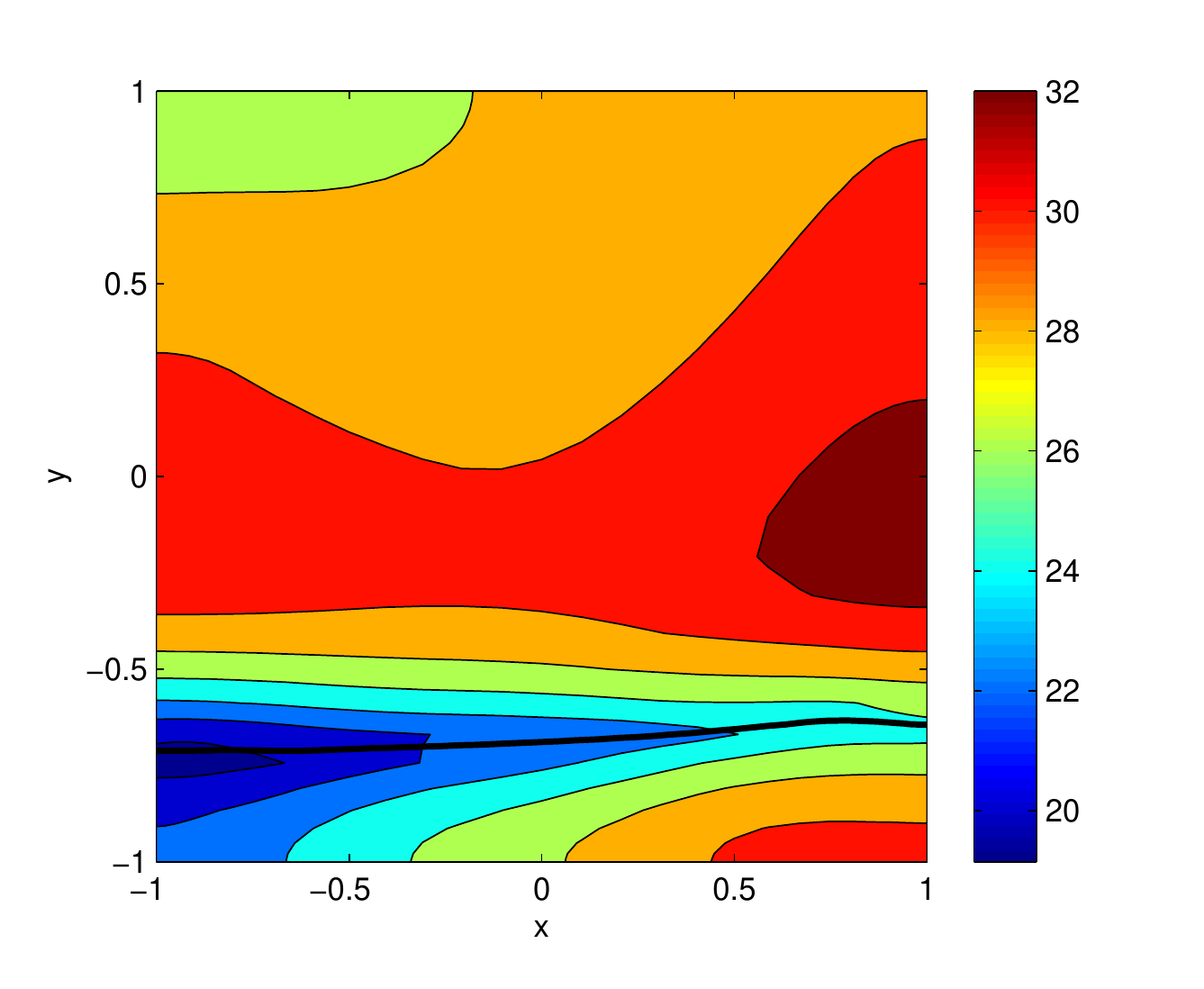}}}
\mbox{
\subfigure[Optimal temperature distribution $u(\tC)$ 
obtained with $T_a=10$, $T_b=19$ and $\alpha=0$]{
\includegraphics[scale=0.45]{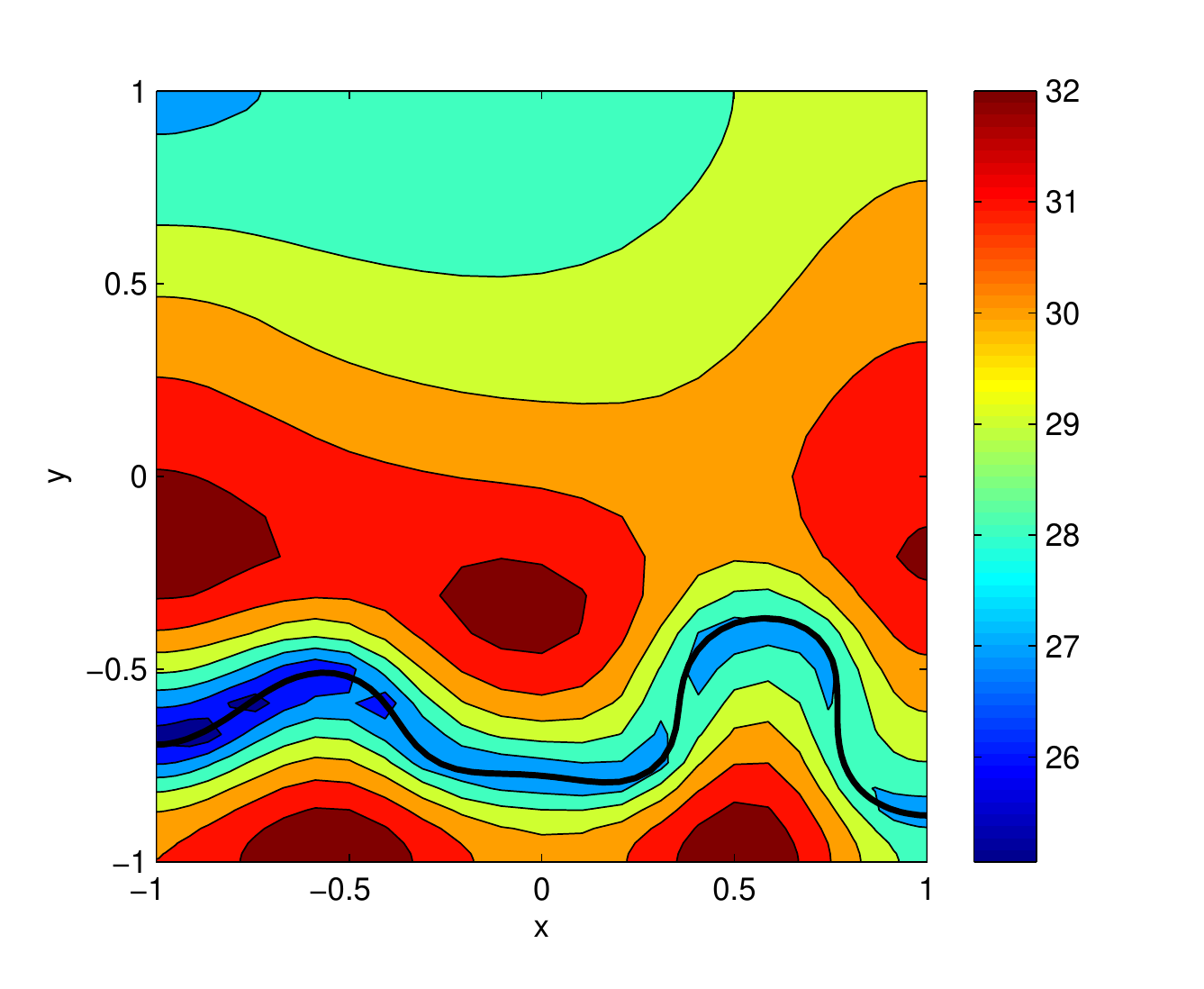}}\quad
\subfigure[Optimal temperature distribution $u(\tC)$ 
obtained with $T_a=10$, $T_b=19$ and $\alpha=1000$]{
\includegraphics[scale=0.45]{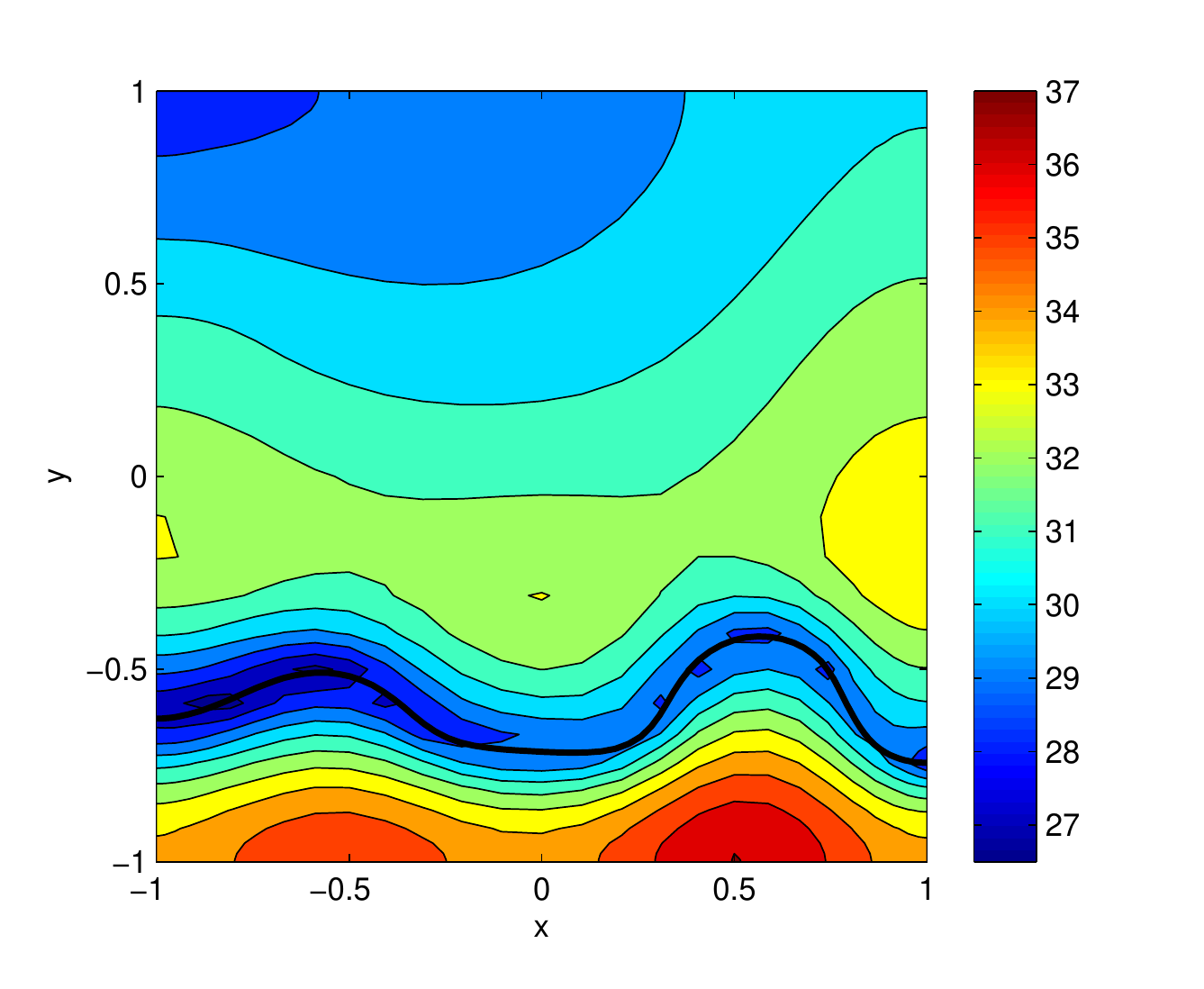}}}
\caption{(continued) Results illustrating solution of the optimization
  problems in CASE \#3 {with different $u_0$,
    cf.~\eqref{eq:u0}.}}
\end{figure}

{The data concerning CASE \#3 is collected in Figure
  \ref{fig:case3}. Using the contour shown in Figure \ref{fig:case3}b
  as the initial guess (cf.~Table \ref{tab:C}), we first solve Problem
  P2 assuming $u_0 = \textrm{Const}$ and with the length constraint
  not enforced ($\alpha=0$). Then, using thus obtained optimal shape
  (marked with the solid line in Figure \ref{fig:case3}e) as the
  initial guess, we solve Problem P2 again, now allowing $u_0$ to vary
  with the arc length $s$. Mimicking changes in the inflow/outflow
  temperature of the coolant liquid, this is achieved by decreasing
  $T_a$ or increasing $T_b$ in \eqref{eq:u0} and corresponds to,
  respectively, dash-dotted and dashed contours in Figure
  \ref{fig:case3}e. In Figure \ref{fig:case3}c we observe that in the
  initial optimization the cost functional drops by over three orders
  of magnitude during less than 10 iterations (in the subsequent
  problems which have better initial guesses this decrease is
  smaller).  Finally, we consider the case with $T_a=10$ and $T_b=19$,
  and solve the optimization problem with the length constraint using
  $L_0 = 2.3$ and increasing values of $\alpha$. The resulting optimal
  contour shapes are shown in Figure \ref{fig:case3}f, whereas Figure
  \ref{fig:case3}d presents the evolution of the contour length
  $L(\C^{(n)})$ with iterations for different values of $\alpha$. As
  expected, we see that for increasing values of $\alpha$ the contour
  length approaches the prescribed value $L_0$ while the contours
  themselves become less deformed. The temperature fields $u(x,y)$
  obtained in the cases with $u_0 = 10 = \textrm{Const}$, $T_a=1$ and
  $T_b=10$, and $T_a=10$ and $T_b=19$, cf.~\eqref{eq:u0}, and without
  the length constraint are shown in Figures \ref{fig:case3}g,h,i.
  This last case with the length constraint and $\alpha=1000$ is shown
  in Figure \ref{fig:case3}j. We see that, as compared to the
  temperature corresponding to the initial guess for the contour
  (Figure \ref{fig:case3}b), the optimal distributions in Figures
  \ref{fig:case3}g,h,i,j have the temperature ranges much closer to
  target field \eqref{eq:ubar}. We also observe that, with the
  exception of the case in which the inflow temperature $T_a$ is quite
  low (Figure \ref{fig:case3}h), the optimal contour shapes tend to
  weave around the two hot spots in the heat source distribution
  (Figure \ref{fig:case3}a) in a complicated manner.}

\section{Conclusions and Future Work}
\label{sec:final}

In this investigation we have addressed the problem of shape
optimization for a {system of elliptic PDEs subject to mixed
  interface boundary conditions which models} the steady-state heat
transfer in 2D. {Our continuous optimization formulation relies
  on Sobolev shape gradients obtained using the shape-differential
  calculus and explicit interface tracking employed to represent the
  optimized contour, all of which are rather well known approaches.
  The key novel contribution of this work is the method we introduced
  to numerically evaluate the shape gradients. By splitting the
  adjoint system into two coupled subproblems we could achieve optimal
  accuracy for each of the subproblems. In particular, the proposed
  boundary-integral formulation exploits the analytical (potential)
  structure of the problem and, as demonstrated by the exhaustive
  validation tests presented in Section \ref{sec:validate}, offers
  high numerical accuracy without the need to construct a
  boundary-fitted mesh at every iteration, as required in other
  approaches based on explicit interface tracking \cite{ws10}. As a
  result, the proposed method is quite efficient from the
  computational point of view and, as shown in Section
  \ref{sec:optim}, can deal with fairly complicated contour shapes in
  an easy and straightforward manner. While boundary-integral
  formulations have been used to study the shape sensitivities of
  elliptic PDEs (e.g., \cite{nr00,hv02,rg07}), the method we
  introduced is designed for higher accuracy than previous
  approaches.}

{As compared to the ``discretize-then-differentiate'' approaches,
  they key advantage of the continuous formulation used here is that,
  as discussed at the end of Section \ref{sec:implement}, it offers
  the freedom to remesh the points discretizing the contour which is
  crucial to achieving spectral accuracy in the solution of the
  boundary integral equation. On the contrary, in the discrete setting
  such remeshing would actually require one to set up a new
  optimization problem (corresponding to the new set of the discrete
  control variables). Moreover, it is also not evident how the
  analytic treatment of singularities described in Section
  \ref{sec:implement} could be employed in the discrete setting. In
  regard to the level-set-based interface capturing methods, the
  present approach arguably offers more flexibility in the
  high-accuracy treatment of the complex interface boundary
  conditions.}

Optimizations performed on three test problems led to rather
nonintuitive optimal shapes of the contour which {differed
  significantly from the initial guesses provided based on the
  ``engineering intuition'', reflecting} the geometric nonlinearity
{and nonlocality} of the governing system. Smoothness of the
contours was enforced by defining the cost functional gradients in a
suitable Sobolev space.  {This, combined with the interpolation
  technique described at the end of Section \ref{sec:implement},
  allowed for an accurate representation of even strongly deformed
  contours using a rather modest number of points ($M=100$ in CASES
  \#1, \#2, and \#3, cf.  Table \ref{tab:cases}).}  Evidence was also
shown for the presence of multiple local minima. {In Problem P1,}
when the length constraint was not imposed, the optimal temperature
distributions were found to capture {the} main features of the
target temperature field $\ubar$. On the other hand, the presence of
the length constraint restricted the ability of the algorithm to
deform the contour, although it was still capable of ``shifting'' the
contour to a different part of the domain $\Omega$ without significant
shape changes. {It ought to be added that extension of the
  proposed approach to three-dimensional (3D) configurations is
  conceptually straightforward. Aside from the need to work with the
  3D fundamental solutions in expressions resulting from ansatz
  \eqref{eq:uh}, some technical complications may arise from the fact
  that the boundary integral equations will be formulated on 2D
  surfaces, rather than on 1D contours which, in particular, may make
  achieving high numerical accuracy more difficult.}

{The formulation developed in this study leads to the following
  open problems of a more fundamental character. Our adjoint system
  \eqref{eq:adj} was derived in the PDE setting \cite{p11} and only
  then recast in terms of the boundary-integral formulation for the
  purpose of the numerical solution. On the other hand, one could
  begin with the boundary-integral formulation of governing system
  \eqref{eq:intr_d} which, after shape differentiation, would give
  rise to an integral expression with more singular, possibly
  hypersingular, kernels.  Assessing the relative advantages and
  disadvantages of such an alternative approach is an interesting open
  question and is left to the future research (we mention that
  sensitivity calculations based on hypersingular integral equations
  have already been discussed in \cite{rg07}). In addition,} our
future work will {also} involve generalizations of the proposed
approach to mathematical models of the battery system more complex
than \eqref{eq:intr_d} {and} accounting for some effects of the
{actual} flow of the coolant fluid {in channels of finite
  thickness} (see Figure \ref{fig:battery_pack}b).  We also intend to
explore optimization of the topology of the {contours
  \cite{sz99}}.

\section*{Acknowledgements}

{The two anonymous referees are acknowledged for providing many
  constructive comments and a number of important references.}  The
authors are {also} grateful to the National Centre of Excellence
AUTO21 (Canada) for generous funding provided for this research
(through grant ED401-EHE ``Multidisciplinary Optimization of Hybrid
and Electric Vehicle Batteries''). The authors also acknowledge many
helpful discussions with the General Motors of Canada R\&{D} Team in
Oshawa, Ontario.

\end{document}